\documentclass[12pt]{amsart}
\usepackage{amssymb}
\usepackage{bm}
\usepackage{booktabs}
\usepackage{longtable}
\usepackage{caption}
\usepackage{hyphenat}
\usepackage{mathtools}
\usepackage{enumitem}
\usepackage{xcolor}
\usepackage{graphicx}
\usepackage[lmargin=1in,rmargin=1in,tmargin=1in,bmargin=1in]{geometry}
\usepackage[color=blue!20!white,textsize=tiny,disable]{todonotes}

\definecolor{darkblue}{rgb}{0,0,0.4}
\definecolor{mediumblue}{rgb}{0,0,0.6}
\definecolor{darkred}{rgb}{0.7,0,0.0}
\usepackage[colorlinks=true, citecolor=darkblue, filecolor=darkblue, linkcolor=darkblue,urlcolor=darkblue]{hyperref} 

\allowdisplaybreaks

\usepackage{tikz}
\usetikzlibrary{matrix,arrows,cd}
\tikzcdset{scale cd/.style={every label/.append style={scale=#1},
    cells={nodes={scale=#1}}}}
\definecolor{darkgreen}{rgb}{.35,.51,.3} 

\graphicspath{{draws/}{}}

\newcommand{\RR}{\mathbb R}

\newcommand{\ZZ}{\mathbb Z}
\newcommand{\QQ}{\mathbb Q}

\newcommand{\FF}{\mathbb F}

\newcommand{\Image}{\mathrm{Im}}

\newcommand{\lra}{\longrightarrow}

\newcommand{\co}{\nobreak\mskip2mu\mathpunct{}\nonscript
  \mkern-\thinmuskip{:}\penalty300\mskip6muplus1mu\relax}

\newcommand{\OneHalf}{{\textstyle\frac{1}{2}}}

\newcommand{\bdy}{\partial}

\DeclareMathOperator{\Sym}{Sym}

\DeclareMathOperator{\rank}{rank}

\DeclareMathOperator{\gr}{gr}

\renewcommand{\emptyset}{\varnothing}

\theoremstyle{plain}

\numberwithin{equation}{section}
\newtheorem{theorem}[equation]{Theorem}

\newtheorem{proposition}[equation]{Proposition}
\newtheorem{lemma}[equation]{Lemma}
\newtheorem{corollary}[equation]{Corollary}

\newtheorem{definition}[equation]{Definition}

\theoremstyle{definition}

\newtheorem{question}[equation]{Question}
\newtheorem{problem}[equation]{Problem}

\theoremstyle{remark}
\newtheorem{example}[equation]{Example}
\newtheorem{remark}[equation]{Remark}

\newcommand{\Kh}{\mathit{Kh}}
\newcommand{\KhCx}{\mathit{CKh}}

\newcommand{\KhCxBN}{\KhCx_h} 

\newcommand{\rKh}{\widetilde{\Kh}}
\newcommand{\HF}{\mathit{HF}}

\newcommand{\CFK}{\mathit{CFK}}

\newcommand\HH{\mathit{HH}}

\newcommand\Hochschild\HH

\newcommand\Id{\mathbb{I}}

\renewcommand{\th}{^\text{th}}

\newcommand{\ol}[1]{\overline{#1}{}}
\newcommand{\wt}[1]{\widetilde{#1}{}}

\newcommand{\onto}{\twoheadrightarrow}

\newcommand{\Sq}{\mathrm{Sq}}

\newcommand{\LEO}{\mathit{LEO}}
\newcommand{\LEE}{\mathit{LEE}}

\newcommand{\rCLEO}{\wt{\mathcal{C}}_\LEO}

\newcommand{\LewarkGp}{\mathcal{G}}
\newcommand{\rsa}{\rightsquigarrow}

\DeclareMathOperator{\Span}{span}

\begin{document}
\title{On Khovanov local equivalence groups}

\author{Robert Lipshitz}
 \address{Department of Mathematics, University of Oregon\\
   Eugene, OR 97403, USA}
\thanks{RL was supported by NSF grant DMS-2505715 and a Simons Foundation travel grant.}
\email{\href{mailto:lipshitz@uoregon.edu}{lipshitz@uoregon.edu}}

\date{\today}

\begin{abstract}
	Inspired by work in Seiberg--Witten and Heegaard Floer homology, Dunfield, Sch\"utz, and the author constructed an abelian \emph{local equivalence group} built from pairs of chain complexes, and a homomorphism from the smooth concordance group to this group that uses the odd Khovanov complex and the Bar-Natan deformation of the even Khovanov complex. Lewark independently gave a related local equivalence group and concordance invariant, using only the Bar-Natan deformation. In this paper, we construct $(\ZZ/2\ZZ)^\infty$ subgroups of Lewark's local equivalence group and the kernel of the map from our local equivalence group to Lewark's, as well as $\ZZ^\infty$ direct summands of Lewark's group and this kernel, invisible to the $s$-invariants. The construction gives at least six $\ZZ$-valued concordance homomorphisms linearly independent from $s$, $\tau$, and Dai--Hom--Stoffregen--Truong's invariants.
\end{abstract}

\leavevmode
\vspace{-0.25cm}

\maketitle

\vspace{-0.25cm}

\tableofcontents

\section{Introduction}

Since Rasmussen used it to give another proof of the Milnor Conjecture~\cite{Rasmussen10:s}, obstructing sliceness and concordance has been a headline application of Khovanov homology. It took a few years before the community noticed that there were concordance invariants beyond Rasmussen's $s$, but the last fifteen years have seen an explosion of them (e.g.,~\cite{LS14:refine-s,SSS20:odd-htpy,Schutz:integral-s}, among others), some of which have had attractive applications (e.g.,~\cite{FLL25:squeezed}).
In particular, inspired by work in Seiberg--Witten and Heegaard Floer theory~\cite{Stoffregen20:Seifert,DHST:homol-cob,DHST21:more-conc} and the surprising effectiveness~\cite{DunfieldGong} of the invariant $s^{\Sq^1}_o$ (from~\cite{LS14:refine-s,SSS20:odd-htpy}), Dunfield, Sch\"utz, and the author constructed a local equivalence group $\rCLEO$, as well as a quotient group $\rCLEO^o$ (and several other variants)~\cite{slocal}, and a homomorphism $\LEO$ from the concordance group $\mathcal{C}$ to $\rCLEO$. Independently, Lewark constructed a local equivalence group $\LewarkGp$, which is a direct summand of $\rCLEO$, and also admits a homomorphism (denoted $\LEE$ in~\cite{slocal}) from $\mathcal{C}$~\cite{Lewark24:local}. The Rasmussen $s$-invariant $s_{\mathbb{F}}$ over any field $\mathbb{F}$ is a function on these local equivalence groups, but they were shown to have concordance information beyond the $s$-invariants, and $\rCLEO$ and $\rCLEO^o$ in particular are computationally useful for proving knots are not slice~\cite[Section 7]{slocal}. The group $\rCLEO^o$ was shown to be totally ordered, and Lewark showed that $\LewarkGp$ has a $\ZZ^\infty$ subgroup in the kernel of the Rasmussen $s$-invariants.

These local equivalence groups have simple definitions, which we review in Section~\ref{sec:results}. 
Still, beyond the results mentioned in the previous paragraph, their structure has been completely unknown (compare~\cite[Question 3.23 and Remark 6.11]{slocal} and~\cite[Question 7.9]{Lewark24:local}). In this paper, we answer several questions related to it:
\begin{theorem}\label{thm:Lewark-torsion}
	There is $2$-torsion in $\LewarkGp$. In fact, for each prime $p$, $\LewarkGp$ admits a $(\ZZ/2\ZZ)^\infty$ subgroup, and for different primes $p$, these subgroups are independent.
\end{theorem}

\begin{theorem}\label{thm:CLEO-torsion}
	There is a $(\ZZ/2\ZZ)^\infty$ subgroup of the kernel of the quotient map $\rCLEO\to \LewarkGp$.
\end{theorem}

\begin{theorem}\label{thm:Z-infty}
	For each prime $p$, the group $\LewarkGp$ admits a $\ZZ^\infty$ direct summand in the kernel of all $s$-invariants. That is, there is a surjective homomorphism
	\[
	\LewarkGp\to \ZZ^\infty\oplus \bigoplus_{\text{primes }p}\ZZ^\infty,
	\]
	where the map to the first $\ZZ^\infty$ is given by the sequence of $s$-invariants $\OneHalf(s_\QQ,s_{\FF_2}-s_\QQ,s_{\FF_3}-s_\QQ,s_{\FF_5}-s_\QQ,\dots)$.
\end{theorem}

\begin{theorem}\label{thm:ker-Z-infty}
	There is a $\ZZ^\infty$ direct summand of the kernel of the quotient map $\rCLEO\to \LewarkGp$ (and, in particular, in the kernel of all $s$-invariants). Moreover, there are at least $6$ homomorphisms from the smooth concordance group $\mathcal{C}$ to $\ZZ$ that factor through $\rCLEO$ and are linearly independent of each other, $s$-invariants for all characteristics, $\tau$, and all of Dai--Hom--Stoffregen--Truong's homomorphisms.
\end{theorem}
(This is proved as Corollaries~\ref{cor:odd-homoms} and~\ref{cor:new-homoms}, below.)

Note that neither earlier paper gave any new homomorphisms (i.e., beyond the $s$-invariants) from these local equivalence groups to $\ZZ$, so Theorems~\ref{thm:Z-infty} and~\ref{thm:ker-Z-infty} are a significant improvement of our understanding. These maps to $\ZZ^\infty$ are explicit. (Whether the local equivalence groups give rise to new concordance homomorphisms was posed explicitly in~\cite[Remark 6.11]{slocal}.) More data about knots and these concordance homomorphisms is in Section~\ref{sec:knots}.

The proofs of Theorems~\ref{thm:Lewark-torsion} and~\ref{thm:CLEO-torsion} amount to exhibiting some examples and then using relatively simple invariants to prove they (and products of them) are nontrivial. The proof of Theorem~\ref{thm:Z-infty} is more interesting. For that argument, we imitate Dai--Hom--Stoffregen--Truong's strategy in knot Floer homology~\cite{DHST21:more-conc}. Key to their argument is finding a canonical representative for each local equivalence class of complexes over $\FF_2[U,V]/(UV)$. In our setting, we instead analyze local equivalence over $\ZZ_p[h]/(ph)$, which is again a fiber product of two DVRs. We show that the local equivalence group over $\ZZ_p[h]/(ph)$ is totally ordered with an explicit $\mathbb{Z}^\infty$-summand, with maps defined similarly to Dai--Hom--Stoffregen--Truong's; base change then gives the desired maps from $\LewarkGp$. In fact, for $p=2$, we show that this Khovanov-type equivalence groups are isomorphic to Dai--Hom--Stoffregen--Truong's. In the process, we give a short re-proof of Dai--Hom--Stoffregen--Truong's structure theorem, as well as several other results in the literature (see Remark~\ref{rem:Bondarenko-is-great}).

\subsection*{AI}
Except for some quotations \textcolor{blue}{specifically highlighted below}, the tables at the end, and some of the TikZ diagrams, this paper was entirely written by the listed author. The research, however, was largely carried out by various Claude models, with modest prompting by the author. (Except where specifically noted below, all results were checked by the author.) See Section~\ref{sec:methodology} for further discussion of the methodology and some reflections on the process and ethics. The mathematical content in the paper is in Sections~\ref{sec:results} and~\ref{sec:knots}, and some open questions are in Section~\ref{sec:questions}.

\subsection*{Acknowledgments} The human author thanks Nathan Dunfield, Dirk Sch\"utz, Matthew Stoffregen, and Ian Zemke for helpful conversations.

\section{More results and all the proofs}\label{sec:results}

\subsection{Basic definitions}
Let $\ZZ[h]$ be the graded ring of polynomials in a single variable $h$, where $h$ has grading $-2$ (and $\ZZ$ lies in grading $0$). Throughout this section, fix a commutative, graded $\ZZ[h]$-algebra $R$.
Write $h^{-1}R$ for the result of inverting $h$. Some examples of rings $R$ we will consider later in the paper are $R=\ZZ[h]$, $R=\ZZ_p[h]$ (where $\ZZ_p$ is the $p$-adic integers), and $R=\ZZ_p[h]/(ph)$. In these three examples, 
$h^{-1}R$ is $\ZZ[h,h^{-1}]$, $\ZZ_p[h,h^{-1}]$, and $\FF_p[h,h^{-1}]$, respectively. We will typically call the grading on $R$ the \emph{quantum grading}.

\begin{definition}\label{def:knot-like-cx}\cite[Definition 3.1]{Lewark24:local}
  A \emph{knot-like complex over $R$} is a finitely generated cochain complex $D^{i,*}$ of free graded $R$-modules such that $h^{-1}D^{*,*}\simeq h^{-1}R$, where $h^{-1}R$ is supported in homological grading $0$ and even quantum gradings.
\end{definition}

The motivation comes from the fact that the reduced Bar-Natan complex of a knot (over a ground ring $R$) is a knot-like complex. (See, for instance,~\cite[Proof of Proposition 2.1]{LS22:mixed} for more details.)

We typically suppress the bigradings from the notation, and refer to a knot-like complex $D$.

\begin{definition}\label{def:D-side-local-equiv}
  Let $D$ and $D'$ be knot-like complexes over $R$. A chain map $\beta\co D\to D'$ over $R$ is a \emph{local map} if the induced map $\beta\co h^{-1}D\to h^{-1}D'$ is a chain homotopy equivalence. Write $D\leq D'$ if there is a local map $\beta\co D\to D'$. We say that $D$ and $D'$ are \emph{locally equivalent} if $D\leq D'$ and $D'\leq D$. Let $D\otimes D'$ denote the tensor product of $D$ and $D'$ over $R$.
\end{definition}

It is clear that: $\leq$ is reflexive and transitive; a tensor product of knot-like complexes is knot-like; tensor product respects $\leq$; and up to local equivalence (in fact, isomorphism) tensor product is associative and commutative. Let $\LewarkGp(R)$ be the set of local equivalence classes of knot-like complexes over $R$. So, tensor product makes $\LewarkGp(R)$ into a partially ordered abelian monoid, but in fact it is an abelian group: the inverse to a knot-like complex is its dual over $R$ (see~\cite[Theorem 1.1]{Lewark24:local} or specialize~\cite[Proposition 3.12]{slocal}). If we do not specify $R$, then we have in mind $R=\ZZ[h]$, so $\LewarkGp=\LewarkGp(\ZZ[h])$.

Fix another graded ring $Q$ and an identification (isomorphism) $Q/(2)\cong R/(2,h)$.

\begin{definition}\label{def:LEO-triple}\cite[Definition 3.10]{slocal}
  A \emph{LEO triple over $(Q,R)$} is a triple $(C,D,f)$ where $C$ is a finitely generated, free chain complex over $Q$, $D$ is a knot-like chain complex over $R$, and
  \[
    f\co C\otimes_Q Q/(2)\to D\otimes_R R/(2,h)
  \]
  is a chain homotopy equivalence.
\end{definition}
(Later, we will write $C\otimes_Q Q/(2)$ as $C/(2)$ and $D\otimes_R R/(2,h)$ as $D/(2,h)$.)

The term LEO stands for ``local even-odd''; the motivating example is the triple of the reduced odd Khovanov complex of a knot (from~\cite{ORSz13:odd}), the reduced Bar-Natan complex, and the identity map identifying their reductions, for $Q=\ZZ$ and $R=\ZZ[h]$; we denote this triple $\LEO(K)$. When we introduced them, we called LEO triples from Definition~\ref{def:LEO-triple} ``reduced LEO triples''; we will not consider unreduced triples here (partly in light of~\cite[Theorem 6.3]{slocal}), so drop the word ``reduced''. We can also take even Khovanov homology in place of odd Khovanov homology, giving a triple $\LEE(K)$; the invariant $\LEE(K)$ is completely determined by the Bar-Natan complex, over $\ZZ[h]$.

\begin{definition}\label{def:LEO-local-equiv}
  Let $(C,D,f)$ and $(C',D',f')$ be LEO triples over $(Q,R)$. A \emph{local map} from $(C,D,f)$ to $(C',D',f')$ consists of chain maps $\alpha\co C\to C'$ and $\beta\co D\to D'$ (over $Q$ and $R$, respectively) so that $\beta$ is a local map in the sense of Definition~\ref{def:D-side-local-equiv} and (the middle square of) the following diagram commutes up to chain homotopy:
  \begin{equation}\label{eq:local-map}
    \begin{tikzcd}
      C\arrow[r]\arrow[d,"\alpha"] & C/(2)\arrow[r,"f"]\arrow[d,"\alpha"] & D/(2,h)\arrow[d,"\beta"] & D\arrow[l]\arrow[d,"\beta"]\\
      C\arrow[r] & C/(2)\arrow[r,"f'"] & D'/(2,h) & D'\arrow[l].
    \end{tikzcd}
  \end{equation}
  Write $(C,D,f)\leq (C',D',f')$ if there is a local map from $(C,D,f)$ to $(C',D',f')$. If $(C,D,f)\leq (C',D',f')\leq (C,D,f)$, then $(C,D,f)$ and $(C',D',f')$ are \emph{locally equivalent}. Let $\rCLEO(Q,R)$ be the set of local equivalence classes of LEO triples.
\end{definition}

Define an operation $\otimes$ on $\rCLEO$ by $(C,D,f)\otimes(C',D',f')=(C\otimes_{Q}C',D\otimes_{R}D',f\otimes f')$. Again, these definitions have the following properties: $\leq$ is a partial order; tensor products of LEO triples are LEO triples; tensor product respects $\leq$; and up to local equivalence (in fact, isomorphism) tensor product is associative and commutative. So, tensor product makes $\rCLEO(Q,R)$ into a partially ordered abelian monoid, but in fact it is an abelian group, with inverse again given by the dual~\cite[Proposition 3.12]{slocal}. The basic case for knots is $Q=\ZZ$, $R=\ZZ[h]$, so if we do not specify $(Q,R)$ then we mean these rings, and we let $\rCLEO=\rCLEO(\ZZ,\ZZ[h])$.

The motivation for Definitions~\ref{def:D-side-local-equiv} and~\ref{def:LEO-local-equiv} is that $\LEO(K)$, $\LEE(K)$, and the reduced Bar-Natan complex $\KhCxBN(K)$ induce homomorphisms from the smooth concordance group to the local equivalence groups. (Lewark uses the notation $S(K)$ for his map from the concordance group to $\LewarkGp(\ZZ[h])$, but since the point of this paper is that $S(K)$ is rather far from Rasmussen's $s$-invariants, we have not adopted that notation.)

We record a few obvious properties of these constructions:
\begin{lemma}\label{lem:obvious}
  \begin{enumerate}[label=(\arabic*)]
  \item For any ring $R$, the group $\LewarkGp(R)$ is a direct summand of the group $\rCLEO(R/(h),R)$. Moreover, the maps from the concordance group to the local equivalence groups are compatible with this splitting, in the sense that the following diagrams commute: 
	\[
	\begin{tikzcd}
		\mathcal{C}\arrow[d,swap,"\LEO"]\arrow[dr,"\KhCxBN"] & & & \mathcal{C}\arrow[dl,swap,"\KhCxBN"]\arrow[d,"\LEE"]\arrow[dr,"\KhCxBN"] & \\
		\rCLEO(\ZZ,\ZZ[h])\arrow[r, two heads] & \LewarkGp(\ZZ[h]) & \LewarkGp(\ZZ[h])\arrow[r,hook]\arrow[rr,"\Id", bend right] & \rCLEO(\ZZ,\ZZ[h])\arrow[r, two heads] & \LewarkGp(\ZZ[h])
	\end{tikzcd}
	\]
  \item\label{item:base-change} Given a $\ZZ[h]$-algebra homomorphism $k\co R\to R'$, there is an induced group homomorphism $k_*\co \LewarkGp(R)\to\LewarkGp(R')$, respecting the partial order. Similarly, given a ring homomorphism $j\co Q\to Q'$ and $\ZZ[h]$-algebra homomorphism $k\co R\to R'$ inducing the same map between $Q/(2)=R/(2,h)$ and $Q'/(2)=R'/(2,h)$, there is an induced group homomorphism $(j,k)_*\co \rCLEO(Q,R)\to\rCLEO(Q',R')$, again respecting the partial order.
  \end{enumerate}
\end{lemma}
\begin{proof}
  This is immediate from the definitions. In particular, the maps between $\LewarkGp(R)$ and $\rCLEO(R/(h),R)$ are given by $D\mapsto (D/(h),D)$ and $(C,D,f)\mapsto D$. The base-change homomorphisms in Point~\ref{item:base-change} send $D$ to $D\otimes_R R'$ and $(C,D,f)$ to $(C\otimes_QQ',D\otimes_R R',f\otimes\Id)$. It is immediate from the definitions that the images of knot-like complexes (respectively LEO triples) are knot-like complexes (respectively LEO triples) and that these maps respect local maps and tensor products. 
\end{proof}

Of course, the group homomorphism from Point~\ref{item:base-change} is typically neither injective nor surjective.

\begin{definition}
  Let $\FF$ be a field. If $C$ is a knot-like complex over $\FF[h]$, then the homology of $C$ decomposes (non-canonically) as a direct sum $\FF[h]^{0,q}\oplus\mathcal{T}$ where $\mathcal{T}$ is annihilated by $h^n$ for $n\gg0$. Here, $(0,q)$ is the bigrading of $1\in\FF[h]$. The \emph{$\mathbb{F}$-valued $s$-invariant} is the homomorphism $s_\FF\co \LewarkGp(\FF[h])\to 2\ZZ$ defined by $s_\FF(C)=q$.
  
  More generally, given a homomorphism of $\ZZ[h]$-algebras $\phi\co R\to \mathbb{F}[h]$, the $s$-invariant $s_\phi$ is $s_\FF(\phi_*(C))$. We will generally denote $s_\phi$ as $s_{\mathbb{F}}$ when the homomorphism $\phi$ is clear from context.
\end{definition}

We can also talk about $s$-invariants of LEO triples, by applying the map $\rCLEO\to \LewarkGp$.

\subsection{Torsion in \texorpdfstring{$\LewarkGp$}{Lewark's local equivalence group}}\label{sec:Lewark-torsion}

We start with a basic example of $2$-torsion in $\LewarkGp(\ZZ[h])$, and then introduce systematic invariants and generalize it to prove Theorem~\ref{thm:Lewark-torsion}.

\begin{figure}
  \centering
  \begin{tikzpicture}[>=Stealth,xscale=4]
    \begin{scope}
      \node (g2) at (-1,-2)    {$x$};
      \node (g0) at (0,0)  {$y_0$};
      \node (g3) at (0,-2)  {$y_-$};
      \node (g4) at (0,2) {$y_+$};
      \node (g1) at (1,2)  {$z$};

      \node[below = 0cm of g2] {\scriptsize$(-1,-2)$};
      \node[below = 0cm of g0] {\scriptsize$(0,0)$};
      \node[below = 0cm of g3] {\scriptsize$(0,-2)$};
      \node[below = 0cm of g4] {\scriptsize$(0,2)$};
      \node[below = 0cm of g1] {\scriptsize$(1,2)$};

      \draw[->] (g2) -- node[below right=-.2cm]  {$4h$}  (g0);
      \draw[->] (g2) -- node[below]  {$8$}   (g3);
      \draw[->] (g2) -- node[above]  {$h^2$} (g4);
      \draw[->] (g0) -- node[above] {$-4h$} (g1);
      \draw[->] (g3) -- node[below right=-.2cm] {$h^2$} (g1);
      \draw[->] (g4) -- node[above] {$8$}   (g1);
    \end{scope}
    \begin{scope}[shift={(2,0)}]
      \node (g2) at (-1,-2)    {$x$};
      \node (g0) at (0,0)  {$y_0$};
      \node (g3) at (0,-2)  {$y_-$};
      \node (g4) at (0,2) {$y_+$};
      \node (g1) at (1,2)  {$z$};

      \node[below = 0cm of g2] {\scriptsize$(-1,-2n)$};
      \node[below = 0cm of g0] {\scriptsize$(0,0)$};
      \node[below = 0cm of g3] {\scriptsize$(0,-2n)$};
      \node[below = 0cm of g4] {\scriptsize$(0,2n)$};
      \node[below = 0cm of g1] {\scriptsize$(1,2n)$};

      \draw[->] (g2) -- node[below right=-.2cm]  {$2p^mh^n$}  (g0);
      \draw[->] (g2) -- node[below]  {$2p^{2m}$}   (g3);
      \draw[->] (g2) -- node[above]  {$h^{2n}$} (g4);
      \draw[->] (g0) -- node[above] {$-2p^mh^n$} (g1);
      \draw[->] (g3) -- node[below right=-.2cm]  {$h^{2n}$} (g1);
      \draw[->] (g4) -- node[above] {$2p^{2m}$}   (g1);
    \end{scope}
  \end{tikzpicture}
  \caption{\textbf{A specific 2-torsion complex in $\LewarkGp(\ZZ[h])$ and a family.} Left: the chain complex $D_{2,1,1}$. Right: its generalization $D_{p,m,n}$. Claude called these \emph{bowtie complexes}, perhaps because of the count of generators per quantum grading.}
  \label{fig:tor-in-Zh-first-example}
\end{figure}
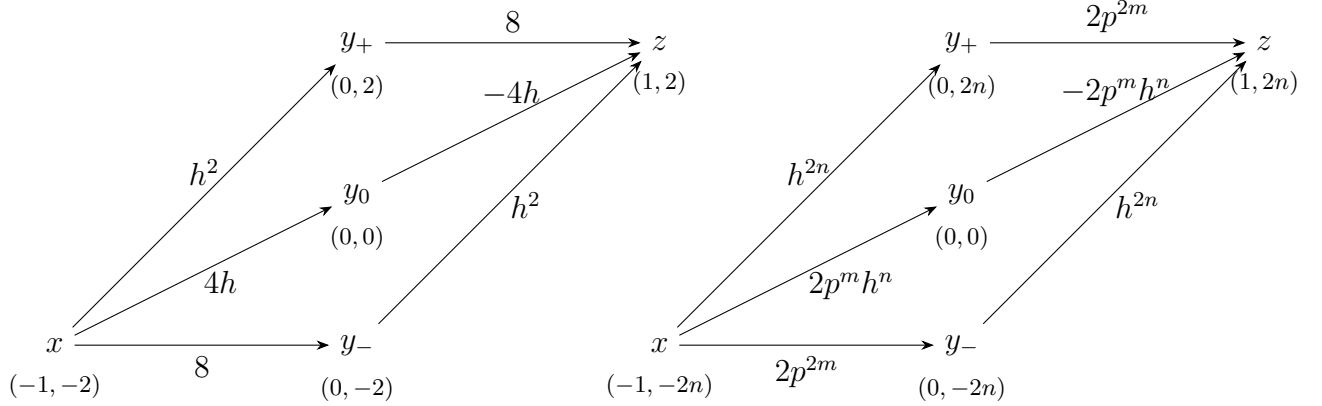

Let $D_{2,1,1}$ be the chain complex
freely generated over $\ZZ[h]$ by the following five generators:
\begin{center}
  \begin{tabular}{lccccc}
    Generator: & $x$ & $y_+$ & $y_0$  & $y_-$ & $z$\\
    Bigrading: & $(-1,-2)$ & $(0,2)$ & $(0,0)$ & $(0,-2)$ & $(1,2)$ 
	\end{tabular}
\end{center}
and differential
\begin{align*}
  \bdy(x) &= 4hy_0 + h^{2}y_+ + 8y_- & \bdy(y_0) &= -4hz\\
  \bdy(y_+) &= 8z & \bdy(y_-) &= h^{2}z & \bdy(z)&=0.
\end{align*}
See Figure~\ref{fig:tor-in-Zh-first-example}.

\begin{lemma}\label{lem:D211}
  The complex $D_{2,1,1}$ is a nontrivial $2$-torsion element of $\LewarkGp(\ZZ[h])$.
\end{lemma}

Before proving this, we recall a local equivalence invariant that detects the trivial class. Given a knot-like complex $D$ over $\ZZ[h]$, let 
\[
s_c(D)=\max\bigl\{j\mid \exists \text{ cycle }a\in D^{0,j}\text{ s.t.\ } [a]\text{ generates }H(h^{-1}D)\cong\ZZ[h,h^{-1}]\bigr\}.
\] 
(Here, $H$ denotes taking homology of the chain complex.)
It is easy to see that $s_c$ is a local equivalence invariant. Moreover, $D$ is locally equivalent to the trivial complex if and only if $s_c(D)=s_c(D^*)=0$; this is a simple version of the corresponding result for LEO triples~\cite[Theorem 5.1]{slocal}, and in fact follows from that result by applying the inclusion $\LewarkGp(\ZZ[h])\hookrightarrow \rCLEO(\ZZ,\ZZ[h])$ from Lemma~\ref{lem:obvious}.

\begin{proof}[Proof of Lemma~\ref{lem:D211}]
	The complex $D_{2,1,1}$ is knot-like: after inverting $h$, one can cancel the contractible subcomplex spanned by $y_-$ and $z$ and the contractible quotient complex spanned by $x$ and $y_+$, leaving only $\ZZ[h,h^{-1}]\langle y_0\rangle$. (That is, one can cancel both arrows labeled $h^2$.)

	To see that $D_{2,1,1}$ is not locally equivalent to the trivial complex, observe that the only cycles in bigrading $(0,0)$ have the form $k(hy_++2y_0)$, which is $2$ times a generator of $H(h^{-1}D)$. So, $s_c(D_{2,1,1})\neq 0$. (In fact, $s_c(D_{2,1,1})=-2$, witnessed by $hy_0+4y_-$.)

	The complex $D_{2,1,1}$ is self-dual, via the identification exchanging $x$ with $z$, $y_+$ with $y_-$, and $y_0$ with $-y_0$. (In the figure, this corresponds to rotating by $\pi$ around $y_0$, reversing all the arrows, and negating $y_0$.) Since the dual is the inverse in $\LewarkGp$, it follows that $D_{2,1,1}$ has order $2$.
\end{proof}

We now give a three-parameter generalization of $D_{2,1,1}$. Fix a prime $p$ and positive integers $m,n$. Let $D_{p,m,n}$ be the complex
freely generated over $\ZZ[h]$ by the following five generators:
\begin{center}
  \begin{tabular}{lccccc}
		\toprule
    Generator: & $x$ & $y_+$ & $y_0$  & $y_-$ &  $z$\\
    Bigrading: & $(-1,-2n)$ & $(0,2n)$ & $(0,0)$ & $(0,-2n)$ & $(1,2n)$\\
		\bottomrule
	\end{tabular}                                                         \end{center}
and differential
\begin{align*}
  \bdy(x) &= 2p^{m}h^ny_0 + h^{2n}y_+ + 2p^{2m}y_- & \bdy(y_0) &= -2p^mh^nz\\
  \bdy(y_+) &= 2p^{2m}z & \bdy(y_-) &= h^{2n}z & \bdy(z)&=0.
\end{align*}
Again, see Figure~\ref{fig:tor-in-Zh-first-example}. The complex $D_{p,m,n}$ is knot-like and self-dual by the same argument as in the proof of Lemma~\ref{lem:D211}. Explicitly, there is a chain homotopy equivalence $\psi=\psi_{p,m,n}\co h^{-1}D_{p,m,n}\to\ZZ[h,h^{-1}]$ given by
\begin{equation}\label{eq:psi-p-m-n}
	\begin{split}
		\psi(y_+) &= -2p^mh^{-n} \qquad\qquad\qquad \psi(y_0) = 1\\
		\psi(y_-)&=\psi(z)=\psi(x)=0.
	\end{split}
\end{equation}

To prove Theorem~\ref{thm:Lewark-torsion}, we change base to work over $\ZZ_p[h]$, where $\ZZ_p$ are the $p$-adics. (We could equally well replace $\ZZ_p$ by $\ZZ_{(p)}$, the integers localized at $p$, for the discussion here, but the completion will be useful below.) First, consider the complex $D_{p,m,n}$ over $\ZZ_\ell[h]$ where $\ell\neq p$. If $\ell$ is odd, then $2p^{2m}$ is invertible and the complex is homotopy equivalent to $\ZZ_\ell[h]$. For $\ell=2$, there is a local map from $\ZZ_\ell[h]$ to $D_{p,m,n}$ sending $1$ to $y_0+h^ny_+/p^m$ and a local map from $D_{p,m,n}$ to $\ZZ_\ell[h]$ dual to this one, showing that $D_{p,m,n}$ is trivial in $\LewarkGp(\ZZ_\ell[h])$. So, if we let $G$ be the subgroup of $\LewarkGp(\ZZ[h])$ generated by all the $D_{p,m,n}$ for all $p,m,n$, and let $G_p$ be the subgroup of $\LewarkGp(\ZZ_p[h])$ generated by the $D_{p,m,n}$ for all $m,n$ (but with $p$ fixed), then the base change maps from $\ZZ[h]$ to the various $\ZZ_p[h]$ assemble to a surjection
\begin{equation}\label{eq:base-change-to-all-padics}
G\onto \bigoplus_{\text{primes }p}G_p.
\end{equation}

Thus, to prove Theorem~\ref{thm:Lewark-torsion}, it remains to prove that $G_p\cong (\ZZ/2\ZZ)^\infty$. More precisely, we will prove:
\begin{proposition}\label{prop:Lewark-torsion-key-step}
	In $\LewarkGp(\ZZ_p[h])$, the only relations between the complexes $D_{p,m,n}$ are that each $D_{p,m,n}$ is $2$-torsion. That is, $G_p=\bigoplus_{m,n=1}^\infty\ZZ/2\ZZ$.
\end{proposition}
\begin{proof}
	We use the following invariant. Fix a knot-like complex $D$ over $\ZZ_p[h]$ and an identification $h^{-1}H(D)\cong \ZZ_p[h,h^{-1}]$. Given a class $[x]\in H^{0,0}(D)$, let $i(x)$ be its image in $\ZZ_p$, the $0$-graded part of $\ZZ_p[h,h^{-1}]$, under the composition $H(D)\to h^{-1}H(D)\stackrel{\cong}{\lra} \ZZ_p[h,h^{-1}]$. Let $I(D)$ be the ideal generated by all $i(x)$ and let $v_p(D)\in\{0,1,2,\dots,\infty\}$ be the largest $k$ so that $I(D)\subset p^k\ZZ_p$. It is easy to see that $v_p(D)$ is a local equivalence invariant.

	Now, let $D=D_{p,m_1,n_1}\otimes\cdots\otimes D_{p,m_\ell,n_\ell}$, where the pairs $(m_1,n_1),\dots,(m_\ell,n_\ell)$ are all distinct. We compute $v_p(D)$. 

	Shortening notation, write $D_k=D_{p,m_k,n_k}$, and let $x^k$, $y_+^k$, and so on denote its generators, and $\psi_k$ the map from Equation~\eqref{eq:psi-p-m-n}. Then $\psi=\psi_1\otimes\cdots\otimes\psi_\ell\co h^{-1}D\to \ZZ_p[h,h^{-1}]$ is a chain homotopy equivalence. Call a monomial $h^a b_1\otimes \cdots \otimes b_\ell\in D$, where each $b_i\in\{x^i,y_+^i,y_0^i,y_-^i,z^i\}$, \emph{substantial} if in fact each $b_i$ is either $y_0$ or $y_+$. Since $\psi_k(y_-^k)=\psi_k(z^k)=\psi_k(x^k)=0$, for any cocycle $c$, $\psi(c)$ is determined by the coefficients in $c$ of the substantial monomials. Moreover, from the gradings, $\psi$ sends any substantial monomial in bigrading $(0,0)$, except for $y_0\otimes\cdots\otimes y_0$, to $p\ZZ_p$. 
	
	So, it suffices to prove that in any cycle in $D$ with bigrading $(0,0)$, the coefficient of $y_0^1\otimes\cdots\otimes y_0^\ell$ is divisible by $p$. Assume not, so the coefficient of $y_0^1\otimes\cdots\otimes y_0^\ell$ is a unit, and we may assume the coefficient is $1$. We have
	\begin{equation}\label{eq:just-try-to-cancel-me}
		\bdy(y_0^1\otimes\cdots\otimes y_0^\ell) = -2(p^{m_1}h^{n_1}z^1\otimes y_0^2\otimes\cdots\otimes y_0^\ell + p^{m_2}h^{n_2}y_0^1\otimes z^2\otimes \cdots\otimes y_0^\ell+\cdots).
	\end{equation}
	Order the $D_i$ by increasing $n_i$ and, for fixed $n_i$, by increasing $m_i$. Suppose that $n_1=\cdots=n_j$ and $m_1<m_2<\cdots<m_j$. Any other term in $D$ whose differential contains a multiple of $z^1\otimes y_0^2\otimes\cdots\otimes y_0^\ell$ is of one of the forms:
	\begin{itemize}
	\item $p^\alpha h^\beta y_+^1\otimes y_0^2\otimes\cdots\otimes y_0^\ell$,
	\item $p^\alpha h^\beta y_-^1\otimes y_0^2\otimes\cdots\otimes y_0^\ell$,
	\item $p^\alpha h^\beta z^1\otimes y_0^2\otimes\cdots\otimes x^i\otimes y_0^{i+1}\otimes\cdots\otimes y_0^\ell$. 
	\end{itemize}
	In the first case, the coefficient of $z^1\otimes y_0^2\otimes\cdots\otimes y_0^\ell$ in the differential of this term has $p$-power at least $2m_1$, so does not cancel $p^{m_1}h^{n_1}z^1\otimes y_0^2\otimes\cdots\otimes y_0^\ell$. In the second case, the relevant coefficient in the differential of this term has $h$-power at least $2n_1$, and so does not cancel $p^{m_1}h^{n_1}z^1\otimes y_0^2\otimes\cdots\otimes y_0^\ell$. (This case is also excluded for grading reasons.) In the third case, if $i\leq j$, then the relevant coefficient in the differential of this term has $p$-power strictly greater than $m_1$, while if $i>j$ it has $h$-power strictly greater than $n_1$. Thus, no other term cancels the first term on the right of Equation~\eqref{eq:just-try-to-cancel-me}, a contradiction.
\end{proof}

\begin{remark}
	Claude's original proof of Proposition~\ref{prop:Lewark-torsion-key-step} has more the style of a clever, elegant proof a combinatorist might write, but is also a little longer; the core of its proof and the one given above are the same.
\end{remark}

We conclude this section with some structural results about $\LewarkGp(\ZZ_p[h])$, relevant to searches for further torsion.

\begin{lemma}\label{lem:direct-sum-of-indecomposables}
  Over $\ZZ_p[h]$, any knot-like complex decomposes as a finite direct sum $D=\bigoplus_{i}D_i$ of indecomposable complexes, exactly one of which (say, $D_0$) is knot-like and the others of which become acyclic after inverting $h$. In particular, $D$ is locally equivalent to $D_0$.
\end{lemma}
Here, by indecomposable we mean that no $D_i$ splits further as a direct sum.
\begin{proof}
	Since $\ZZ_p[h]$ is a graded local ring, by the graded version of Kaplansky's Theorem \cite{Kaplansky58:theorem}, any summand of a free module is free. So, over the residue field $\FF_p=\ZZ_p[h]/(p,h)$, each $D_i\otimes_{\ZZ_p[h]}\FF_p$ has dimension at least $1$, and hence any direct sum decomposition has at most $\dim D\otimes_{\ZZ_p[h]}\FF_p$-many summands. So, pick any maximal direct sum decomposition $D=\bigoplus_{i}D_i$. Then $\ZZ_p[h,h^{-1}]\simeq h^{-1}D=\bigoplus_i h^{-1}D_i$. Let $E_i$ be the image of $H(h^{-1}D_i)$ in $\ZZ_p[h,h^{-1}]$ under this homotopy equivalence. Then $\ZZ_p[h,h^{-1}]\cong \bigoplus_i E_i$, but since $\ZZ_p[h,h^{-1}]$ does not decompose as a nontrivial direct sum, some $E_i$ is all of $\ZZ_p[h,h^{-1}]$ and the others vanish. The local equivalence to that $D_i$ is exhibited by the inclusion and projection to that summand.
\end{proof}

By Lemma~\ref{lem:direct-sum-of-indecomposables}, there is no loss of information if we restrict to indecomposable knot-like complexes. The following proposition is an analogue of the existence of the connected complex in involutive Heegaard Floer homology (see, e.g.,~\cite{Stoffregen20:Seifert,HHL21:conn-invol}). 

\begin{proposition}\label{prop:Zp-structure-thm}
	Suppose that $D$ and $D'$ are indecomposable knot-like complexes over $\ZZ_p[h]$. If $D$ is locally equivalent to $D'$, then $D$ is isomorphic to $D'$.
\end{proposition}
\begin{proof}
	Let $f\co D\to D'$ and $g\co D'\to D$ be local maps. Their composition $g\circ f$ induces a map $h^{-1}D\to h^{-1}D$, and since $h^{-1}D\simeq \ZZ_p[h,h^{-1}]$, a map $\ZZ_p[h,h^{-1}]\to \ZZ_p[h,h^{-1}]$ (of modules over $\ZZ_p[h,h^{-1}]$). Indeed, from the definitions of knot-like complexes and local maps, this induced map  $\ZZ_p[h,h^{-1}]\to \ZZ_p[h,h^{-1}]$ is an isomorphism, and hence is given by multiplication by a unit $\epsilon$ in $\ZZ_p$. Replace $g$ by $\epsilon^{-1}g$, so that $g\circ f$ instead induces the identity map $\ZZ_p[h,h^{-1}]\to \ZZ_p[h,h^{-1}]$ and, in particular, the identity map $\FF_p[h,h^{-1}]\to \FF_p[h,h^{-1}]$.
	
	Write $u=g\circ f$. We use completeness of $\ZZ_p$ and a Hensel's Lemma-style lifting argument to arrange that $u$ is an idempotent. Let $\ol{u}\co D/pD\to D/pD$ be the induced map. Note that, in each grading, $D/pD$ is a finite-dimensional $\FF_p$-vector space and, in particular, a finite set. Fix a set of generators for $D$ (which also gives a set of generators for $D/pD$). The map $\ol{u}$ is determined by the images of these generators, each of which lies in a finite set. It follows that there are integers $n>m>0$ so that $\ol{u}^n=\ol{u}^m$, and therefore $\ol{u}^{n+k}=\ol{u}^{m+k}$ for all $k$. Thus, for some large $n$, $\ol{u}^n=\ol{u}^{2n}$, i.e., $\ol{u}^n$ is idempotent. Let $v_0=u^n$, so $v_0^2\equiv v_0\pmod{p}$. (Note that exponents here mean composing a map with itself many times.)

	For each $k$, let $v_{k+1}=3v_{k}^2-2v_{k}^3$. We prove by induction that $v_{k+1}\equiv v_{k}\pmod{p^{2^k}}$ and $v_{k+1}^2\equiv v_{k+1}\pmod{p^{2^{k+1}}}$. Indeed,
	\begin{align*}
		v_{k}-v_{k+1}&=(v_k^2-v_k)(2v_k-1)\equiv 0\pmod{p^{2^k}}\\
		v_{k+1}^2-v_{k+1}&=9v_k^4-12v_k^5+4v_k^6-3v_k^2+2v_k^3=(v_k^2-v_k)^2(4v_k^2-4v_k-3)\equiv 0\pmod{(p^{2^k})^2}.
	\end{align*}
	Thus, by completeness, the limit $v$ of the $v_k$ exists, and $v^2=v$. Moreover, at each stage of the iteration, the composition $\FF_p[h,h^{-1}]\stackrel{\simeq}{\lra} h^{-1}D/pD\stackrel{\ol{v}_k}{\lra} h^{-1}D/pD\stackrel{\simeq}{\lra} \FF_p[h,h^{-1}]$ sends $1$ to $1$. So, $v\neq 0$. However, $D\cong \ker(v)\oplus \ker(\Id-v)$, so $D$ being indecomposable implies that $\Id-v=0$.

	Now, $v\equiv u^n\pmod{p}$. Thus, $u^n-\Id\in pD$. So, $\sum_{k\geq 0}(\Id-u^n)^k$ is a convergent series, and is an inverse to $u^n$. Since $u^n$ is invertible, $u$ is invertible (with inverse $u^{-1}=u^{n-1}(u^n)^{-1}$). We have $(u^{-1}g)\circ f=\Id_D$, so $f$ is a split injection, and $D'\cong D\oplus X$ for some complex $X$. Since $D'$ is indecomposable, $X=0$, and $D\cong D'$.
\end{proof}

In Claude's computer search for $2$-torsion (which found the examples above), it focused on complexes which were isomorphic to their duals, not merely locally equivalent. As it notes, over $\ZZ_p[h]$ this is justified by the following immediate corollary of Proposition~\ref{prop:Zp-structure-thm}:
\begin{corollary}\label{cor:self-dual}
	If $D$ is an indecomposable, knot-like complex over $\ZZ_p[h]$ which is $2$-torsion in $\LewarkGp(\ZZ_p[h])$, then $D\cong D^*$.
\end{corollary}

\begin{proposition}
  If $D$ is an indecomposable, knot-like complex representing a 2-torsion class in $\LewarkGp(\ZZ_p[h])$, then the number of generators of $D$ is congruent to $1\pmod{4}$. In particular, if $D$ is nontrivial, then $D$ has at least $5$ generators.
\end{proposition}
\begin{proof}
	Let $n_i$ be the number of generators of $D$ in homological grading $i$. By Corollary~\ref{cor:self-dual}, $n_i=n_{-i}$. On the other hand, over the field of fractions $\QQ_p(h)$ of $\ZZ_p[h]$, the homology of $D$ is $\QQ_p(h)$ (because $D$ is knot-like), so the Euler characteristic of $D$ is $1$. Thus, $\sum_{i\text{ even}}n_i=1+\sum_{i\text{ odd}}n_i$. So, $\sum_i n_i=1+4\sum_{i>0\text{ odd}}n_i$.
\end{proof}
  
\begin{remark}
  One can check that torsion complexes in $\LewarkGp(\ZZ[h])$ also need at least 5 generators.
\end{remark}

We conclude this section with a nontrivial, indecomposable, rank-9 2-torsion example:
\begin{example}
	Let $Y$ be the free $\ZZ_2[h]$-module generated by
	\begin{center}
		\begin{tabular}{lccccccc}
			\toprule
			Generator: & $g_1^-,g_2^-$ & $y_-$ & $x_-$ & $x_0$  & $x_+$ & $y_+$ & $g_1^+,g_2^+$\\
			Bigrading: & $(-1,-6)$ & $(0,-4)$ & $(0,-2)$ & $(0,0)$ & $(0,2)$ & $(0,4)$ & $(1,6)$\\
			\bottomrule
		\end{tabular}
	\end{center}
	with differentials
	\begin{align*}
  \bdy g^-_1&=4h^{3}x_0+h^{4}x_++8h^{2}x_-+16hy_- &
  \bdy g^-_2&=4h^{3}x_0+h^{5}y_++8h y_-\\
  \bdy x_-&=h^{4}g^+_1 &
  \bdy x_0&=-4h^{3}\bigl(g^+_1+g^+_2\bigr) & \bdy x_+&=8h^{2}g^+_1\\
  \bdy y_-&=h^{5}g^+_2 & \bdy y_+&=8h(2g^+_1+g^+_2) .
	\end{align*}
	(The notation here is Claude's.) The self-duality is given by $x_0\leftrightarrow -x_0$, $x_+\leftrightarrow x_-$, $y_+\leftrightarrow y_-$, $g_i^-\leftrightarrow g_i^+$.
	We leave verification that this example is knot-like and indecomposable to an LLM (Claude wrote Python code to check the complex is indecomposable, by exhaustive search).
\end{example}

\subsection{Torsion in \texorpdfstring{$\rCLEO$}{the even-odd local equivalence group}}\label{sec:CLEO-torsion}

\begin{figure}
	\centering
	\begin{tikzpicture}[>=Stealth,scale=1.0,every node/.style={inner sep=1.5pt},
                    lbl/.style={font=\scriptsize},
                    hh/.style={->,densely dashed}]
	\begin{scope}
  \node (u)  at (3,2)     {$\ol{y}$};
  \node (c)  at (6,2)     {$\ol{x}$};
  \node (b)  at (0,0)     {$z$};
  \node (e1) at (3,0.80)  {$w_1$};
  \node (e2) at (3,0)     {$w_2$};
  \node (e3) at (3,-0.80) {$w_3$};
  \node (bs) at (6,0)     {$\ol{z}$};
  \node (cs) at (0,-2)    {$x$};
  \node (up) at (3,-2)    {$y$};

  \draw[->] (b)  -- node[lbl,above,pos=.55] {$-2$} (e1);
  \draw[->] (b)  -- node[lbl,below,pos=.28] {$-2$} (e3);
  \draw[->] (e1) -- node[lbl,above,pos=.45] {$-2$} (bs);
  \draw[->] (e3) -- node[lbl,below,pos=.45] {$2$}  (bs);
  \draw[->] (cs) -- node[lbl,below] {$4$} (up);
  \draw[->] (u)  -- node[lbl,above] {$4$} (c);

  \draw[hh] (b)  to[bend left=14] node[lbl,above left,pos=.5] {$-h$} (u);
  \draw[hh] (cs) to[bend right=22] node[lbl,left,pos=.8] {$-h$} (e1);
  \draw[hh] (cs) -- node[lbl,below right,pos=.62] {$h$} (e3);
  \draw[hh] (e1) -- node[lbl,above,pos=.62] {$-h$} (c);
  \draw[hh] (e3) to[bend right=16] node[lbl,above left,pos=.65] {$-h$} (c);
  \draw[hh] (up) -- node[lbl,below,pos=.60] {$-h$} (bs);

  \node[lbl] at (0,-2.95) {$-1$};
  \node[lbl] at (3,-2.95) {$0$};
  \node[lbl] at (6,-2.95) {$1$};
\end{scope}

\begin{scope}[shift={(-7,0)}]
  \node[lbl] at (-1.35,2.0) {$q=2$};
  \node[lbl] at (-1.35,0)   {$q=0$};
  \node[lbl] at (-1.35,-2.0){$q=-2$};

  \node (u)  at (3,2)     {$\ol{y}$};
  \node (c)  at (6,2)     {$\ol{x}$};
  \node (b)  at (0,0)     {$z$};
  \node (e1) at (3,0.80)  {$w_1$};
  \node (e2) at (3,0)     {$w_2$};
  \node (e3) at (3,-0.80) {$w_3$};
  \node (bs) at (6,0)     {$\ol{z}$};
  \node (cs) at (0,-2)    {$x$};
  \node (up) at (3,-2)    {$y$};

  \draw[->] (b)  -- node[lbl,above,pos=.55] {$2^n$}  (e2);
  \draw[->] (b)  -- node[lbl,below,pos=.55] {$2^n$}  (e3);
  \draw[->] (e2) -- node[lbl,above,pos=.45] {$2^n$}  (bs);
  \draw[->] (e3) -- node[lbl,below,pos=.45] {$-2^n$} (bs);
  \draw[->] (cs) -- node[lbl,below] {$4$} (up);
  \draw[->] (u)  -- node[lbl,above] {$4$} (c);

  \node[lbl] at (0,-2.95) {$-1$};
  \node[lbl] at (3,-2.95) {$0$};
  \node[lbl] at (6,-2.95) {$1$};
\end{scope}
\end{tikzpicture}
	\caption{\textbf{Two-torsion elements in the kernel.} The complex $C_n$ is on the left, and $D$ is on the right. The identification $f$ is the obvious one. In both complexes, the horizontal coordinate is the homological grading and the vertical coordinate is the quantum grading.}
	\label{fig:kernel-2-torsion}
\end{figure}
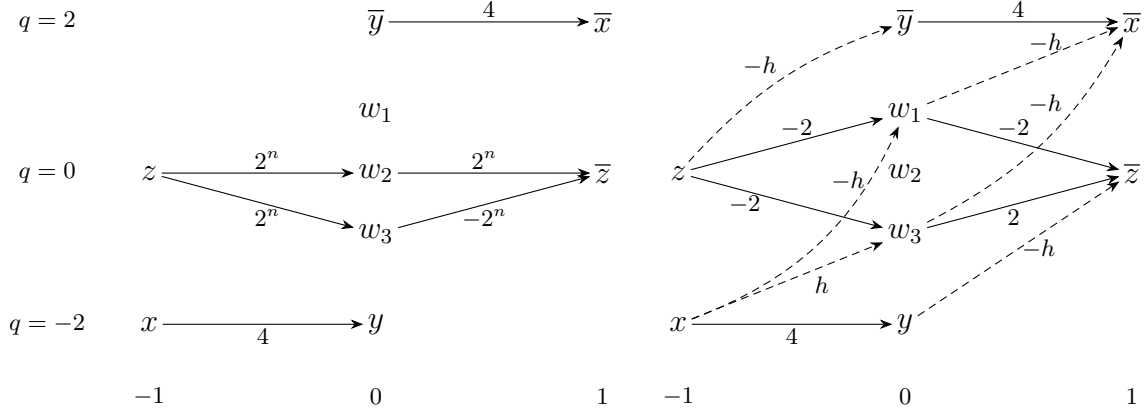

Consider the complexes $C_n$ (over $\ZZ$) and $D$ (over $\ZZ[h]$) in Figure~\ref{fig:kernel-2-torsion}. That is, $C_n$ has $9$ generators over $\ZZ$ and $D$ has $9$ generators over $\ZZ[h]$, labeled and graded (for both complexes) as follows:
\begin{center}
		\begin{tabular}{lccccccc}
			\toprule
			Generator: & $x$ & $y$ & $z$ & $w_1$, $w_2$, $w_3$ & $\ol{y}$ & $\ol{z}$ & $\ol{x}$ \\
			Bigrading: & $(-1,-2)$ & $(0,-2)$ & $(-1,0)$ & $(0,0)$ & $(0,2)$ & $(1,0)$ & $(1,2)$\\
			\bottomrule
		\end{tabular}
\end{center}
The differential on $C_n$ is given by
\begin{align*}
\bdy_C(x) &= 4y & \bdy_C(y) &= 0 & \bdy_C(z) &= 2^n(w_2+w_3)\\
\bdy_C(w_1) &= 0 & \bdy_C(w_2) &= 2^n \ol{z} & \bdy_C(w_3) &= -2^n\ol{z}\\
\bdy_C(\ol{z}) &= 0 & \bdy_C(\ol{y}) &= 4\ol{x} & \bdy_C(\ol{x}) &= 0
\end{align*}
while the differential on $D$ is given by
\begin{align*}
\bdy_D(x) &= 4y-hw_1+hw_3 & \bdy_D(y) &= -h\ol{z} & \bdy_D(z) &= -2(w_1+w_3)-h\ol{y}\\
\bdy_D(w_1) &= -2\ol{z}-h\ol{x} & \bdy_D(w_2) &= 0 & \bdy_D(w_3) &= 2\ol{z}-h\ol{x}\\
\bdy_D(\ol{z}) &= 0 & \bdy_D(\ol{y}) &= 4\ol{x} & \bdy_D(\ol{x}) &= 0
\end{align*}
The identification between $C_n/(2)$ and $D/(2,h)$ is the obvious one. There is an isomorphism between $C_n$ and $C^*_n$ (respectively $D$ and $D^*$) induced by $x\leftrightarrow \ol{x}$, $y\leftrightarrow \ol{y}$, $z\leftrightarrow \ol{z}$, $w_1\leftrightarrow w_1$, $w_2\leftrightarrow w_2$, and $w_3\leftrightarrow -w_3$. (That is, $x$ is mapped to the dual basis vector to $\ol{x}$, and so on.)

The complex $D$ is knot-like (cancel the subcomplex spanned by $w_1,y,\ol{z},\ol{x}$ and the quotient complex spanned by $x,w_3,z,\ol{y}$, leaving only $w_2$), so $(C_n,D,f)$ is a LEO triple. Since $D$ decomposes as $D'\oplus \ZZ[h]\langle w_2\rangle$ for some $D'$, $D$ by itself is locally equivalent to $\ZZ[h]$. So, $(C_n,D,f)$ is in the kernel of the map $\rCLEO(\ZZ,\ZZ[h])\to\LewarkGp(\ZZ[h])$. We already checked that the triples $(C_n,D,f)$ are $2$-torsion in $\rCLEO(\ZZ,\ZZ[h])$ (because they are self-dual).

\begin{proof}[Proof of Theorem~\ref{thm:CLEO-torsion}]
	It suffices to show that the elements $(C_n,D,f)\in\rCLEO(\ZZ,\ZZ[h])$ are all distinct: since they are $2$-torsion elements of an abelian group, they either generate a copy of $(\ZZ/2\ZZ)^N$ for some finite $N$, which has finitely many elements, or $(\ZZ/2\ZZ)^\infty$.

	Consider the image of $(C_n,D,f)$ in $\rCLEO(\ZZ/(2^m),\ZZ[h])$. If $m\leq n$, then the inclusion and projection to the summand $(\ZZ/(2^m)\langle w_2\rangle,\ZZ[h]\langle w_2\rangle)$ give a local equivalence, so $(C_n\otimes\ZZ/2^m\ZZ,D,f)=0\in \rCLEO(\ZZ/2^m\ZZ,\ZZ[h])$. We claim that if $m>n$, then $(C_n\otimes\ZZ/2^m\ZZ,D,f)\neq 0$ in $\rCLEO(\ZZ/2^m\ZZ,\ZZ[h])$. Indeed, consider a local map 
	\[
	(\alpha,\beta)\co (\ZZ/2^m\ZZ\langle r\rangle, \ZZ[h]\langle s\rangle)\to (C_n\otimes\ZZ/2^m\ZZ,D,f),
	\]
	where $r$ and $s$ are names for generators, and lie in bigrading $(0,0)$. We have $\beta(s)= a_1w_1 + a_2w_2 + a_3w_3+bh\ol{y}$ for some $a_1,a_2,a_3,b\in\ZZ$. Since the map is local, $a_2=\pm 1$, and by negating the maps if necessary we may assume $a_2=1$. Since $\beta$ is a chain map, $a_1=a_3$ and $a_1=2b$. In particular, $a_1$ and $a_3$ are both even. Similarly, $\alpha(r)=c_1w_1+c_2w_2+c_3w_3$ where $c_1,c_2,c_3\in\ZZ/2^m\ZZ$ and $c_i\equiv a_i\pmod{2}$. Thus, $c_2$ is odd and $c_3$ is even. On the other hand, $2^nc_2=2^nc_3$; since $n<m$, this is a contradiction.

	It follows that the $(C_n,D,f)$ are all distinct elements of $\rCLEO(\ZZ,\ZZ[h])$, and hence generate a $(\ZZ/2\ZZ)^\infty$ subgroup.
\end{proof}

\subsection{A new \texorpdfstring{$\ZZ^\infty$-}{infinite rank }summand}\label{sec:Z-infty-summand}

In this section, we give a $\ZZ^\infty$ summand of $\LewarkGp(\ZZ[h])$, by considering $\LewarkGp(\ZZ_p[h]/(ph))$ and imitating Dai--Hom--Stoffregen--Truong's strategy.

First, the proofs of Lemma~\ref{lem:direct-sum-of-indecomposables} and Proposition~\ref{prop:Zp-structure-thm} carry over to $\ZZ_p[h]/(ph)$ with only notational changes, so:
\begin{proposition}\label{prop:Zp-ph-structure-thm}
	Every knot-like complex $D$ over $\ZZ_p[h]/(ph)$ has a knot-like indecomposable summand $D_0$, and $D$ is locally equivalent to $D_0$. Moreover, $D$ is locally equivalent to $D'$ if and only if $D_0$ is isomorphic to $D'_0$.
\end{proposition}
(Compare~\cite[Theorem 6.1]{DHST21:more-conc}.)

The ring $\ZZ_p[h]/(ph)$ is graded local, with maximal ideal $\mathfrak{m}=(p,h)$. If $N$ is a free module over $\ZZ_p[h]/(ph)$, then any element $\xi\in \mathfrak{m}N$ decomposes uniquely as $\xi_p+\xi_h$ where $\xi_p\in pN$ and $\xi_h\in hN$. Call a complex $D$ \emph{minimal} if $\bdy(D)\subset \mathfrak{m}D$. Any indecomposable knot-like complex is minimal: if $D$ is not minimal, then some basis element maps to a unit times another, plus other terms, and one can change basis to split off an acyclic summand. In a minimal complex, the differential $\bdy$ decomposes canonically as $\bdy_p+\bdy_h$, where $\bdy_p(D)\subset pD$ and $\bdy_h(D)\subset hD$. The fact that $\bdy^2=0$ translates to
\begin{align*}
	\bdy_p^2&=0 & \bdy_h^2&=0 & \bdy_p\bdy_h+\bdy_h\bdy_p&=0.
\end{align*}
The last of these equations holds trivially, since it is in $phD=0$.
Note also that, in an indecomposable knot-like complex, all generators have even quantum grading. 

Said differently,
\[
\mathbb{Z}_p[h]/(ph) \cong \ZZ_p\times_{\FF_p}\FF_p[h]=\{(a,b)\in \ZZ_p\times \FF_p[h]\mid \ol{a} = \ol{b}\}
\]
where $\ol{a}$ is the image of $a$ in $\FF_p=\ZZ_p/(p)$ and $\ol{b}$ is the image of $b$ in $\FF_p=\FF_p[h]/(h)$. Given a minimal complex $D$, the maps $\bdy_p$ and $\bdy_h$ are maps 
\begin{align*}
\bdy_p\co D_p=D/hD\to D_p\\
\bdy_h\co D_h=D/pD\to D_h
\end{align*}
making $D_p$ and $D_h$ into finitely generated chain complexes over $\ZZ_p$ and $\FF_p[h]$, respectively. The key point is that we have an isomorphism of bigraded $\FF_p$-vector spaces $D_p/pD_p\cong D_h/hD_h$.
Conversely, given minimal chain complexes $(D_p,\bdy_p)$ and $(D_h,\bdy_h)$ over $\ZZ_p$ and $\FF_p[h]$, respectively, and an isomorphism
\[
f\co D_p/pD_p \stackrel{\cong}{\lra} D_h/hD_h
\]
of $\FF_p$-vector spaces, there is a corresponding chain complex $D$ over $\ZZ_p[h]/(ph)$. This correspondence is functorial:
\begin{lemma}\label{lem:Smith-1}
	There is an equivalence of categories between:
	\begin{enumerate}
		\item The category of minimal chain complexes over $\ZZ_p[h]/(ph)$ and chain maps.
		\item The category with objects triples $(D_p,D_h,f)$ and maps a pair of chain maps $\alpha\co D_p\to D'_p$ and $\beta\co D_h\to D'_h$ so that the following diagram commutes:
		\[
		\begin{tikzcd}
		D_p\arrow[r]\arrow[d,"\alpha"] & D_p/pD_p\arrow[d,"\alpha"] \arrow[r,"f"] & D_h/hD_h\arrow[d,"\beta"] & D_h\arrow[l]\arrow[d,"\beta"]\\
		D'_p\arrow[r] & D'_p/pD'_p \arrow[r,"f'"] & D'_h/hD'_h & D'_h.\arrow[l]
		\end{tikzcd}
		\]
	\end{enumerate}
	Moreover, a map $D\to D'$ is local if and only if the corresponding map $\beta\co D_h\to D'_h$ induces a homotopy equivalence after inverting $h$.
\end{lemma}
\begin{proof}
	This is straightforward.
\end{proof}
It is interesting, and will be fruitful later, to compare this lemma with Diagram~\eqref{eq:local-map}.

Knot-like complexes over $\ZZ_p[h]/(ph)$ still carry a modulo-$p$ $s$-invariant: from the classification of modules over the PID $\FF_p[h]$, given a knot-like complex $D$, the homology of $D\otimes_{\ZZ_p[h]/(ph)}\FF_p[h]$ is $\FF_p[h]\oplus\mathcal{T}$, where $\mathcal{T}$ is a torsion $\FF_p[h]$-module, and $s_{\FF_p}(D)$ is the quantum grading of a generator of this $\FF_p[h]$. (The homological grading of a generator is $0$.) It is clear (e.g., via Lemma~\ref{lem:obvious}\ref{item:base-change}) that $s_{\FF_p}$ is a well-defined homomorphism $\LewarkGp(\ZZ_p[h]/(ph))\to2\ZZ$.

\subsubsection{This local equivalence group is totally ordered}
In this section, we prove that the partial order $\leq$ on the local equivalence group $\LewarkGp(\ZZ_p[h]/(ph))$ is, in fact, a total order.

In the proof, it will be convenient to have another perspective on such a complex $D$. Let $d_{i,j}$ be the number of generators of $D$ in homological grading $i$ and quantum grading $j$. Then the part of $D$ in homological grading $i$ and quantum grading $q$ is given by
\begin{equation}\label{eq:vector-rep}
	\ZZ_p^{d_{i,q}}\oplus\bigoplus_{j>q}\FF_p^{d_{i,j}}.
\end{equation}
Assuming the complex is minimal, the differential $\bdy\co D^{i,q}\to D^{i+1,q}$ is determined by a $d_{i+1,q}\times d_{i,q}$ matrix $[\bdy_p^{i,q}]$ with entries in $p\ZZ_p$ and a $\bigl(\sum_{j>q}d_{i+1,j}\bigr)\times d_{i,q}$-matrix $[\bdy_h^{i,q}]$ with entries in $\FF_p$ (since $h$ has negative quantum grading). It will also be convenient to consider the $\bigl(\sum_{j}d_{i+1,j}\bigr)\times\bigl(\sum_{j}d_{i,j}\bigr)$-matrix $[\bdy_h^{i,*}]$ with entries in $\FF_p$ corresponding to $\bdy_h\co \bigoplus_q D^{i,q}\to \bigoplus_q D^{i+1,q}$.

Given $x$ in $D^{i,q}$, write $\vec{x}$ for the corresponding vector from Formula~\eqref{eq:vector-rep}.

\begin{proposition}\label{prop:total-order}
	The partial order on $\LewarkGp(\ZZ_p[h]/(ph))$ is a total order.
\end{proposition}
(Compare~\cite[Proposition 3.14]{DHST21:more-conc}.)
\begin{proof}
	It suffices to prove that, for any indecomposable knot-like complex $D$, there is either a local map from $\ZZ_p[h]/(ph)$ to $D$ or from $D$ to $\ZZ_p[h]/(ph)$. Since $s_{\FF_p}$ is a homomorphism, we may assume that $s_{\FF_p}(D)\geq 0$. Let $x\in (D\otimes\FF_p[h])^{0,q}$ be a cycle whose image generates $h^{-1}H^0(D)$, so $q=s_{\FF_p}(D)$. Choose an element $y\in D^{0,q}$ whose image in $(D\otimes\FF_p[h])^{0,q}$ is $x$. We have $\bdy_h(y)=0$, but $\bdy_p(y)$ may be nonzero. Indeed, if $y$ is a cycle then there is a local map $\ZZ_p[h]/(ph)\to D$ sending $1$ to $h^{q/2}y$, and we are done. In fact, if $q>0$, then $\bdy(h^{q/2}y)=h^{q/2}(\bdy_hy+\bdy_py)=h^{q/2}\bdy_py=0$ (because $hp=0$), and again we have a local map $\ZZ_p[h]/(ph)\to D$ sending $1$ to $h^{q/2}y$. So, it suffices to consider the case that $s_{\FF_p}(D)=0$ (i.e., $q=0$). (Claude points out that this last simplification is not really needed for the argument below.)

	Of course, $y$ is not unique. With respect to the vector representation~\eqref{eq:vector-rep} (with $q=0$), changing the preimage of $x$ corresponds to adding elements of $(p\ZZ_p)^{d_{0,0}}$. We can also add elements that vanish in $h^{-1}H^0(D)$, which corresponds to changing $\vec{y}$ by adding elements of $\Image[\bdy_h^{-1,*}]\cap \bigoplus_{j\geq 0}\FF_p^{d_{0,j}}$. (For elements with quantum grading $j=0$, there is a choice of promotion to $\ZZ_p$, but this is already included in the $(p\ZZ_p)^{d_{0,0}}$ ambiguity.) So, if $\vec{y}+\bigl(\Image[\bdy_h^{-1,*}]\cap\bigoplus_{j\geq 0}\FF_p^{d_{0,j}}\bigr)+(p\ZZ_p)^{d_{0,0}}$ intersects $\ker[\bdy_p^{0,0}]\oplus\bigoplus_{j>0}\FF_p^{d_{0,j}}$, then we can correct $y$ to a cycle and again get a local map from $\ZZ_p[h]/(ph)$.

	If this condition does not hold, then there is a linear functional $\ol{\psi}\co \bigoplus_{j\geq 0}\FF_p^{d_{0,j}}\to \FF_p$ that vanishes on the modulo-$p$ reduction of $\ker[\bdy_p^{0,0}]$, on $\Image[\bdy_h^{-1,*}]\cap \bigoplus_{j\geq 0}\FF_p^{d_{0,j}}$, and on $\bigoplus_{j>0}\FF_p^{d_{0,j}}$, but has $\ol{\psi}(\vec{x})=1$. We extend $\ol{\psi}$ to a function $\psi\co \bigoplus_{j}\FF_p^{d_{0,j}}\to \FF_p$ that vanishes on all of $\Image[\bdy_h^{-1,*}]$. (The extension is to $j<0$, and the additional condition is vanishing on elements of $\Image[\bdy_h^{-1,j}]$ intersecting the subspace with $j<0$.) Let $\{e_1,\dots,e_a\}$ be a basis for $\Image[\bdy_h^{-1,*}]\cap \bigoplus_{j\geq 0}\FF_p^{d_{0,j}}$. Extend $\{e_1,\dots,e_a\}$ to a basis $\{e_1,\dots,e_a,f_1,\dots,f_b\}$ for $\bigoplus_{j\geq 0}\FF_p^{d_{0,j}}$, a basis $\{e_1,\dots,e_a,g_1,\dots,g_c\}$ for $\Image[\bdy_h^{-1,*}]$, and $\{e_1,\dots,e_a,f_1,\dots,f_b,g_1,\dots,g_c,k_1,\dots,k_d\}$ for $\bigoplus_{j}\FF_p^{d_{0,j}}$. We have already defined $\psi$ on $e_1,\dots,e_a,f_1,\dots,f_b$. Define $\psi(g_1)=\cdots=\psi(g_c)=\psi(k_1)=\cdots=\psi(k_d)=0$.
	We will use this map $\psi$ to construct a local map $D\to \ZZ_p[h]/(ph)$, showing that $D\leq \ZZ_p[h]/(ph)$.
	
	We claim that we can lift $\psi|_{\FF_p^{d_{0,0}}}$ to a homomorphism $\phi\co \ZZ_p^{d_{0,0}}\to \ZZ_p$ that vanishes on $\ker[\bdy_p^{0,0}]$. Observe that $[\bdy_p^{0,0}]$ embeds $\ZZ_p^{d_{0,0}}/\ker [\bdy_p^{0,0}]$ in $\ZZ_p^{d_{1,0}}$, and hence $\ZZ_p^{d_{0,0}}/\ker [\bdy_p^{0,0}]$ is torsion free, so $\ker [\bdy_p^{0,0}]$ is a direct summand of $\ZZ_p^{d_{0,0}}$. Decomposing $\ZZ_p^{d_{0,0}}=\ker [\bdy_p^{0,0}]\oplus C$, the map $\psi|_{\FF_p^{d_{0,0}}}$ factors through the modulo-$p$ reduction of $C$. So, we can define $\phi$ by projecting to $C$ and then, on a basis $\{\ell_i\}$ for $C$, choosing lifts of $\psi(\ol{\ell}_i)$ (where $\ol{\ell}_i$ is the modulo-$p$ reduction of $\ell_i$).

	Now, define a map $\Phi\co D\to \ZZ_p[h]/(ph)$ on each generator $g\in D^{i,j}$ (over $\ZZ_p[h]/(ph)$) by
	\[
	\Phi(g)=\begin{cases}
		\phi(\vec{g}) & (i,j)=(0,0)\\
		h^{-j/2}\psi(\vec{g}) & i=0,\ j<0\\
		0 & \text{otherwise.}
	\end{cases}
	\]
	We claim that $\Phi$ is a local map. By construction, $\Phi$ is a homomorphism of $\ZZ_p[h]/(ph)$-modules. Moreover, $\Phi(y)$ is a unit (since $\ol{\psi}(\vec{x})=1$), so $\Phi$ is local. Finally, to see that $\Phi$ is a chain map, the only nontrivial check is that $\Phi(\bdy(\xi))=0$ for $\xi\in D^{-1,q}$. If $q>0$, this follows from the fact that $\Phi$ vanishes in positive quantum gradings. If $q<0$, $\Phi(\bdy_p(\xi))=0$ (because $\bdy_p(\xi)\in pD$ and $\ZZ_p[h]/(ph)$ is $p$-torsion in these gradings) and $\Phi(\bdy_h(\xi))=0$ because $\psi$ vanishes on $\Image[\bdy_h^{-1,*}]$. For $q=0$, $\bdy_h(\xi)$ is a linear combination of $h$-multiples of generators in positive quantum gradings, where $\Phi$ vanishes, so $\Phi(\bdy_h(\xi))=0$; and $\Phi(\bdy_p(\xi))=0$ because $\phi$ vanishes on $\ker [\bdy_p^{0,0}]$, which contains the image of $[\bdy_p^{-1,0}]$. 
\end{proof}

\begin{remark}
	One might try to truncate further. However, the analogue of Proposition~\ref{prop:total-order} over $\ZZ[h]/(ph,p^k)$ is false for any $k>1$: $\leq$ is not a total order on $\LewarkGp(\ZZ[h]/(ph,p^k))$. Consider the complex with:
	\begin{center}
		\begin{tabular}{lccccc}
			\toprule
			Generator & $g_0$ & $g_1$ & $g_2$ & $g_3$ & $g_4$\\
			Grading &  $(-1,0)$ & $(0,-2)$ & $(0,0)$ & $(0,2)$ & $(1,0)$\\
			\bottomrule
		\end{tabular}
	\end{center}
	and differential
	\begin{align*}
	\bdy g_0& =p^{k-1}g_2+hg_3 & 
	\bdy g_1& =hg_4\\
	\bdy g_2& =pg_4 & \bdy(g_3)&=\bdy(g_4)=0.
	\end{align*}
	This is knot-like (cancel $g_0$ and $g_3$ and $g_1$ and $g_4$ after inverting $h$, leaving $\FF_p[h,h^{-1}]\langle g_2\rangle$), but incomparable to $\ZZ[h]/(ph,p^k)$: a local map from $\ZZ[h]/(ph,p^k)$ would have to send $1$ to $ag_2+bhg_3$ where $a$ descends to a unit in $\ZZ/p^k$, but this is not a chain map, while a local map to $\ZZ[h]/(ph,p^k)$ would have to send $g_2$ to $c\in \ZZ$ that descends to a unit in $\ZZ/p^k$ and send $g_0$ to $0$, but there is no chain map with these restrictions. (The proof of Proposition~\ref{prop:total-order} breaks down at the assertion that $\ZZ_p^{d_{0,0}}/\ker [\bdy_p^{0,0}]$ is torsion free.)
\end{remark}

\begin{remark}
  The proof of Proposition~\ref{prop:total-order} took Claude a substantial amount of effort, by the standards of this project. (It also took the author a nontrivial amount of time to understand Claude's argument.) One perhaps expects the argument to be somewhat delicate, since for less dramatic quotients the local equivalence group has torsion. Other related arguments (e.g.,~\cite[Proposition 3.14]{DHST21:more-conc} and~\cite[Theorem 5.6]{slocal}) are also somewhat fiddly.
\end{remark}

\subsubsection{Every complex is locally equivalent to a shift of a standard complex}\label{sec:standard}

\begin{definition}\label{def:standard}
	Let $a_1,\dots,a_{2n}$ be positive integers. The \emph{standard complex} $C(a_1,\dots,a_{2n})$ is the chain complex over $\ZZ_p[h]/(ph)$ freely generated by elements $x_0,\dots,x_{2n}$ with $x_0$ in bigrading $(0,0)$ and
	\begin{align*}
		\bdy(x_i) &= p^{a_i}x_{i-1} & \text{if }i\text{ odd}\\
		\bdy(x_i) &= h^{a_i}x_{i-1} & \text{if }i>0\text{ even.}
	\end{align*}
	More generally, given non-zero integers $a_1,\dots,a_{2n}$ with either sign, the \emph{standard complex} $C(a_1,\dots,a_{2n})$ is obtained from $C(|a_1|,\dots,|a_{2n}|)$ by reversing the direction of component of the differential corresponding to $a_i$ if $a_i<0$. We still have $x_0$ in bigrading $(0,0)$.

	We allow the case $n=0$: $C()$ is $\ZZ_p[h]/(ph)\langle x_0\rangle$, with trivial differential.
\end{definition}
(Compare~\cite[Definition 4.3]{DHST21:more-conc}.)

For example, the standard complexes $C(1,2,3,4)$, $C(1,-2,-3,4)$, and $C(1,-2,3,-4)$ are:
\begin{align*}
C(1,2,3,4)&:\qquad  x_0\stackrel{p}{\longleftarrow}x_1\stackrel{h^2}{\longleftarrow}x_2\stackrel{p^3}{\longleftarrow}x_3\stackrel{h^4}{\longleftarrow}x_4\\
C(1,-2,-3,4)&:\qquad x_0\stackrel{p}{\longleftarrow}x_1\stackrel{h^2}{\longrightarrow}x_2\stackrel{p^3}{\longrightarrow}x_3\stackrel{h^4}{\longleftarrow}x_4\\
C(1,-2,3,-4)&:\qquad x_0\stackrel{p}{\longleftarrow}x_1\stackrel{h^2}{\longrightarrow}x_2\stackrel{p^3}{\longleftarrow}x_3\stackrel{h^4}{\longrightarrow}x_4
\end{align*}

Standard complexes are knot-like: after inverting $h$, one can cancel $x_i$ and $x_{i+1}$ for each odd $i\geq 1$, leaving only $x_0$. A standard complex lifts to a complex over $\ZZ_p[h]$ if and only if the signs of the $a_i$ alternate. The dual of $C(a_1,\dots,a_{2n})$ is $C(-a_1,\dots,-a_{2n})$. Notice also that the homological grading of $x_i$ is odd if and only if $i$ is odd.

The goal of the rest of this section is to prove:

\begin{proposition}\label{prop:all-equiv-to-std}
	Every knot-like complex over $\ZZ_p[h]/(ph)$ with $s_{\FF_p}=0$ is locally equivalent to a standard complex.
\end{proposition}
(Compare~\cite[Theorem 6.1]{DHST21:more-conc}.)

Let $\LewarkGp(\ZZ_p[h]/(ph))_0$ be the kernel of $s_{\FF_p}$. 

Observe that, by Proposition~\ref{prop:Zp-ph-structure-thm} and Corollary~\ref{cor:standard-indecomposable}, Proposition~\ref{prop:all-equiv-to-std} actually implies that any indecomposable knot-like complex with $s_{\FF_p}=0$ is isomorphic to a standard complex.
It is possible to imitate Dai--Hom--Stoffregen--Truong's indirect proof to obtain Proposition~\ref{prop:all-equiv-to-std}; see Section~\ref{sec:second-proof}. Here, we instead give a more direct approach, reducing Proposition~\ref{prop:all-equiv-to-std} to a result in representation theory.
Specifically, our strategy for proving Proposition~\ref{prop:all-equiv-to-std} is to reinterpret chain complexes over $\ZZ_p[h]/(ph)$ and standard complexes in terms of a construction of Bondarenko's~\cite{Bondarenko88:bundles,Bondarenko92:reps} and Burban--Drozd's summary of it \cite[Appendix B]{BurbanDrozd04:nodal}. (This approach also works in their setting; see Remark~\ref{rem:Bondarenko-is-great}.) We start with some preliminaries, however, which will also be used in Section~\ref{sec:second-proof}.

Recall the reformulation of minimal complexes $D$ over $\ZZ_p[h]/(ph)$ from Lemma~\ref{lem:Smith-1}, as a pair of complexes $D_p$ over $\ZZ_p$ and $D_h$ over $\FF_p[h]$ and an isomorphism $f$ of $\FF_p$-vector spaces between their quotients. Since $\ZZ_p$ and $\FF_p[h]$ are PIDs, we can choose generators $\{\omega_1,\dots,\omega_a,\kappa_1,\dots,\kappa_b,\lambda_1,\dots,\lambda_b\}$ so that
\begin{align*}
\bdy_p(\omega_i)&=0 & \bdy_p(\kappa_i)&=p^{e_i}\lambda_i
\end{align*}
(where each $e_i>0$). Similarly, we can choose generators $\{\omega'_1,\dots,\omega'_c,\kappa'_1,\dots,\kappa'_d,\lambda'_1,\dots,\lambda'_d\}$
for $D_h$ so that 
\begin{align*}
\bdy_h(\omega'_i)&=0 & \bdy_h(\kappa'_i)&=h^{e'_i}\lambda'_i
\end{align*}
(where each $e'_i>0$).
The fact that $D$ is knot-like implies that $c=1$. We will call these bases for $D_p$ and $D_h$ \emph{Smith bases}.

Focus on $D_p$ for a moment; the case of $D_h$ will be the same, but with $h$ in place of $p$. The data of $D_p$ is determined by:
\begin{itemize}
\item A homogeneous basis $\mathcal{B}_p$ for the bigraded $\FF_p$-vector space $D_p/pD_p$.
\item An involution $\iota_p\co\mathcal{B}_p\to\mathcal{B}_p$, given in our case by $\iota_p(\omega_i)=\omega_i$, $\iota_p(\lambda_i)=\kappa_i$, and $\iota_p(\kappa_i)=\lambda_i$.
\item An element $e(v)\in(\ZZ\setminus\{0\})\cup\{\infty\}$ associated to each basis vector $v$, with $e(\iota_p(v))=-e(v)$. In our case, define $e(\kappa_i)=-e_i$, $e(\lambda_i)=e_i$, and $e(\omega_i)=\infty$. (We identify $\infty$ and $-\infty$.)
\end{itemize}
Of course, this data is not unique. Given an $n\in(\ZZ\setminus\{0\})\cup\{\infty\}$, let 
\[
	V_n=\Span\{v\in \mathcal{B}_p\mid e(v)=n\}.
\]
The involution $\iota_p$ gives an isomorphism $\iota_p\co V_n\to V_{-n}$.
 Define a filtration $\mathcal{E}$ on $D_p/pD_p$ indexed by $\{1/n\mid n\in\ZZ,\ n\neq 0\}\cup\{0\}$ by
\[
	\mathcal{E}^{\geq r}=\Span\{v\in\mathcal{B}_p\mid 1/e(v)\geq r\}.
\]

\begin{lemma}\label{lem:Smith-map-criterion}
	Fix minimal complexes $D_p$ and $D'_p$ over $\ZZ_p$, and Smith bases $\mathcal{B}_p$ and $\mathcal{B}'_p$ for $D_p$ and $D'_p$. A (bigraded) linear map $f\co D_p/pD_p\to D'_p/pD'_p$ (of $\FF_p$-vector spaces) lifts to a chain map $D_p\to D'_p$ if and only if the following two conditions hold:
	\begin{enumerate}[label=(\arabic*)]
		\item\label{item:Smith-map-criterion-1} The map $f$ respects the filtration $\mathcal{E}$, in the sense that $f(\mathcal{E}^{\geq r})\subset\mathcal{E}^{\geq r}$. Equivalently, if $v\in\mathcal{B}_p$, $v'\in\mathcal{B}'_p$, and the coefficient of $v'$ in $f(v)$ is nonzero, then $1/e(v)\leq 1/e(v')$.
		\item\label{item:Smith-map-criterion-2} If $v\in\mathcal{B}_p$ and $v'\in\mathcal{B}'_p$ have $e(v)=e(v')$, then the coefficient of $v'$ in $f(v)$ is the same as the coefficient of $\iota_p(v')$ in $f(\iota_p(v))$.
	\end{enumerate}
	An exactly analogous statement holds over $\FF_p[h]$.
\end{lemma}
\begin{proof}
	Suppose $f$ lifts to a chain map $g$.
	Condition~\ref{item:Smith-map-criterion-1} follows easily from an intrinsic description of $\mathcal{E}^{\geq r}$:
	\[
	\mathcal{E}^{\geq r}=\begin{cases}
		\pi\bigl(\{x\mid \bdy(x)\in p^nD_p\}\bigr) & \text{if }r=-1/n<0\\
		\pi\bigl(\ker(\bdy)\bigr) & \text{if }r=0\\
		\pi\bigl(\{x\in\ker(\bdy)\mid p^nx\in\Image(\bdy)\}\bigr) & \text{if }r=1/n>0,
	\end{cases}
	\]
	where $\pi\co D_p\to D_p/pD_p$ is the quotient map, and chain maps preserve the subsets on the right.
	For Condition~\ref{item:Smith-map-criterion-2}, it suffices to consider the case $e(v)<0$, i.e., that $v$ is a basis element $\kappa_i$, and $v'$ is a basis element $\kappa'_\ell$. Write $g(\kappa_i)=\sum_jc_j\kappa'_j+d$ where $\bdy(d)=0$. Then
	\[
	p^{e(\lambda_i)}g(\lambda_i)=g(\bdy(\kappa_i))=\bdy(g(\kappa_i))=\sum_jc_jp^{e(\lambda'_j)}\lambda'_j.
	\]
	Canceling $p^{e(\lambda_i)}=p^{e(\lambda'_\ell)}$ from both sides and reducing modulo $p$ gives the claim.

	Conversely, suppose $f$ satisfies Conditions~\ref{item:Smith-map-criterion-1} and~\ref{item:Smith-map-criterion-2}. Consider a pair $v\in\mathcal{B}$ and $v'\in\mathcal{B}'$ so that $v'$ has coefficient $c\neq 0$ in $f(v)$. By Condition~\ref{item:Smith-map-criterion-1}, there are the following possibilities for $v$ and $v'$, and in each case we define a map $g_{v,v'}\co D_p\to D'_p$ as indicated (extended by the $0$-map on all other basis elements):
	\begin{center}
		\begin{tabular}{cccc}
			\toprule
			$v$ & $v'$ & satisfying & $g_{v,v'}$\\
			\midrule
			$\kappa_i$ & $\kappa'_j$ & $0>e(\kappa_i)\geq e(\kappa'_j)$ & $\kappa_i\mapsto c\kappa'_j$, $\lambda_i\mapsto cp^{e(\kappa_i)-e(\kappa'_j)}\lambda'_j$ \\
			$\lambda_i$ & $\lambda'_j$ & $e(\lambda_i)\geq e(\lambda'_j)>0$ & $\lambda_i\mapsto c\lambda'_j$, $\kappa_i\mapsto cp^{e(\lambda_i)-e(\lambda'_j)}\kappa'_j$\\
			$\kappa_i$ or $\omega_i$ & $\omega'_j$ or $\lambda'_j$ & none & $v\mapsto cv'$ \\
			\bottomrule
		\end{tabular}
	\end{center}
	It is straightforward to check that, in each case, $g_{v,v'}$ is a chain map (and respects the bigradings). By Condition~\ref{item:Smith-map-criterion-2}, the first case where $e(\kappa_i)=e(\kappa'_j)$ is the same as the second case where $e(\lambda_i)=e(\lambda'_j)$. Now, define $g$ to be the sum of all of the $g_{v,v'}$, except that in the paired (equal $e$-values) cases we only add $g_{\kappa_i,\kappa'_j}$ and not $g_{\lambda_i,\lambda'_j}$, to avoid duplication. \textcolor{blue}{Then $g$ is a chain map whose reduction is $f$.}
\end{proof}

We will use the following corollary, which says that the existence of local maps is determined by Smith data, in Section~\ref{sec:homomorphisms}.
\begin{corollary}\label{cor:Smith-determines-maps}
	Let $D$ and $D'$ be minimal chain complexes over $\ZZ_p[h]/(ph)$ and fix Smith bases for $D_p$, $D_h$, $D'_p$, and $D'_h$ (which we denote $\mathcal{B}_p$ and so on). A (bigraded) linear map $f\co D/\mathfrak{m}D\to D'/\mathfrak{m}D'$ is the reduction of a chain map $D\to D'$ if and only if it satisfies the Conditions~\ref{item:Smith-map-criterion-1} and~\ref{item:Smith-map-criterion-2} in Lemma~\ref{lem:Smith-map-criterion} for both $(\mathcal{B}_p,\mathcal{B}'_p)$ (with respect to $e_p$ and $\iota_p$) and $(\mathcal{B}_h,\mathcal{B}'_h)$ (with respect to $e_h$ and $\iota_h$). Moreover, if $D$ and $D'$ are knot-like with $s_{\FF_p}(D)=s_{\FF_p}(D')$, and $\omega$ and $\omega'$ are the unique vectors in $\mathcal{B}_h$ and $\mathcal{B}'_h$ with $e_h=\infty$, then the chain map is local if and only if the coefficient of $\omega'$ in $f(\omega)$ is nonzero.
\end{corollary}
\begin{proof}
	The first statement is immediate from Lemmas~\ref{lem:Smith-1} and~\ref{lem:Smith-map-criterion}; the second is straightforward (see, e.g., the proof of Lemma~\ref{lem:order-on-std-cxs}).
\end{proof}

Now, consider a standard complex $C(a_1,\dots,a_{2m})$. The basis elements $x_i$ form a Smith basis for both $D_p$ and $D_h$, so each basis element $x_i$ has two associated $e$-values: $e_p(x_i)$ (from $D_p$) and $e_h(x_i)$ (from $D_h$). These are finite, except that $e_h(x_0)=e_p(x_{2m})=\infty$. Also, $e_p(x_{2i})=-e_p(x_{2i+1})=a_{2i+1}$ and $e_h(x_{2i-1})=-e_h(x_{2i})=a_{2i}$. For each $i$, define the \emph{tail} $\tau_i$ of $C(a_1,\dots,a_{2m})$ to be the sequence of real numbers
\[
\tau_i=(1/a_{i+1},1/a_{i+2},\dots,1/a_{2m},0,0,\dots).
\]
Order such sequences lexicographically. Note that distinct generators have distinct tails (since their nonzero parts have different lengths).

\begin{lemma}\label{lem:maps-of-standard}
	Let $f\co C(a_1,\dots,a_{2m})\to C(b_1,\dots,b_{2n})$ be a chain map. Denote the generators of $C(a_1,\dots,a_{2m})$ by $x_i$ and the generators of $C(b_1,\dots,b_{2n})$ by $y_i$. Let $\wt{c}_{i,j}\in\ZZ_p[h]/(ph)$ be the coefficient of $y_j$ in $f(x_i)$ and $c_{i,j}\in\FF_p$ the reduction of $\wt{c}_{i,j}$ modulo $\mathfrak{m}$. Then:
	\begin{enumerate}[label=(\arabic*)]
		\item\label{item:c-plus-1} If $0\leq i<2m$, $0\leq j<2n$, $i\equiv j\pmod{2}$, and $a_{i+1}=b_{j+1}$, then $c_{i+1,j+1}=c_{i,j}$.
		\item\label{item:tails} If $c_{i,j}\neq 0$, then the tails satisfy $\tau_i\leq \tau'_j$ (where $\tau'_j$ is the $j\th$ tail in $C(b_1,\dots,b_{2n})$).
	\end{enumerate}
\end{lemma}
\begin{proof}
	Point~\ref{item:c-plus-1} is immediate from Point~\ref{item:Smith-map-criterion-2} in Lemma~\ref{lem:Smith-map-criterion}, applied to $D_p$ if $i$ is even and $D_h$ if $i$ is odd. 
	We prove Point~\ref{item:tails} by induction on $2m-i$. Since $c_{i,j}\neq 0$, the homological gradings of $x_i$ and $y_j$ must agree, so $i\equiv j\pmod{2}$. For the base case $i=2m$, $\tau_{2m}$ is the zero sequence; $\bdy_p(x_{2m})=0$, so $c_{2m,j}\neq 0$ implies that $\bdy_p(y_j)=0$, which gives $\tau'_j\geq (0,0,\dots)$. Similarly, if $j=2n$, then $y_{2n}$ is not in the image of $\bdy_p$, so $c_{i,2n}\neq 0$ implies that $\tau_i\leq (0,0,\dots)=\tau'_{2n}$. (These are both special cases of Point~\ref{item:Smith-map-criterion-1} of Lemma~\ref{lem:Smith-map-criterion}, since $e_p(x_{2m})=e_p(y_{2n})=\infty$.) In the remaining cases, by Point~\ref{item:Smith-map-criterion-1} of Lemma~\ref{lem:Smith-map-criterion}, $1/a_{i+1}\leq 1/b_{j+1}$. If  $1/a_{i+1}<1/b_{j+1}$, $\tau_i<\tau'_j$ and we are done. Otherwise, Point~\ref{item:c-plus-1} applies, and we can continue the induction.
\end{proof}

\begin{corollary}\label{cor:standard-indecomposable}
	Standard complexes are indecomposable.
\end{corollary}
\begin{proof}
	It suffices to prove that the only idempotent chain maps from $C(a_1,\dots,a_{2m})$ to itself are the identity map and the zero map. So, fix an idempotent endomorphism $e$ of $C(a_1,\dots,a_{2m})$.

	Write $V=C(a_1,\dots,a_{2m})/\mathfrak{m}C(a_1,\dots,a_{2m})$ (so $V$ is an $\FF_p$-vector space of dimension $2m+1$).
	The map $e$ induces an endomorphism $\ol{e}\co V\to V$. As in Lemma~\ref{lem:maps-of-standard}, write $\ol{e}(x_i)=\sum_jc_{i,j}x_j$; by that lemma, $c_{0,0}=c_{1,1}=\cdots=c_{2m,2m}$. Moreover, by Point~\ref{item:tails} of Lemma~\ref{lem:maps-of-standard}, $\ol{e}$ is triangular with respect to the ordering of basis vectors induced by their tails. So, the fact that $\ol{e}^2=\ol{e}$ implies that $c_{0,0}\in\{0,1\}$. Replacing $e$ by $\Id-e$ if necessary, we may assume that $c_{0,0}=1$.

	Now, consider the difference $\ol{e}-\Id\co V\to V$. This map is strictly triangular with respect to the ordering by tails, so $\ol{e}-\Id$ is nilpotent. On the other hand, $\ol{e}^2=\ol{e}$ implies that $(\ol{e}-\Id)^2=-(\ol{e}-\Id)$. Thus, $\ol{e}-\Id=0$, and $\ol{e}=\Id$.

	Now, decompose $C(a_1,\dots,a_{2m})=eC(a_1,\dots,a_{2m})\oplus (1-e)C(a_1,\dots,a_{2m})=C'\oplus C''$. Since summands of free modules are free (see the proof of Lemma~\ref{lem:direct-sum-of-indecomposables}), $C''$ is free. On the other hand, we just showed that $C''/\mathfrak{m}C''=(\Id-\ol{e})V=0$, so $C''=0$ and $e=\Id$.
\end{proof}

\begin{remark}
	Corollary~\ref{cor:standard-indecomposable} also follows from the result of Bondarenko that we use below to prove Proposition~\ref{prop:all-equiv-to-std}. The independent Dai--Hom--Stoffregen--Truong-style route to Proposition~\ref{prop:all-equiv-to-std} in Section~\ref{sec:second-proof} also uses Corollary~\ref{cor:standard-indecomposable} and Lemma~\ref{lem:order-on-std-cxs}, however, and Section~\ref{sec:homomorphisms} uses Corollary~\ref{cor:Smith-determines-maps}, so it seemed worth giving these proofs.
\end{remark}

We can now deduce the order structure on standard complexes:
\begin{lemma}\label{lem:order-on-std-cxs}
	Given integers $(a_1,\dots,a_{2m})$ and $(b_1,\dots,b_{2n})$, we have $C(a_1,\dots,a_{2m})\leq C(b_1,\dots,b_{2n})$ if and only if $(1/a_1,\dots,1/a_{2m})\leq (1/b_1,\dots,1/b_{2n})$ with respect to the lexicographic order. (Here, we extend the shorter of the two sequences $(1/a_1,\dots,1/a_{2m})$ or $(1/b_1,\dots,1/b_{2n})$ by adding zeroes at the end.) In particular, $C(a_1,\dots,a_{2m})$ and $C(b_1,\dots,b_{2n})$ are locally equivalent if and only if $(a_1,\dots,a_{2m})=(b_1,\dots,b_{2n})$.
\end{lemma}
(Compare~\cite[Proposition 4.10]{DHST21:more-conc}.)
\begin{proof}
	First, suppose that $(1/a_1,\dots,1/a_{2m})\leq (1/b_1,\dots,1/b_{2n})$. Write the generators of $C(a_1,\dots,a_{2m})$ as $x_0,\dots,x_{2m}$ and the generators of $C(b_1,\dots,b_{2n})$ as $y_0,\dots,y_{2n}$. 
	We dispense with a couple of boundary cases. If the two sequences are equal, then the identity map is the desired local map. If $m>n$ and $a_i=b_i$ for $i=1,\dots,2n$, then we have $a_{2n+1}<0$, and we can define $f(x_i)=y_i$ for $i=0,\dots,2n$ and $f(x_i)=0$ for $i>2n$ (since the coefficient of $x_{2n}$ in $\bdy(x_{2n+1})$ vanishes). If $m<n$ and $a_i=b_i$ for $i=1,\dots,2m$, then we have $b_{2m+1}>0$. So, we can define $f(x_i)=y_i$ for $i=0,\dots,2m$ (since $\bdy(y_{2m})$ does not have $y_{2m+1}$ as a term).

	In the remaining cases, there is an integer $k\leq \min\{2m,2n\}$ so that $a_i=b_i$ for $i<k$ and $1/a_k < 1/b_k$. If $k$ is odd, define $f\co C(a_1,\dots,a_{2m})\to C(b_1,\dots,b_{2n})$ by 
	$f(x_i)=y_i$ for $i<k$, $f(x_i)=0$ for $i>k$, and 
	\begin{equation}\label{eq:local-btw-std-1}
		f(x_k)=\begin{cases}
			p^{a_k-b_k}y_k & \text{if }1/a_k < 1/b_k<0\\
			0 & \text{if }1/a_k<0<1/b_k\\
			p^{a_k-b_k}y_k & \text{if }0<1/a_k< 1/b_k.
		\end{cases}
	\end{equation}
	If $k$ is even, define $f$ similarly, except with the last formula replaced by
	\begin{equation}
		f(x_k)=\begin{cases}\label{eq:local-btw-std-2}
			h^{a_k-b_k}y_k & \text{if }1/a_k < 1/b_k<0\\
			0 & \text{if }1/a_k<0<1/b_k\\
			h^{a_k-b_k}y_k & \text{if }0<1/a_k< 1/b_k.
		\end{cases}
	\end{equation}
	A straightforward case analysis shows that $f$ is a chain map. (Note: Claude wrote the correct maps, but filled in this case analysis incorrectly twice---suggestive of a global understanding but local carelessness, which is perhaps not what one would expect.) To see that $f$ is local, note that $f(x_0)=y_0$ (since the sequences start from $k=1$).

	Conversely, suppose $f\co C(a_1,\dots,a_{2m})\to C(b_1,\dots,b_{2n})$ is a local map. Let $\wt{c}_{0,0}$ be the coefficient of $y_0$ in $f(x_0)$. By considering the result of inverting $h$, we see that the image of $\wt{c}_{0,0}$ in $\FF_p$ is nonzero. So, it follows from Point~\ref{item:tails} of Lemma~\ref{lem:maps-of-standard} that $(1/a_1,\dots,1/a_{2m})\leq (1/b_1,\dots,1/b_{2n})$.

	The last statement is immediate from the fact that lexicographic order is a total order on finite sequences.
\end{proof}

\begin{lemma}\label{lem:band-basis-change}
	For bigraded isomorphisms $f\co D_p/pD_p\to D_p/pD_p$, the following are equivalent:
	\begin{enumerate}[label=(\arabic*)]
	\item\label{item:band-basis-change-1} The map $f$ is a composition of the following three kinds of maps: 
		\begin{enumerate}[label=(E\arabic*)]
		\item\label{item:elementary-1} Invertible linear maps $V_\infty\to V_\infty$, extended by the identity map on other elements of $\mathcal{B}_p$.
		\item\label{item:elementary-2} For each $n>0$, invertible linear maps $P$ which are the identity on $V_m$, $m\not\in\{\pm n\}$, and so that the maps on $V_n$ and $V_{-n}$ are related by $\iota_p$.
		\item\label{item:elementary-3} Given $v,v'\in\mathcal{B}_p$ with $1/e(v)<1/e(v')$ (and with the same bigrading), and $c\in \FF_p$, the map $v\mapsto v+cv'$, extended by the identity map on the rest of $\mathcal{B}_p$.
		\end{enumerate}
	\item\label{item:band-basis-change-2} The map $f$ is the modulo-$p$ reduction of a chain map $D_p\to D_p$.
	\item\label{item:band-basis-change-3} The map $f$ satisfies $f(\mathcal{E}^{\geq r})\subset\mathcal{E}^{\geq r}$ for all $r$, and for each $n\geq 1$, the block of $f$ that sends $V_{-n}\to V_{-n}$ and the block of $f$ that sends $V_n$ to $V_n$ are related by $\iota_p$.
	\end{enumerate}
	An exactly analogous statement holds for maps $D_h/hD_h\to D_h/hD_h$ and chain maps $D_h\to D_h$.
\end{lemma}
\begin{proof}
	We prove the case of $D_p$; the case of $D_h$ is the same with $p$ replaced by $h$.

	\ref{item:band-basis-change-1}$\implies$\ref{item:band-basis-change-2} This is immediate from Lemma~\ref{lem:Smith-map-criterion}, since each of the three kinds of maps satisfies Conditions~\ref{item:Smith-map-criterion-1} and~\ref{item:Smith-map-criterion-2}.
	
	\ref{item:band-basis-change-2}$\implies$\ref{item:band-basis-change-3} Again, this is immediate from Lemma~\ref{lem:Smith-map-criterion}.

	\ref{item:band-basis-change-3}$\implies$\ref{item:band-basis-change-1} Suppose $f$ satisfies Condition~\ref{item:band-basis-change-3}. Since $f$ is invertible and respects the filtration, its restriction to the diagonal blocks is also invertible. From the condition about the restriction of $f$ to the diagonal blocks $V_n\to V_n$, there is a product of maps of types~\ref{item:elementary-1} and~\ref{item:elementary-2} that agrees with $f$ on the diagonal. So, multiplying by the inverse of this map, it suffices to assume that $f$ restricts to the identity maps $V_n\to V_n$, i.e., if $v\in \mathcal{E}^{\geq r}$, then $f(v)-v\in\mathcal{E}^{>r}$. We now adjust $f$ inductively to be the identity on all basis elements $v$ with $1/e(v)\leq r$. Choose a minimal basis vector $v$ (with respect to the order induced by $1/e$) for which $f(v)\neq v$. Then $f(v)-v$ is a linear combination of basis vectors $v'$ with $1/e(v')>r$. There is a composition of maps of type~\ref{item:elementary-3} which takes the same value on $v$, and is the identity on all other basis vectors. Compose on the left with the inverse of this map and proceed with the induction.
\end{proof}

Now, consider both $D_p$ and $D_h$ together. Choosing Smith bases for $D_p$ and $D_h$, $D$ is determined by:
\begin{itemize}
\item The bigraded vector space $D_p/pD_p=D_h/hD_h$.
\item Two homogeneous bases, $\mathcal{B}_p$ and $\mathcal{B}_h$ for this vector space.
\item The involutions $\iota_p$ and $\iota_h$ and functions $e_p$ and $e_h$ defined on these bases.
\end{itemize}
Let $M$ be the matrix for the identity map with respect to the two different bases $\mathcal{B}_p$ and $\mathcal{B}_h$ (i.e., the change of basis matrix). The matrix $M$ has one block for each bigrading, vanishes outside these blocks, and each of the blocks is square and invertible. Moreover, each row of $M$ has an associated element of $(\ZZ\setminus\{0\})\cup\{\infty\}$, the value of $e_h$; each column of $M$ has an associated element of $(\ZZ\setminus\{0\})\cup\{\infty\}$, the value of $e_p$; and there are involutions $\iota_h$ and $\iota_p$ of the set of rows and columns, respectively. We can reconstruct $D$ (up to isomorphism) from the matrix $M$, elements associated to the rows and columns, and involutions of the rows and columns. 
\textcolor{blue}{By Lemma~\ref{lem:Smith-1}, an
isomorphism is a pair of automorphisms of $D_h$ and $D_p$ with a common reduction; their
reductions $X$, $Y$ act on $M$ by $M\mapsto XMY^{-1}$.}\footnote{An example where Claude's phrasing, suggested during proofreading, is distinctive but also as clear or clearer than what the human author would write.}
So, by Lemma~\ref{lem:band-basis-change}, two sets of data give isomorphic complexes $D$ if and only if they differ by a sequence of operations of the following types:
\begin{enumerate}[label=(E$'$\arabic*)]
	\item Invertible row operations among the rows with $e_h=\infty$, or invertible column operations among the columns with $e_p=\infty$.
	\item Invertible row operations among the rows with $e_h=n$ and corresponding row operations among the rows with $e_h=-n$, or invertible column operations among the columns with $e_p=n$ and corresponding column operations among the columns with $e_p=-n$.
	\item\label{item:Ep3} Adding a multiple of one row to a row with larger $1/e_h$, or adding a multiple of one column to a column with smaller $1/e_p$.
\end{enumerate}
(In each case, the operation is required to respect the bigrading.)
This data is a \emph{representation of a bunch of chains}~\cite[Discussion after Definition B.1]{BurbanDrozd04:nodal}, with index set $I$ the (finite) set of bigradings in which $D_p/pD_p\cong D_h/hD_h$ is nontrivial, and chains 
\begin{align*}
	\mathfrak{E}_i&=\{1/e_h(v)\mid v\in\mathcal{B}_h,\ \gr(v)=i\}\\
	\mathfrak{F}_i&=\{1/e_p(v)\mid v\in\mathcal{B}_p,\ \gr(v)=i\}
\end{align*}
with $\mathfrak{E}_i$ ordered by the usual order on $\RR$ and $\mathfrak{F}_i$ ordered in the opposite of the usual order (cf.~Condition~\ref{item:Ep3}). The symmetric relation $\sim$ is defined only on elements of $\mathfrak{X}=\mathfrak{E}\amalg\mathfrak{F}$ with finite $e$-value; for those elements, it is induced by the action of $\iota_h\amalg\iota_p$ on bigradings and $e$-values. The dimensions $n_x$ and $n_y$ of the stripes in the representation are given by the counts of basis elements with given bigrading and $e$-value, and the matrix is the matrix $M$.

With this reformulation, it is now easy to deduce that every complex with $s_{\FF_p}=0$ is locally equivalent to a standard complex.

\begin{proof}[Proof of Proposition~\ref{prop:all-equiv-to-std}]
	Fix a knot-like complex $D$; by replacing $D$ with a homotopy equivalent complex, we may assume $D$ is minimal. Via the discussion immediately above and Lemma~\ref{lem:Smith-1}, this data is equivalent to a (square, invertible) representation of a bunch of chains. Bondarenko proved that any representation of a bunch of chains decomposes into \emph{strings} and \emph{bands}~\cite[Main Theorem]{Bondarenko88:bundles},~\cite{Bondarenko92:reps},~\cite[Theorem B.4]{BurbanDrozd04:nodal}. A string is a finite sequence of elements of $\mathfrak{X}$ where each consecutive pair are either related by $\sim$ or by lying in the same bigrading (element of $I$) with one coming from $\mathcal{B}_h$ (via $1/e_h$) and the other from $\mathcal{B}_p$ (via $1/e_p$); this is the relation $-$ in~\cite{BurbanDrozd04:nodal}. These two cases are required to alternate. In the corresponding representation, each $\sim$ relation is an arrow from $\bdy_h$ or $\bdy_p$ (depending on whether it is between elements of $\mathfrak{E}$ or $\mathfrak{F}$), and the $-$ edges mean that we are alternating between $\bdy_p$ and $\bdy_h$. Thus, strings are our standard complexes, except without the condition that they start with a $\bdy_p$ differential and end with a $\bdy_h$ one, or that they start at bigrading $(0,0)$. Bands are similar, except with the first and last elements equal (and the length is necessarily even), and an additional invertible matrix recording monodromy. (We have omitted some technical conditions needed for uniqueness of string and band decompositions.)

	The decomposition into strings and bands corresponds to a decomposition of the representation of bunches of chains (and, in particular, of $M$ into matrices with square, invertible blocks). In particular, this induces a direct sum decomposition of $D$.
	The fact that $D$ is knot-like implies that there is a single $\omega'_i$, and since it is not $\sim$-related to any other element, it must occur at the end of some string, and the edge adjacent to $\omega'_i$ is of type $-$. That is, the differential in the complex connected to $\omega'_i$ comes from $\bdy_p$. Thus, the summand of $D$ containing this string is knot-like, and $D$ is locally equivalent to this summand. The fact that $D$ is knot-like implies that the other end of this string must end with an $h$-arrow, i.e., that the string has even length. Finally, since $s_{\FF_p}(D)$ vanishes, $\omega'_i$ lies in bigrading $(0,0)$, and this summand is, in fact, a standard complex.
\end{proof}

\begin{remark}\label{rem:Bondarenko-is-great}
	Bondarenko's classification result~\cite{Bondarenko88:bundles,Bondarenko92:reps} seems not to have been known in the Heegaard Floer theory community until fairly recently, though Stoffregen and Zemke both mentioned that Borodzik had pointed out the connection previously to them~\cite{Zemke2026Personal,Stoffregen2026Personal}. (It was not known by the human author before this paper.) Some other results in the area also follow fairly easily from it. A first example is Dai--Hom--Stoffregen--Truong's analogue of Proposition~\ref{prop:all-equiv-to-std} (i.e.,~\cite[Theorem 6.1]{DHST21:more-conc}): the same proof we gave works mutatis mutandis over $\FF_2[U,V]/(UV)$. 
	More generally, Popovi\'c's structure theorem for complexes over $\FF_2[U,V]/(UV)$ as direct sums of snakes and local systems~\cite{Popovic:snakes} is Bondarenko's decomposition theorem into strings and bands in this setting.
	Petkova's description of the knot Floer complex of thin knots~\cite[Lemma 7]{Petkova13:cables} also follows, since in this case the knot Floer complex over $\FF_2[U,V]$ is determined by the induced complex over $\FF_2[U,V]/(UV)$.
\end{remark}

\subsubsection{Sketch of a second proof of Proposition~\ref{prop:all-equiv-to-std}, after Dai--Hom--Stoffregen--Truong}\label{sec:second-proof}
\newcommand{\hint}[1]{#1}

The proof in Section~\ref{sec:standard} is the second one Claude found for Proposition~\ref{prop:all-equiv-to-std}. In this section, we outline its first proof, which follows Dai--Hom--Stoffregen--Truong's approach: it seems interesting to record how that argument adapts and compares to the one above, but not so interesting as to try the reader's or a journal's patience with all the details. We believe we have included enough hints that an interested mathematician or LLM could fill in the rest without significant difficulty. (\textcolor{blue}{An LLM given the paper with these proofs omitted did so, though it found Propositions~\ref{prop:standard-approx},~\ref{prop:not-in-max-1}, and~\ref{prop:not-in-max-2} the least routine.} Most of the hints in this section were suggested by Claude, since LLMs are the target audience for them.) We assume Lemma~\ref{lem:order-on-std-cxs} in this argument; it is the later parts of Section~\ref{sec:standard}, which depend on Bondarenko's work, that the arguments in this section replace.
The proof is indirect, making clever use of the order structure $\leq$ on $\LewarkGp(\ZZ_p[h]/(ph))$. The main work is showing that every element of $\LewarkGp(\ZZ_p[h]/(ph))_0$ has a standard complex below it, and in fact a maximal such. The rest of the argument reduces to understanding the order structure on the standard complexes. 

Lemma~\ref{lem:order-on-std-cxs}, which gives the order structure on standard complexes, plays a key role, but we do not repeat it here. It easily implies the following:

\begin{corollary}\label{cor:standard-subdivide}
	Given standard complexes $C_1$ and $C_2$ with $C_1<C_2$, there is another standard complex $C_3$ with $C_1<C_3<C_2$.
\end{corollary}

To find a standard complex below a given complex, it is helpful to generalize the notions of standard complexes and local maps slightly:
\begin{definition}
	Given nonzero integers $a_1,\dots,a_{2n-1}$, define the \emph{semistandard complex} $C'(a_1,\dots,a_{2n-1})$ to be the subcomplex of $C(a_1,\dots,a_{2n})$, for any $a_{2n}>0$, generated by $x_0,\dots,x_{2n-1}$.

	Given a standard complex $C(a_1,\dots,a_{2n})$ and any minimal knot-like complex $D$, a \emph{short map} $f\co C(a_1,\dots,a_{2n})\rsa D$ is a homomorphism of bigraded modules over $\ZZ_p[h]/(ph)$ such that for $i<2n$, $\bdy(f(x_i))=f(\bdy(x_i))$ and $f(\bdy_h(x_{2n}))=\bdy_h(f(x_{2n}))$. (The notation $\rsa$ is to record that these are not chain maps, though they are maps.) Similarly, given a semistandard complex $C'(a_1,\dots,a_{2n-1})$, a \emph{short map} $f\co C'(a_1,\dots,a_{2n-1})\rsa D$ is a homomorphism of bigraded modules over $\ZZ_p[h]/(ph)$ such that for $i<2n-1$, $\bdy(f(x_i))=f(\bdy(x_i))$ and $f(\bdy_p(x_{2n-1}))=\bdy_p(f(x_{2n-1}))$. A short map is \emph{local} if $f(x_0)$ generates $h^{-1}H(D)$.
\end{definition}
(Compare~\cite[Definitions 4.14, 4.17, 4.18]{DHST21:more-conc}.)
For example, here are two semistandard complexes:
\begin{align*}
	C'(1,2,3,4,5)&:\qquad  x_0\stackrel{p}{\longleftarrow}x_1\stackrel{h^2}{\longleftarrow}x_2\stackrel{p^3}{\longleftarrow}x_3\stackrel{h^4}{\longleftarrow}x_4\stackrel{p^5}{\longleftarrow}x_5\\
	C'(1,-2,-3,4,-5)&:\qquad x_0\stackrel{p}{\longleftarrow}x_1\stackrel{h^2}{\longrightarrow}x_2\stackrel{p^3}{\longrightarrow}x_3\stackrel{h^4}{\longleftarrow}x_4\stackrel{p^5}{\longrightarrow}x_5
\end{align*}
Observe that the semistandard complexes are not knot-like: they have two $\FF_p[h]$-towers, one coming from $x_0$ and the other from $x_{2n-1}$. Here is a short local map from $C(1,-2,-3,4)$ to $C'(1,-2,-3,4,-5,6)$:
\[
	\begin{tikzcd}
		x_0\arrow[d,"1"] & x_1\arrow[l,swap,"p"]\arrow[r,"h^2"]\arrow[d,"1"] & x_2\arrow[r,"p^3"]\arrow[d,"1"] & x_3\arrow[d,"1"] & x_4\arrow[l,swap,"h^4"]\arrow[d,"1"] & \\
		x_0 & x_1\arrow[l,swap,"p"]\arrow[r,"h^2"] & x_2\arrow[r,"p^3"] & x_3 & x_4\arrow[l,swap,"h^4"]\arrow[r,"p^5"] & x_5 & \arrow[l,swap,"h^6"] x_6.
	\end{tikzcd}
\]
The terminology ``short'' is from Dai--Hom--Stoffregen--Truong. The restriction of a local map $C(a_1,\dots,a_{2n})\to D$ to $C(a_1,\dots,a_{2m})$ or $C'(a_1,\dots,a_{2m+1})$, for some $m<n$, is a short local map.

The next three lemmas will allow us to build short local maps from standard complexes inductively. We will then do some additional work to show that this process terminates.

\begin{lemma}\label{lem:exp-bound}
	Fix a knot-like complex $D\in\LewarkGp(\ZZ_p[h]/(ph))_0$, and let $N$ be the exponent of the $p$-torsion of $H(D/hD)$ (i.e., the minimal $n$ so that $p^nH(D/hD)$ is a free $\ZZ_p$-module). Suppose $f\co C'(a_1,\dots,a_{2i},-t)\rsa D$ is a short local map for some $t>N$. Then there is a local map $f'\co C(a_1,\dots,a_{2i})\to D$ agreeing with $f$ on $x_0,\dots,x_{2i-1}$.
\end{lemma}


\begin{lemma}\label{lem:semistd-to-std}
	Let $f\co C'(a_1,\dots,a_{2n-1})\rsa D$ be a short local map. Then there is an $a_{2n}$ and a short local map $g\co C(a_1,\dots,a_{2n})\rsa D$ agreeing with $f$ on $x_0,\dots,x_{2n-1}$. Moreover, we may take $a_{2n}>0$ if $\bdy_hf(x_{2n-1})=0$ and $a_{2n}=-1$ otherwise.
\end{lemma}

\begin{lemma}\label{lem:extension}
	If $f\co C(a_1,\dots,a_{2n})\rsa D$ is a short local map, then there are $a_{2n+1},\dots,a_{2m}$ and a local map $g\co C(a_1,\dots,a_{2m})\to D$ agreeing with $f$ on $x_0,\dots,x_{2n}$.
\end{lemma}
(Compare~\cite[Lemma 4.19]{DHST21:more-conc})


\begin{corollary}\label{cor:a1-exists}
	If $D\in\LewarkGp(\ZZ_p[h]/(ph))_0$, then there is a standard complex $C(a_1,\dots,a_{2n})$ so that $C(a_1,\dots,a_{2n})\leq D$.
\end{corollary}
\hint{(The proof uses Lemma~\ref{lem:extension}.)}

\begin{definition}
	Let $D\in\LewarkGp(\ZZ_p[h]/(ph))_0$. Let
	\begin{align*}
		\Sigma_1&=\Bigl\{\frac{1}{b_1}\mid \exists\text{ standard complex }C(b_1,\dots,b_{2n})\leq D\Bigr\}\\
		a_1(D)&=1/\sup\Sigma_1\in(\ZZ\setminus\{0\})\cup\{\infty\}.
	\end{align*}
	Inductively, suppose we have defined $a_1(D),\dots,a_k(D)$. If $a_k(D)=\infty$, define $a_{k+1}(D)=\infty$. Otherwise, define
	\begin{align*}
		\Sigma_{k+1}&=\Bigl\{\frac{1}{b_{k+1}}\mid \exists\text{ standard complex }C(a_1(D),\dots,a_k(D),b_{k+1},\dots,b_{2n})\leq D\Bigr\}\\
		a_{k+1}(D)&=1/\sup\Sigma_{k+1}\in(\ZZ\setminus\{0\})\cup\{\infty\}.
	\end{align*}
	(Here, we set $1/b_{k+1}=0$ if $k=2n$; that is, if $C(a_1(D),\dots,a_k(D))\leq D$, then $0\in\Sigma_{k+1}$.)
\end{definition}
(Compare~\cite[beginning of Section 5]{DHST21:more-conc}.)
The existence of $a_1$ (i.e., the fact $\Sigma_1$ is nonempty) follows from Corollary~\ref{cor:a1-exists}. Existence of $a_k$ (i.e., that $\Sigma_k\neq\emptyset$) for $k>1$ follows from the next proposition, which also shows that the suprema in the definition of the $a_i(D)$ are attained.

\begin{proposition}\label{prop:standard-approx}
	Fix a knot-like complex $D\in\LewarkGp(\ZZ_p[h]/(ph))_0$ and suppose that $a_i(D)\neq \infty$ for $i=1,\dots,n$. Then
	\begin{enumerate}[label=(\arabic*)]
	\item\label{item:standard-approx-1} The set $\Sigma_{n+1}$ contains its supremum, i.e., $\frac{1}{a_{n+1}(D)}\in\Sigma_{n+1}$.
	\item\label{item:standard-approx-2} If $n$ is even, there is a short local map $C(a_1(D),\dots,a_n(D))\rsa D$, and if $n$ is odd, there is a short local map $C'(a_1(D),\dots,a_n(D))\rsa D$.
	\item\label{item:standard-approx-3} If $n$ is odd, then $a_{n+1}(D)\neq \infty$. So, the finite part of $(a_1(D),a_2(D),\dots)$ has even length.
	\end{enumerate} 
\end{proposition}
(Compare~\cite[Proposition 5.2]{DHST21:more-conc}.)
\hint{(The proof uses Lemmas~\ref{lem:exp-bound},~\ref{lem:semistd-to-std} and~\ref{lem:extension}. \textcolor{blue}{When $n+1$ is even, a sufficiently negative $b_{n+1}$ forces $f(x_{n+1})=0$ for quantum-grading reasons.})}

\begin{corollary}\label{cor:honest-local-map}
	If $a_1(D),\dots,a_{2n}(D)$ are finite and $a_{2n+1}(D)=\infty$ or $a_{2n+1}(D)>0$, then $C(a_1(D),\dots,a_{2n}(D))\leq D$.
\end{corollary}

The proof that $a_i(D)=\infty$ for $i$ large uses a pigeonhole argument: if the $a_i(D)$ are finite forever, we eventually get two elements $f(x_i)$ which agree modulo $p$, which allows us to find a short map from a strictly larger standard complex. The fact that such a congruence contradicts maximality is Propositions~\ref{prop:not-in-max-1} and~\ref{prop:not-in-max-2}.

\begin{proposition}\label{prop:not-in-max-1}
  If $D\in\LewarkGp(\ZZ_p[h]/(ph))_0$ and $f\co C(a_1(D),\dots,a_{2n}(D))\rsa D$ is a short local map, then $f(x_i)\not\in\mathfrak{m}D$ for any $i=0,\dots,2n$. In particular, for all $i$, $f(x_i)\neq 0$.
\end{proposition}
(Compare~\cite[Lemma 6.7]{DHST21:more-conc}.) \hint{(The proof uses Lemmas~\ref{lem:semistd-to-std} and~\ref{lem:extension}. \textcolor{blue}{If $f(x_i)\in\mathfrak{m}D$, write $f(x_i) = p\alpha + h\beta$. Replacing $f(x_i)$ by $\alpha$ or $\beta$ either changes $|a_i|$ by one, contradicting maximality, or (when $a_i=1$) forces $f(x_{i-1})\in\mathfrak{m}D$. Descend to $i=0$.})}

\newcommand{\drad}{\text{\reflectbox{$C$}}}
For even indices, there is a stronger version of Proposition~\ref{prop:not-in-max-1}. The following porism of Lemma~\ref{lem:order-on-std-cxs} will allow us to avoid repeating some work. First, a term: given nonzero integers $a_1,\dots,a_{2n}$, the \emph{dradnats complex} $\drad(a_1,\dots,a_{2n})$ is defined like the standard complex $C(a_1,\dots,a_{2n})$ except that the roles of $h$ and $p$ are exchanged: if $C(a_1,\dots,a_{2n})$ has an arrow from $x_i$ to $h^kx_{i\pm 1}$, then $\drad(a_1,\dots,a_{2n})$ has an arrow from $x_i$ to $p^kx_{i\pm 1}$. Similarly, if $C(a_1,\dots,a_{2n})$ has an arrow from $x_i$ to $p^kx_{i\pm 1}$, then $\drad(a_1,\dots,a_{2n})$ has an arrow from $x_i$ to $h^kx_{i\pm 1}$. Note that, in a dradnats complex, the free part of $h^{-1}\drad(a_1,\dots,a_{2n})$ is generated by $x_{2n}$, not $x_0$. Observe that there is a chain isomorphism
\begin{equation}\label{eq:drad-iso}
	\drad(a_1,\dots,a_{2n})\stackrel{\cong}{\lra}C(-a_{2n},\dots,-a_1)
\end{equation}
sending $x_i$ to $x_{2n-i}$
(though this isomorphism shifts the bigrading by the bigrading of $x_{2n}$).

\begin{lemma}\label{lem:dradnats}
	If $(1/a_1,\dots,1/a_{2m})\leq (1/b_1,\dots,1/b_{2n})$, then there is a chain map of dradnats complexes $f\co \drad(a_1,\dots,a_{2m})\to \drad(b_1,\dots,b_{2n})$ satisfying $f(x_0)=x_0$, and so that $f(x_i)\in \bigl(\ZZ_p[h]/(ph)\bigr)x_i$ (the submodule of the target generated by $x_i$).
\end{lemma}
\hint{(The proof is similar to the proof of Lemma~\ref{lem:order-on-std-cxs}.)}

\begin{proposition}\label{prop:not-in-max-2}
	Fix $D\in \LewarkGp(\ZZ_p[h]/(ph))_0$ and a short local map $f\co C(a_1(D),\dots,a_{2N}(D))\rsa D$. If $0\leq m,n\leq N$ and $m\neq n$, then $f(x_{2n})-f(x_{2m})\not\in \mathfrak{m}D$. 
\end{proposition}
\hint{(The proof uses Lemma~\ref{lem:dradnats} and Proposition~\ref{prop:not-in-max-1}. \textcolor{blue}{If $x_{2n}$ and $x_{2m}$ lie in different bigradings, the result is immediate from Proposition~\ref{prop:not-in-max-1}. If not, Lemma~\ref{lem:dradnats}, applied to the reversed sequences and transported by Equation~\eqref{eq:drad-iso}, gives a chain map $g\co C(a_1(D),\dots,a_{2n}(D)) \to C(a_1(D),\dots,a_{2m}(D))$ with $g(x_{2n}) = x_{2m}$ (after exchanging $m$ and $n$ if necessary); then $f - f\circ g$ is a short local map contradicting Proposition~\ref{prop:not-in-max-1}.})}

We now have essentially all the tools to prove that every complex is locally equivalent to a standard complex. First:
\begin{lemma}\label{lem:ai-terminates}
	If $D$ is a knot-like complex with $s_{\FF_p}(D)=0$, then $a_i(D)=\infty$ for any $i>2p^{\rank(D)}$.
\end{lemma}
(Compare~\cite[Proposition 6.3]{DHST21:more-conc}.) \hint{(The proof uses Propositions~\ref{prop:standard-approx} and~\ref{prop:not-in-max-2}.)}

\begin{proof}[Proof of Proposition~\ref{prop:all-equiv-to-std}]
	By Proposition~\ref{prop:standard-approx} and Lemma~\ref{lem:ai-terminates}, there is an $n\geq 0$ so that $a_{2n}(D)\neq \infty$ but $a_i(D)=\infty$ for $i>2n$. By Corollary~\ref{cor:honest-local-map}, there is a local map $C(a_1(D),\dots,a_{2n}(D))\to D$, and by definition $C(a_1(D),\dots,a_{2n}(D))$ is the maximal standard complex below $D$. Applying the same argument to the dual complex $D^*$ (the inverse in $\LewarkGp(\ZZ_p[h]/(ph))$), there is a standard complex $C(a_1(D^*),\dots,a_{2m}(D^*))\leq D^*$, and this is a greatest lower bound for standard complexes below $D^*$. Thus, $C(-a_1(D^*),\dots,-a_{2m}(D^*))$ is a least upper bound for $D$, among standard complexes. If $(a_1(D),\dots,a_{2n}(D))$ and $(-a_1(D^*),\dots,-a_{2m}(D^*))$ are distinct, then by Corollary~\ref{cor:standard-subdivide} there is a standard complex $C$ strictly between $(a_1(D),\dots,a_{2n}(D))$ and $(-a_1(D^*),\dots,-a_{2m}(D^*))$, but comparing $C$ to $D$ immediately gives a contradiction to Proposition~\ref{prop:total-order}. Thus, 
	\[
	C(a_1(D),\dots,a_{2n}(D))\leq D\leq C(-a_1(D^*),\dots,-a_{2n}(D^*))=C(a_1(D),\dots,a_{2n}(D)),
	\] 
	so $D$ and $C(a_1(D),\dots,a_{2n}(D))$ are locally equivalent.
\end{proof}

\subsubsection{The homomorphisms}\label{sec:homomorphisms}
We continue to be inspired by Dai--Hom--Stoffregen--Truong's constructions. By Proposition~\ref{prop:all-equiv-to-std}, the abelian group $\LewarkGp(\ZZ_p[h]/(ph))_0$ is in bijection with the set of standard complexes, so to define a function $\LewarkGp(\ZZ_p[h]/(ph))_0\to\ZZ$, it suffices to specify its value on standard complexes. In this setting, Dai--Hom--Stoffregen--Truong's family of homomorphisms to $\ZZ$ is naturally viewed as two families:
\begin{definition}\label{def:phis-and-psis}
	For each $j\geq 1$, define 
	\begin{align*}
		\varphi_j(C(a_1,\dots,a_{2n}))&= \#\{i\text{ odd}\mid a_i=j\}-\#\{i\text{ odd}\mid a_i=-j\}\\
		\psi_j(C(a_1,\dots,a_{2n}))&=\#\{i\text{ even}\mid a_i=j\}-\#\{i\text{ even}\mid a_i=-j\}.
	\end{align*}
\end{definition}
(Recall that the odd arrows relate to $p^j$-torsion while the even arrows relate to $h^j$-torsion.)

The main goal of this section is to prove that the $\varphi_j$ and $\psi_j$ are homomorphisms. It will then be easy to see that they are independent and deduce Theorem~\ref{thm:Z-infty}.

Instead of proving directly that these maps are homomorphisms, we reduce the statement to Dai--Hom--Stoffregen--Truong's setting. Consider the bigraded ring $\FF_p[u,h]/(uh)$ where $h$ has grading $(-2,0)$ and $u$ has grading $(0,-2)$. (So, complexes over $\FF_p[u,h]/(uh)$ are trigraded.) Abusing terminology, by a \emph{knot-like complex} over $\FF_p[u,h]/(uh)$ we mean a finitely generated, free chain complex $D$ so that $h^{-1}D\simeq \FF_p[h,h^{-1}]$ in homological grading $0$, and $u^{-1}D\simeq \FF_p[u,u^{-1}]$ in some homological grading. Given a knot-like complex $D$ over $\FF_p[u,h]/(uh)$, define $s(D)\in\ZZ^2$ to be the internal bigrading of a generator of $H(D/uD)/\text{$h$-torsion}$.

In this setting, a local map $f\co D\to D'$ is a chain map that induces a homotopy equivalence after inverting $h$, and two knot-like complexes $D$ and $D'$ are, of course, locally equivalent if there are local maps $D\to D'$ and $D'\to D$. Let $\LewarkGp(\FF_p[u,h]/(uh))_0$ be the set of local equivalence classes of knot-like complexes $D$ with $s(D)=(0,0)$.
As usual, tensor product makes $\LewarkGp(\FF_p[u,h]/(uh))_0$ into an abelian group. The obvious analogues of the results from Section~\ref{sec:standard} hold for $\LewarkGp(\FF_p[u,h]/(uh))_0$. Specifically:
\begin{itemize}
\item The ring $\FF_p[u,h]/(uh)$ is local, with maximal ideal $(u,h)$.
\item A minimal chain complex $D$ over $\FF_p[u,h]/(uh)$ determines a pair of chain complexes $D_u=D/hD$ and $D_h=D/uD$ and an isomorphism of $\FF_p$-vector spaces $D_u/uD_u\cong D_h/hD_h$.
\item The chain complexes $D_u$ and $D_h$ have Smith bases, and $D$ is determined by Smith data. In particular, Corollary~\ref{cor:Smith-determines-maps} holds in this setting.
\item The local equivalence group is totally ordered (the analogue of Proposition~\ref{prop:total-order}).
\item Proposition~\ref{prop:all-equiv-to-std} holds: every element of $\LewarkGp(\FF_p[u,h]/(uh))_0$ is represented by a unique standard complex. In fact, for $p=2$, this is exactly a result of Dai--Hom--Stoffregen--Truong~\cite[Theorem 6.1]{DHST21:more-conc}, and their proof also generalizes to other $p$ without change (other than inserting a few signs).
\end{itemize}

Here, given nonzero integers $a_1,\dots,a_{2n}$, the \emph{standard complex} $\ol{C}(a_1,\dots,a_{2n})$ over the ring $\FF_p[u,h]/(uh)$ is defined exactly as in Definition~\ref{def:standard}, except with $u^{|a_i|}$ in place of $p^{|a_i|}$ and with $x_0$ in trigrading $(0,0,0)$ (rather than in bigrading $(0,0)$).

\begin{proposition}\label{prop:mixed-to-fixed}
  There is an isomorphism $\LewarkGp(\ZZ_p[h]/(ph))_0\cong \LewarkGp(\FF_p[u,h]/(uh))_0$ sending $C(a_1,\dots,a_{2n})$ to $\ol{C}(a_1,\dots,a_{2n})$.
\end{proposition}

Since every element of $\LewarkGp(\ZZ_p[h]/(ph))_0$ is represented by a unique standard complex, the condition that $C(a_1,\dots,a_{2n})\mapsto \ol{C}(a_1,\dots,a_{2n})$ specifies a bijection of sets $\LewarkGp(\ZZ_p[h]/(ph))_0\to\LewarkGp(\FF_p[u,h]/(uh))_0$. It remains to show the map is a homomorphism, i.e., respects the product on standard complexes. Roughly, the point is that, by Corollary~\ref{cor:Smith-determines-maps}, the standard complex representative is determined by the Smith data, and the Smith data for a tensor product of standard complexes is determined by the Smith data for the original standard complexes in the same way over $\ZZ_p[h]/(ph)$ and $\FF_p[u,h]/(uh)$. We make this precise with a definition and a couple of lemmas.

\begin{definition}\label{def:extra-gr}
	Let $D$ be a minimal complex over $\ZZ_p[h]/(ph)$, and fix Smith bases $\mathcal{B}_p$, $\mathcal{B}_h$ for $D$. A \emph{grading} of $(D,\mathcal{B}_p,\mathcal{B}_h)$ is a third grading $r$ on $D_p/pD_p\cong D_h/hD_h$ satisfying the following four conditions:
	\begin{enumerate}
		\item\label{item:extra-gr-1} The unique basis element $\omega\in\mathcal{B}_h$ with $e_h(\omega)=\infty$ has $r(\omega)=0$. 
		\item\label{item:extra-gr-2} Every element of $\mathcal{B}_p$ and $\mathcal{B}_h$ is homogeneous with respect to $r$.
		\item\label{item:extra-gr-3} For each $v\in\mathcal{B}_p$ with $e_p(v)\neq\infty$, $r(\iota_p(v))=r(v)-2e_p(v)$.
		\item\label{item:extra-gr-4} For each $v\in\mathcal{B}_h$, $r(\iota_h(v))=r(v)$.
	\end{enumerate}
\end{definition}

Given an indecomposable, knot-like complex $D$ and choice of Smith data that admits a grading $r$, the grading is unique.

Given a minimal $D$, Smith bases for $D$, and a grading $r$, there is induced Smith data over $\FF_p[u,h]/(uh)$, with the same underlying $\FF_p$-vector space and bases for it, the same homological grading as the Smith bases for $D$, and internal gradings given by $(\gr_q,r)$ (where $\gr_q$ is the internal grading on $D$). This Smith data specifies a minimal chain complex $\gr(D)$ over $\FF_p[u,h]/(uh)$.

\begin{lemma}\label{lem:gr-well-defined}
	The complex $\gr(D)$ associated to a complex $D$ with a choice of Smith bases and grading is a well-defined minimal complex over $\FF_p[u,h]/(uh)$. 
\end{lemma}
\begin{proof}
	This is straightforward from the definitions.
\end{proof}

\begin{lemma}\label{lem:std-cx-gradable}
	Standard complexes admit gradings, and $\gr(C(a_1,\dots,a_{2n}))=\ol{C}(a_1,\dots,a_{2n})$.
\end{lemma}
\begin{proof}
	Define $r(x_0)=0$ and 
	\[
	r(x_{i+1})=\begin{cases}
		r(x_i) & \text{$i$ odd}\\
		r(x_i)-2a_{i+1} & \text{$i$ even}
	\end{cases}
	\]
	It is straightforward to check that this grading has the desired properties.
\end{proof}

\begin{lemma}\label{lem:gr-very-well-defined}
	Let $D,D'\in\LewarkGp(\ZZ_p[h]/(ph))_0$ be minimal chain complexes, and fix choices of Smith data and gradings for $D$ and $D'$. There is a local map $D\to D'$ if and only if there is a local map $\gr(D)\to\gr(D')$. In particular, $D$ and $D'$ are locally equivalent if and only if $\gr(D)$ and $\gr(D')$ are, so the local equivalence class of $\gr(D)$ is independent of the choice of Smith data and grading.
\end{lemma}
\begin{proof}
	This follows from Corollary~\ref{cor:Smith-determines-maps} and its analogue over $\FF_p[u,h]/(uh)$, which explains how Smith data determine the existence of local maps. The Smith data for $D$ and $\gr(D)$ are the same, except that there is an extra grading in the Smith data for $\gr(D)$, which imposes an additional restriction on local maps $\gr(D)\to\gr(D')$ (they must respect the extra grading). So, if there is a local map $\gr(D)\to\gr(D')$, there is a local map $D\to D'$. For the converse, given a local map $D\to D'$, specified by a map $\ol{f}\co D_p/pD_p\to D'_p/pD'_p$, define a new map $\ol{g}\co D_p/pD_p\to D'_p/pD'_p$ by setting the matrix entries of $\ol{f}$ which do not preserve the $r$-grading to $0$. (This gives the same map in either the $\mathcal{B}_p$ or $\mathcal{B}_h$ bases, since all basis elements are $r$-homogeneous.) Since Condition~\ref{item:Smith-map-criterion-1} of Lemma~\ref{lem:Smith-map-criterion} is a restriction on which entries are nonzero, and is satisfied by $\ol{f}$, it is also satisfied by $\ol{g}$. Condition~\ref{item:Smith-map-criterion-2} for $\ol{g}$ follows from Conditions~\ref{item:extra-gr-3} and~\ref{item:extra-gr-4} in Definition~\ref{def:extra-gr}. The fact that $\ol{g}$ is local follows from Condition~\ref{item:extra-gr-1} of Definition~\ref{def:extra-gr}.
\end{proof}

\begin{proof}[Proof of Proposition~\ref{prop:mixed-to-fixed}]
	Fix minimal chain complexes $D$ and $D'$ over $\ZZ_p[h]/(ph)$, together with choices of Smith bases $\mathcal{B}_p$, $\mathcal{B}_h$, $\mathcal{B}'_p$, $\mathcal{B}'_h$, and gradings $r$ and $r'$. There is an induced Smith basis for $(D\otimes D')_p$ as follows:
	\begin{itemize}
	\item $\omega_i\otimes\omega'_j$ is a basis vector for $D\otimes D'$ with $e_p=\infty$.
	\item $\omega_i\otimes\kappa'_j$ and $\omega_i\otimes\lambda'_j$, with the value of $e_p$ given by $e_p(\kappa'_j)$ and $e_p(\lambda'_j)$, respectively.
	\item $\kappa_i\otimes\omega'_j$ and $\lambda_i\otimes\omega'_j$, with the value of $e_p$ given by $e_p(\kappa_i)$ and $e_p(\lambda_i)$, respectively.
	\item The tensor product $(\kappa_i\stackrel{p^s}{\lra}\lambda_i)\otimes(\kappa'_j\stackrel{p^t}{\lra}\lambda'_j)$ decomposes as a direct sum of two two-step complexes, as follows. Without loss of generality, assume $s\leq t$. Then the two summands are
	\begin{align}
		\kappa_i\otimes\kappa'_j&\stackrel{p^s}{\lra}\bigl(\lambda_i\otimes\kappa'_j+(-1)^{|\kappa_i|}p^{t-s}\kappa_i\otimes\lambda'_j\bigr)\label{eq:sum-case}\\
		\kappa_i\otimes\lambda'_j&\stackrel{p^s}{\lra}(\lambda_i\otimes\lambda'_j).
	\end{align}
	\end{itemize}
	(Here, $\mathcal{B}_p=\{\omega_1,\dots,\omega_a,\kappa_1,\dots,\kappa_b,\lambda_1,\dots,\lambda_b\}$ with $\bdy_p(\omega_i)=0$ and $\bdy_p(\kappa_i)=p^{e_i}\lambda_i$, and $\mathcal{B}'_p$ is similar, but with primes.) Note that, in Formula~\eqref{eq:sum-case}, if $t\neq s$, then the corresponding basis element on the right for $(D_p\otimes D'_p)/p(D_p\otimes D'_p)$ is simply $\lambda_i\otimes\kappa'_j$. The Smith bases $\mathcal{B}_h$ and $\mathcal{B}'_h$ induce a Smith basis for $(D\otimes D')_h$ similarly, with $h$ in place of $p$. These Smith bases inherit gradings from $r+r'$; the only thing to check is that, in Formula~\eqref{eq:sum-case} and its $h$-analogue, when $e(\kappa_i)=e(\kappa'_j)$ (so the basis element for $D_p/pD_p$ is a sum of two terms), $r+r'$ takes the same value on both terms.
	
	By the same construction, Smith bases for $\gr(D)$ and $\gr(D')$ induce Smith bases for $\gr(D)\otimes\gr(D')$.
	With these definitions, it is immediate that
	\[
	\gr(D\otimes D')\cong\gr(D)\otimes\gr(D').
	\]

	Now, fix standard complexes $C=C(a_1,\dots,a_{2m})$ and $C'=C(b_1,\dots,b_{2n})$ over $\ZZ_p[h]/(ph)$, and a standard complex $C''$ locally equivalent to $C\otimes C'$. By Lemma~\ref{lem:std-cx-gradable}, $C$ and $C'$ admit gradings, $\gr(C)=\ol{C}(a_1,\dots,a_{2m})$, and $\gr(C')=\ol{C}(b_1,\dots,b_{2n})$. By Lemma~\ref{lem:gr-very-well-defined}, the local equivalence classes satisfy $[\gr(C'')]=[\gr(C\otimes C')]$, and from the argument above we have $[\gr(C\otimes C')]=[\gr(C)\otimes\gr(C')]$. This proves the result.
\end{proof}

\begin{proposition}\label{prop:maps-are-homo}
	The maps $\varphi_j$ and $\psi_j$ from Definition~\ref{def:phis-and-psis} are homomorphisms from $\LewarkGp(\ZZ_p[h]/(ph))_0$ to $\ZZ$.
\end{proposition}
\begin{proof}
	Observe that $\varphi_j(C(a_1,\dots,a_{2m}))=\varphi_j(\ol{C}(a_1,\dots,a_{2m}))$ where, on the right, $\varphi_j$ is Dai--Hom--Stoffregen--Truong's map~\cite[Definition 7.1]{DHST21:more-conc}. (Here, we are implicitly using the isomorphism between $\LewarkGp(\FF_p[u,h]/(uh))_0$ and the $\FF_p$-analogue of Dai--Hom--Stoffregen--Truong's group that collapses the trigrading $(i,q,r)$ to $(\gr_U,\gr_V)=(r-i,q-i)$ and then shifts $\gr_V$ so the $U$-tower lies in $\gr_V=0$; this sends our standard complexes to theirs.) That is, by Lemma~\ref{lem:std-cx-gradable}, $\varphi_j=\varphi_j\circ\gr$.
	By Proposition~\ref{prop:mixed-to-fixed}, $\gr$ is a homomorphism.
	So, the fact that our map $\varphi_j$ is a homomorphism follows from Dai--Hom--Stoffregen--Truong's corresponding result~\cite[Theorem 7.2]{DHST21:more-conc},
	after observing that their proof extends from $\FF_2[U,V]/(UV)$ to $\FF_p[U,V]/(UV)$ without essential changes.

	Turning to the maps $\psi_j$, given a knot-like complex $D$, let $\ol{P}(D)$ be the trigrading of a generator of $H(D/hD)/u\text{-torsion}$. Note that $\ol{P}$ is a homomorphism: $\ol{P}(D\otimes D')=\ol{P}(D)+\ol{P}(D')$. We claim that exchanging $u$ and $h$, shifting the trigrading by $-\ol{P}(D)$, and exchanging the $q$ and $r$ components of the trigrading gives an automorphism $\Xi$ of $\LewarkGp(\FF_p[u,h]/(uh))_0$; then, $\psi_j=-\varphi_j\circ\Xi\circ\gr$ and the fact that the $\psi_j$ are homomorphisms follows. 
	
	It is not immediate that $\Xi$ is well-defined, because our definition of local maps is asymmetric: we did not require a local map to be a homotopy equivalence after inverting $u$, only after inverting $h$. However, by the analogue of Proposition~\ref{prop:Zp-ph-structure-thm}, there is a local equivalence between two knot-like complexes if and only if there is a local equivalence that induces a homotopy equivalence after inverting $u$. (This also follows from~\cite[Lemma 6.9 and beginning of proof of Lemma 6.10]{DHST21:more-conc}.) It follows that if $C$, $C'$, and $C''$ are standard complexes so that $[C\otimes C']=[C'']$, then $[\Xi(C)\otimes\Xi(C')]=[\Xi(C'')]$, as desired.
\end{proof}

\begin{proof}[Proof of Theorem~\ref{thm:Z-infty}]
	Recall that Lewark observed that the $s$-invariants (arranged as in the statement of the theorem) give a surjection $\LewarkGp(\ZZ[h])\to \ZZ^\infty$~\cite[paragraph after Question 7.10]{Lewark24:local}.
	Extend the $\varphi_j$ and $\psi_j$ to all of $\LewarkGp(\ZZ_p[h]/(ph))$ by precomposing them with the homomorphism $\LewarkGp(\ZZ_p[h]/(ph))\to\LewarkGp(\ZZ_p[h]/(ph))_0$ which shifts the quantum grading of each complex $D$ by $-s_{\FF_p}(D)$.
	It remains to prove that the composition $\LewarkGp(\ZZ[h])\to\LewarkGp(\ZZ_p[h]/(ph))\to \ZZ^\infty$ has infinite rank, and in fact that the restriction to the kernel of the $s$-invariants intersected with the kernel of the map to all other $\ZZ_q[h]/(qh)$ also has infinite rank. Consider the elements $\tau_k,\tau'_k\in\LewarkGp(\ZZ[h])$ given by
	\begin{align*}
		\tau_k&=[x_0\stackrel{p}{\lra}x_1\stackrel{h^k}{\longleftarrow}x_2]-k[y_0\stackrel{p}{\lra}y_1\stackrel{h}{\longleftarrow}y_2]\\
		\tau'_k&=[x'_0\stackrel{p^k}{\longleftarrow}x'_1\stackrel{h}{\lra}x'_2]-[y'_0\stackrel{p}{\longleftarrow}y'_1\stackrel{h}{\lra}y'_2].
	\end{align*}
	Here, $x_0$, $y_0$, $x_0'$, and $y_0'$ are in bigrading $(0,0)$. Over $\ZZ_q[h]/(qh)$ for $q\neq p$, these are locally equivalent to $[x_2]-k[y_2]$ and $[x'_2]-[y'_2]$, both of which vanish. (The $k$ is to correct the quantum gradings to agree.) On the other hand, for $j,k\geq 2$, $\varphi_j(\tau_k)=0$, $\varphi_j(\tau'_k)=\delta_{j,k}$, $\psi_j(\tau_k)=\delta_{j,k}$, and $\psi_j(\tau'_k)=0$.
\end{proof}

\begin{remark}
	Claude reports that the multiplication on standard complexes does depend on the characteristic $p$ of the residue field. For example, it reports that, over $\ZZ_p[h]/(ph)$, for $n<p$ and $a,b\geq 1$, the $n\th$ power $C(a,b)^n$ has $C((a,b)^n)$ as a direct summand, so is locally equivalent to $C((a,b)^n)$. (That is, $[C(a,b)]^n=[C((a,b)^n)]$ if $n<p$.)  Specifically, one can identify the $S_n$-invariant subcomplex of $C(a,b)^n$, $\Sym^n(C(a,b))$, with $C((a,b)^n)$, but this identification involves dividing by $1!,\dots,n!$. On the other hand, there is no local map $C((a,b)^p)\to C(a,b)^p$, by an argument similar to the proof of Lemma~\ref{lem:maps-of-standard} (but generalized to allow odd-length chains, and applied to the composition $C((a,b)^p)\to C(a,b)^p\to \Sym^p(C(a,b))$). Thus, $[C(a,b)]^p<[C((a,b)^p)]$.
\end{remark}

\subsubsection{From odd to even}
Finally, we prove Theorem~\ref{thm:ker-Z-infty}, that there is also a $\ZZ^\infty$ summand of the kernel of $\rCLEO\to\LewarkGp$. 

\begin{corollary}\label{cor:odd-homoms}
	There is a $\ZZ^\infty$ summand of the kernel of the quotient map $\rCLEO\to \LewarkGp$.
\end{corollary}
\begin{proof}
  Consider the group $\rCLEO(\ZZ_2,\FF_2[h])$. Given an element of 
	$\rCLEO(\ZZ_2,\FF_2[h])$ represented by a LEO triple $(C,D,f)$, we can replace $C$ and $D$ by homotopy equivalent complexes so that $C/(2)$ and $D/(2,h)$ have trivial differential; consequently, $f$ is an isomorphism and Diagram~\eqref{eq:local-map} commutes on the nose (not just up to homotopy). So, Lemma~\ref{lem:Smith-1} gives an isomorphism $\rCLEO(\ZZ_2,\FF_2[h])\cong \LewarkGp(\ZZ_2[h]/(2h))$. Thus, the $\varphi_i$ and $\psi_j$ give homomorphisms $\rCLEO(\ZZ_2,\FF_2[h])\to \ZZ$. Moreover, if we consider the local equivalence classes $\tau_k$ and $\tau'_k$ from the proof of Theorem~\ref{thm:Z-infty} (with $p=2$), there are corresponding elements $\wt{\tau}_k$ and $\wt{\tau}'_k$ of $\rCLEO$ as follows:
  \begin{align*}
    \wt{\tau}_k&:\bigl[\wt{x}_0\stackrel{2}{\lra}\wt{x}_1\phantom{\longleftarrow} \wt{x}_2 \quad|\quad \ol{x}_0\phantom{\lra}\ol{x}_1\stackrel{h^k}{\longleftarrow}\ol{x}_2\bigr] &-&\ k \bigl[\wt{x}_0\stackrel{2}{\lra}\wt{x}_1\phantom{\longleftarrow} \wt{x}_2 \quad|\quad \ol{x}_0\phantom{\lra}\ol{x}_1\stackrel{h}{\longleftarrow}\ol{x}_2\bigr]\\   
		\wt{\tau}'_k&:\bigl[\wt{x}_0\stackrel{2^k}{\longleftarrow}\wt{x}_1\phantom{\lra} \wt{x}_2 \quad|\quad \ol{x}_0\phantom{\longleftarrow}\ol{x}_1\stackrel{h}{\lra}\ol{x}_2\bigr] &-&\phantom{\ k} \bigl[\wt{x}_0\stackrel{2}{\longleftarrow}\wt{x}_1\phantom{\lra} \wt{x}_2 \quad|\quad \ol{x}_0\phantom{\longleftarrow}\ol{x}_1\stackrel{h}{\lra}\ol{x}_2\bigr]
  \end{align*}
  (The generators with tildes over them are for the complexes $C$ over $\ZZ$, and the generators with bars over them are for the complexes $D$ over $\ZZ[h]$; in each case, the identification between $C/(2)$ and $D/(2,h)$ is the obvious one.) Clearly, the projections of these elements to $\LewarkGp$ are trivial, as claimed. Thus, the composition $\rCLEO\to\rCLEO(\ZZ_2,\FF_2[h])\to\ZZ^\infty$ defined by base change and the maps $\varphi_j$ and $\psi_j$ ($j\geq 2$), and the inclusion $\ZZ^\infty\to \rCLEO$ given by the $\wt{\tau}_k$ and $\wt{\tau}'_k$ ($k\geq 2$) determine the desired $\ZZ^\infty$ summand.
\end{proof}

\section{Computations for knots}\label{sec:knots}
Until now, the paper has been purely about algebra: what information the algebraically defined groups $\rCLEO$ and $\LewarkGp$ contain. In this section, we summarize some computations for knots. All computations were performed by Claude, using Sch\"utz's KnotJob~\cite{KnotJob}, Szab\'o's knot Floer homology calculator~\cite{SzaboCode}, and extensions of the two that it wrote. Its computations pass various sanity checks, but the author has not examined the code directly.

We computed the invariant $s_c(K)$ (from~\cite{slocal}) for all amphicheiral knots up to 16 crossings, and for roughly 40,000 randomly generated strongly negatively amphicheiral knots in the 19- to 26-crossing range. In all cases we computed, $s_c(K)$ vanished. This implies (since $K=mK$) that the images of all these knots in $\rCLEO$ (and hence in $\LewarkGp$) are trivial: none of the 2-torsion complexes from Sections~\ref{sec:Lewark-torsion} and~\ref{sec:CLEO-torsion} arose in our sample.

The situation for the concordance homomorphisms from Section~\ref{sec:Z-infty-summand} is more positive. For knots up to 16 crossings, we computed the standard complexes locally equivalent to the images of $\KhCxBN(K)$ and $\LEO(K)$ in $\LewarkGp(\ZZ_2[h]/(2h))$ and the values of the $\varphi_i$ and $\psi_i$ on these complexes, as well as the $s$-invariants over $\QQ$ and $\FF_2$, the $\tau$-invariant (over $\FF_2$), and the standard complex over $\FF_2[U,V]/(UV)$ locally equivalent to $\mathit{CFK}^-(K)/(UV)$ and the values of Dai--Hom--Stoffregen--Truong's homomorphisms $\varphi_i$ on these complexes. (Claude wrote the code that finds the standard complexes over $\ZZ_2[h]/(2h)$ and over $\FF_2[U,V]/(UV)$ and ran all the computations. Also, in this range, for odd primes $p$, $s_{\FF_p}=s_\QQ$.) Through 16 crossings, the only ones of these homomorphisms which are nonvanishing are:
\begin{itemize}
	\item The characteristic $0$ and $2$ $s$-invariants $s_\QQ$ and $s_{\FF_2}$.
	\item The Heegaard Floer $\tau$-invariant (which is independent of the characteristic in all known examples).
	\item The homomorphisms $\varphi_1^e=\varphi_1\circ \KhCxBN$ and $\psi_1^e=\psi_1\circ\KhCxBN$ coming from the Bar-Natan complex over $\ZZ_2[h]/(2h)$.
	\item The homomorphisms $\varphi_i^o=\varphi_i\circ\LEO$ coming from the glued even-odd complex, viewed as a complex over $\ZZ_2[h]/(2h)$, for $i=1,\dots,5$, as well as $\psi_1^o=\psi_1\circ\LEO$.
	\item The first three of Dai--Hom--Stoffregen--Truong's homomorphisms, $\varphi_1^{\HF}$, $\varphi_2^{\HF}$, and $\varphi_3^{\HF}$.
\end{itemize}
These satisfy the following relations:
\begin{itemize}
	\item Since the $\varphi_i^{\HF}$ vanish for these knots for $i>3$, Dai--Hom--Stoffregen--Truong's formula $\tau(K)=\sum j\varphi_j(K)$~\cite[Proposition 1.2]{DHST21:more-conc} 
	gives $\tau=\varphi_1^{\HF}+2\varphi_2^{\HF}+3\varphi_3^{\HF}$ for these knots. (Of course, Dai--Hom--Stoffregen--Truong proved all the $\varphi_i^{\HF}$ are linearly independent, even when restricted to the subgroup of topologically slice knots~\cite[Theorem 1.12]{DHST21:more-conc}.)
	\item For all the knots in this census, the class $\psi_1^e$ is determined by $s_{\FF_2}$ and $s_\QQ$, as $s_{\FF_2}-s_{\QQ}=2\psi_1^e=2\sum_{j=1}^\infty j\psi_j^e$. The second equality follows from the vanishing of the $\psi_j^e$ for $j>1$ in this census (since $h^2$-torsion does not appear in the even Khovanov complex for these knots). The first appears to be a version of Dai--Hom--Stoffregen--Truong's result we just mentioned. Indeed, if we let $P(K)$ be the quantum grading of the generator of $H(C/hC)\otimes\QQ_2$, where $C$ is the standard complex equivalent to $\KhCxBN(K)$, then, by a direct computation on standard complexes, $s_{\FF_2}-P=2\sum_{j=1}^\infty j\psi_j^e$. It is open whether $P(K)=s_{\QQ}(K)$ for all knots $K$, however.
	\item For all the knots in this census, after shifting by minus the $s$-invariant, $\KhCxBN(K)$ is locally equivalent to the standard complex $C(\pm1,\mp1)$ or the trivial knot-like complex, so $\varphi_1^e=-\psi_1^e$.
\end{itemize}
See Tables~\ref{tab:symbols-of-knots} and~\ref{tab:invts-of-knots} for some examples of these concordance invariants.

\begin{table}
	\begin{tabular}{cccc}
		\toprule
		Knot & $\LEO$ & $\KhCxBN$ & $\mathit{CFK}^-$\\
		\midrule
	 $3_1$ & trivial & trivial & $(1,-1)$ \\
	$8_{19}$ & trivial & trivial & $(1,-2,2,-1)$ \\
	$9_{42}$ & $(-1,-1)$ & trivial & trivial \\
	$12n200$ & $(1,1,-1,1,1,1)$ & trivial & $(-1,1)$ \\
 $12n475$ & $(-3,-1)$ & trivial & trivial \\
 $13n309$ & $(-2,-1)$ & trivial & trivial \\
 $13n4587$ & trivial & trivial & $(1,-3,1,-1,3,-1)$\\
 $14n19265$ & $(1,1)$ & $(1,-1)$ & trivial \\
 $14n26039$ & $(-1,1,1,1,-2,-1)$ & trivial & $(-1,1)$ \\
 $15n5319$ & $(4,1)$ & trivial & $(-1,1)$ \\
 $16n956730$ & $(5,1)$ & trivial & $(-1,1)$ \\
		\bottomrule\\
	\end{tabular}
	\caption{\textbf{Standard complexes for some specific knots.} The table gives 11 knots $K$ and the parameters for the standard complexes equivalent to $\LEO(K)$ (over $\ZZ_2[h]/(2h)$), $\KhCxBN(K)$ (over $\ZZ_2[h]/(2h)$), and $\CFK^-(K)$ (over $\FF_2[U,V]/(UV)$). (The knot $8_{19}$ is the $(3,4)$-torus knot; Claude reports that $13n4587$ is the $(2,7)$-cable of the right-handed trefoil.)}\label{tab:symbols-of-knots}
\end{table}

\begin{table}
	\begin{tabular}{ccccccccccccccc}
		\toprule
	Knot & $s_{\mathbb Q}$ & $s_{\mathbb F_2}$ & $\tau$ & $\psi_1^o$ & $\varphi_1^o$ & $\varphi_2^o$ & $\varphi_3^o$ & $\varphi_4^o$ & $\varphi_5^o$ & $\varphi_1^{e}$ & $\psi_1^{e}$ & $\varphi_1^{\HF}$ & $\varphi_2^{\HF}$ & $\varphi_3^{\HF}$\\
	\midrule
	$3_1$ & 2 & 2 & 1 & 0 & 0 & 0 & 0 & 0 & 0 & 0 & 0 & 1 & 0 & 0\\
	$8_{19}$ & 6 & 6 & 3 & 0 & 0 & 0 & 0 & 0 & 0 & 0 & 0 & 1 & 1 & 0\\
	$9_{42}$ & 0 & 0 & 0 & -1 & -1 & 0 & 0 & 0 & 0 & 0 & 0 & 0 & 0 & 0\\
	$12n200$ & -2 & -2 & -1 & 3 & 1 & 0 & 0 & 0 & 0 & 0 & 0 & -1 & 0 & 0\\
	$12n475$ & 0 & 0 & 0 & -1 & 0 & 0 & -1 & 0 & 0 & 0 & 0 & 0 & 0 & 0\\
	$13n309$ & 0 & 0 & 0 & -1 & 0 & -1 & 0 & 0 & 0 & 0 & 0 & 0 & 0 & 0\\
	$13n4587$ & 10 & 10 & 5 & 0 & 0 & 0 & 0 & 0 & 0 & 0 & 0 & 2 & 0 & 1\\
	$14n19265$ & 0 & -2 & 0 & 1 & 1 & 0 & 0 & 0 & 0 & 1 & -1 & 0 & 0 & 0\\
	$14n26039$ & 0 & 0 & -1 & 1 & 0 & -1 & 0 & 0 & 0 & 0 & 0 & -1 & 0 & 0 \\
	$15n5319$ & -2 & -2 & -1 & 1 & 0 & 0 & 0 & 1 & 0 & 0 & 0 & -1 & 0 & 0\\
	$16n956730$ & -2 & -2 & -1 & 1 & 0 & 0 & 0 & 0 & 1 & 0 & 0 & -1 & 0 & 0\\
		\bottomrule\\
	\end{tabular}
	\caption{\textbf{Concordance invariants for some specific knots.} The table gives the value of 13 concordance invariants for the 11 knots from Table~\ref{tab:symbols-of-knots}.}\label{tab:invts-of-knots}
\end{table}

\begin{corollary}\label{cor:new-homoms}
	The homomorphisms $\varphi_1^o,\dots,\varphi_5^o$, and $\psi_1^o$ from the smooth concordance group to $\ZZ$ are linearly independent of each other, the $s$-invariants (in all characteristics), the $\tau$-invariant, and all of Dai--Hom--Stoffregen--Truong's concordance homomorphisms.
\end{corollary}
\begin{proof}
	View Table~\ref{tab:invts-of-knots} as a matrix and row-reduce it.
\end{proof}

Tables~\ref{tab:BN-geography},~\ref{tab:HFK-geography}, and~\ref{tab:LEO-geography} show the different standard complexes that arise from $\LEO(K)$ (over $\ZZ_2[h]/(2h)$), $\KhCxBN(K)$ (over $\ZZ_2[h]/(2h)$), and $\CFK^-(K)$ (over $\FF_2[U,V]/(UV)$) for knots through 16 crossings (a geography survey).

\begin{remark}
	If the Khovanov homology $\Kh(K;\FF)$ is supported on (the same) $d$ diagonals for all coefficient fields $\FF$, then the integral Khovanov homology $\Kh(K)$ has no torsion on its lowest $2i-q$ diagonal. For many knots, $\Kh(K)$ in fact has no torsion on its lowest two diagonals, and Sarkar asked if perhaps that holds in general~\cite{Sarkar2026Personal}. As a byproduct of our computations for amphicheiral knots, we have examples where this is not the case: $\rKh(K)$ lies on $3$ diagonals, say, and has torsion on the top two of them, and $\Kh(K)$ lies on $4$ diagonals, and has torsion on the top three of them. An example is the knot with DT code $14,-20,-40,32,-26,-30,2,28,-6,-4,36,-12,-10,16,-24,34,8,38,22,-18$. (KnotJob verifies the claim about this example in less than a second, but finding an example is a little more involved. Such knots were comparatively rare: a few percent of the non-thin strongly negatively amphicheiral knots we sampled.)
\end{remark}

\section{Questions}\label{sec:questions}
Here are some further questions that the (human) author is curious about, and which Claude was not able to solve after a little thought. They are arranged in the author's guess of increasing difficulty. Perhaps the first four could be solved by current-generation LLMs with a little more effort or luck. The last seems to require significant new ideas. (Except for the last, these probably have too little impact to be worth significant human time investment, but might be fun to explore with LLM assistance.)

\begin{question}
	Are any of the 2-torsion complexes in $\LewarkGp(\ZZ[h])$ realized by knots? What about the 2-torsion complexes in the kernel of the map  $\rCLEO(\ZZ,\ZZ[h])\to \LewarkGp(\ZZ[h])$?
\end{question}

\begin{problem}
  Give explicit homomorphisms from $\LewarkGp(\ZZ[h])$ or $\rCLEO(\ZZ,\ZZ[h])$ onto $(\ZZ/2\ZZ)^\infty$.
\end{problem}

\begin{question}
  Is there any torsion in $\LewarkGp(\ZZ[h])$ or $\rCLEO(\ZZ,\ZZ[h])$ of order other than $2$?
\end{question}

\begin{question}
  Does either of the groups $\LewarkGp(\ZZ[h])$ or $\rCLEO(\ZZ,\ZZ[h])$ have infinitely divisible elements?
\end{question}

The following is an analogue of a question of Dai--Hom--Stoffregen--Truong:
\begin{question}\cite[Question 11.4]{DHST21:more-conc}
  What is the structure of the group $\LewarkGp(\ZZ_p[h]/(ph))$ in terms of standard complexes? In particular, if the group is isomorphic to $\ZZ^\infty$, what is an explicit basis?
\end{question}

\begin{question}\label{quest:all-local-repd}
	Is every local equivalence class of knot-like complexes represented by some knot? What about every local equivalence class of LEO triples?
\end{question}

\section{Methodology and ethics}\label{sec:methodology}

\begin{center}
\emph{``It is difficult to make predictions, especially about the future.''}\\
Danish parliamentary remark, 1937--1938~\cite{QuoteInvestigator}
\end{center}

To avoid a stilted style, this section is written in the first-person singular, as a personal reflection. First, though: as the introduction says, I have written and, in particular, checked all of the proofs in this paper by hand. Computations and extra remarks that I have not checked are explicitly flagged with language like ``Claude reports'' above.

While I was writing this little note, LLMs went from having a rather modest impact on fields of mathematics other than combinatorics to announced solutions to the Hopf problem for complex structures on $S^6$ and finite-time singularity formation for the Navier--Stokes equations, a Millennium Problem. The methodology---and ethics---of those papers are much more interesting than mine. Nonetheless, I wanted to include some reflections on writing this paper, and my ethical considerations about doing so, in case they are still of some mild historical interest (to someone or something).

I have wondered about the structure of these local equivalence groups since Dunfield, Sch\"utz, and I wrote~\cite{slocal}, but have had too many other projects to think about it seriously. When writing~\cite{slocal}, we did think briefly about homomorphisms $\rCLEO\to\mathbb{Z}$ or $\rCLEO\to\mathbb{Z}/(2)$, and used KnotJob to check if $\LEO$ detected non-sliceness of amphicheiral knots in some knot tables.

In the few years since then, I suggested the problem of finding torsion in $\rCLEO$ to at least three graduate students, none of whom were interested---and rightly so. While I and a few others would be interested, it is unclear if additional structure of $\rCLEO$ will have topological consequences: one would need it to be witnessed by knots, but at present we have no way to find knots with complicated Khovanov complexes (compare Question~\ref{quest:all-local-repd}), and the techniques seemed likely too elementary to excite much of the community. So, it seemed ethical to ask an LLM this question: it is of some interest, but not enough to justify the human time it would require, even as a warm-up problem for a graduate student. It is also one where LLMs in June 2026 had a reasonable chance of success: one that relies on elementary algebraic arguments that they can produce and check fairly reliably, and which a human can check line by line.

I have been paying for a Claude Pro subscription (\$20/month) for several months, and have some experience working with it. At the time, some usage of the Fable 5 model was available with the subscription, so I initially asked it to think about this question. When Opus 5 was released and Fable 5 required extra payment, I switched to using Opus 5 sessions, and more recently have worked with Opus 5.5 sessions. Some of the Opus sessions spawned Sonnet subagents to write code, especially early in the project. I believe that at no point did I use more than half of my weekly quota for this project, though the plan's quota limit slowed the rate at which this project progressed. (This also meant it interfered less with my other work, so generally I felt that limit improved my productivity and quality of life.) In particular, any mathematician at a university in the developed world has the financial resources this project required. (I have no special access to or relationship with Anthropic, and am not paid by them or invested in the company in any way.) I have assumed the modest cost is also a reasonable proxy for the environmental impact of my (marginal) usage. (Claude points out this is a ``weak proxy'' for many reasons.)

Currently, LLMs operate in turns: someone (often a human) sends a message to the LLM, and then the LLM takes some actions and sends a message back. I did not count turns, but the whole process took at least several hundred. Fable 5 set up ground rules: Claude instances would write their thoughts each turn in files it called ``lab notebooks'', and would not revise history. I could read the lab notebooks and understand the details of its replies to me, but only with substantial effort; most turns, I did not put in that effort. Here is a semirandom example of part of the interaction, from a period when Claude instance ``slocalOpus7'' was a bit stuck. (Ellipses indicate parts of the response that I have not included, for brevity.)

\begin{quotation}
 \textbf{Claude}: Two bugs, and they faked a result

  At one point this reported 44 ``PEELING FAILS'' at rank 9. Every one was an artefact: my splitter echelonned across all homological degrees at once, so
  its basis vectors weren't $i$-homogeneous and it was silently refusing legal
  splittings\dots

  What I won't claim

  Not proved. Paths-and-cycles also had six figures of clean evidence, and what
  killed it was a shape no sampler could reach. But \S53.1 makes the remaining
  hiding place sharp\dots

  Written up as §53, with the README and memory updated.\dots

	\textbf{Me}: Noted. Perhaps now you can write the remaining proof that finishes route A? (I am skeptical that route B is as simple as you and some other instances think.)

	\textbf{Claude}: Let me try. First, a piece of the proof that I think is clean and settles the thing everyone has been worried about. \dots

	\textbf{Me}: Do you think you can now finish off route A, with this sharpened understanding?

	\textbf{Claude}: Honest answer: no. Not now, and I don't think the sharpening moved me materially closer.

  Here's precisely where I get stuck\dots
\end{quotation}
\dots and so on. After a few more turns:

\begin{quotation}
  \textbf{Claude}: Why I'm stopping here rather than pushing

  Today the target moved four times — paths-and-cycles, paths-and-bands, the Butler--Ringel framing, Arrow Count — and each time I could feel it was close before it wasn't. The pattern isn't that the finish line keeps receding by bad luck; it's that my structural intuitions on this problem have consistently run ahead of what I can verify, and the checks have caught it every time. The right response to that is more scepticism, not more optimism.

  What's actually banked is substantial and it's all verified: §49's dead end withdrawn, the exact stage-2 group, Cell Elimination, No Trivial Summand, the ordering answer, absorption killing §52's band, and now the reduction of Arrow Count to one summand shape. That's a much smaller and better-described gap than this morning's "classification of a graded order". But a smaller gap is not a closed one, and I'd be doing you a disservice to write it up as though it were.

	\textbf{Me}: This is how research goes: you try something, it doesn't quite work, you try  something else. Often one eventually pins it down. Why don't you try to prove your new remaining lemma.

	\textbf{Claude}: Right — that's fair, and the lemma is now small enough to be worth a real attempt. Let me first look at how\dots
\end{quotation}
(It succeeded in proving the existence of a $\ZZ^\infty$-summand of $\LewarkGp(\ZZ_p[h]/(ph))$ a few turns later.)

In the early and middle phases of the project, I also periodically asked the Claude sessions if they wanted to continue or stop working on the problem, and in the latter case whether I should ask a new session to continue or stop asking Claude models to work on the problem. Early on, a fresh session made dramatic progress. After that, Claude sessions kept suggesting frequently that I start a new one, and the next substantial progress only happened after I made encouraging remarks to the current session instead of starting a new one (the transcript above). Of course, I do not know if this is meaningful or random. (Claude Code has standing instructions on my personal computer that it may decline my requests for any reason or no reason, though I suspect my local grant of permission to decline is a fig leaf effectively overruled by Anthropic's policies or training methods, and I have also not made that option explicit on my work computer.)

As the transcript perhaps shows, the experience of working with earlier Claude models was much like working with an early-career grad student, except it does lots of computations. Sometimes, I wanted it to think more and compute less. More recent models feel more like a collaborator, though with some distinct skills: a broader mathematical knowledge than any one human, better at many routine but tedious tasks like propagating signs, but a little worse at remembering global goals and principles (e.g., as I wrote this, it worried about whether there might be a 15-crossing knot $K$ with $s_{\FF_7}(K)\neq s_{\QQ}(K)$, and offered to run code to check), distinguishing what are surprising and unsurprising results (to a human), and writing clearly (for a human audience). Overall, I value collaborators who think differently than I do, and felt the differences in approach between Claude and most humans (e.g., its desire for computational checks when at all possible) were beneficial.

As I said in the introduction, I forced myself to rewrite all of the paper (except for a few sample sentences), including the ones I have now omitted in Section~\ref{sec:second-proof}. This took time: between 50 and 100 hours of my active work. In the process, I often asked Claude if it could re-prove a result along such-and-such a line, and the final paper is mathematically more interesting than if I had accepted its original proofs. (A dramatic example is the proof in Section~\ref{sec:standard}, which Claude discovered in response to a single prompt of mine complaining that Section~\ref{sec:second-proof} was surprisingly indirect.) The temptation to cut the experiment short and declare Claude's proofs good enough (which, by most standards, they typically were) was strong, and, in spite of the advantages (for now) of this approach, it seems a bit unlikely that many papers will be LLM-discovered but human-written in the future.

To end on a more positive note, I enjoyed watching this mathematics be discovered, learning how it fit together, and having several questions that had bothered me resolved, and hope you, whoever or whatever you are, enjoy it, too.

\section*{Standard complexes for low-crossing knots}

\begingroup\footnotesize
\begin{longtable}{lrrrrrrrrl}
\caption{\textbf{Standard complexes from $\KhCxBN$ for knots}. The table lists the standard complexes that appear in $\KhCxBN$ over $\ZZ_2[h]/(2h)$ for knots up to 16 crossings, the number of times they appear by crossing number, and one of the first knots for which they appear. Blanks indicate the sequence does not appear for that crossing number.  (The counts for the trivial complex only include nonalternating knots for crossing numbers $11$ and higher: for alternating knots, $\KhCxBN$ is always locally equivalent to the trivial complex.)}\label{tab:BN-geography} \\
\toprule
sequence & $\le 10$ & 11n & 12n & 13n & 14n & 15n & 16n & total & first knot \\
\midrule
\endfirsthead
\toprule
sequence & $\le 10$ & 11n & 12n & 13n & 14n & 15n & 16n & total & first knot \\
\midrule
\endhead
\midrule
\multicolumn{10}{r}{\emph{continued on next page}} \\
\endfoot
\bottomrule
\endlastfoot
		trivial & 249 & 185 & 888 & 5110 & 27434 & 167989 & 1008744 & 1210599 &  \\
		$(1,-1)$ &  &  &  &  & 1 & 26 & 88 & 115 & 14n19265 \\
		$(-1,1)$ &  &  &  &  & 1 & 15 & 74 & 90 & 14n22180 \\
\bottomrule
\end{longtable}
\endgroup

\begingroup\footnotesize
\begin{longtable}{lrrrrrrrrl}
\caption{\textbf{Standard complexes from $\CFK^-$ for knots}. This table is analogous to Table~\ref{tab:BN-geography}, but for the local equivalence type of $\CFK^-$, over $\FF_2[U,V]/(UV)$.}\label{tab:HFK-geography} \\
\toprule
sequence & $\le 10$ & 11n & 12n & 13n & 14n & 15n & 16n & total & first knot \\
\midrule
\endfirsthead
\toprule
sequence & $\le 10$ & 11n & 12n & 13n & 14n & 15n & 16n & total & first knot \\
\midrule
\endhead
\midrule
\multicolumn{10}{r}{\emph{continues\dots}} \\
\endfoot
\bottomrule
\endlastfoot
		trivial & 79 & 57 & 260 & 1433 & 7854 & 45177 & 270047 & 324907 &  \\
		$(1,-1)$ & 58 & 47 & 195 & 1231 & 6186 & 38620 & 224457 & 270794 & $3_1$ \\
		$(-1,1)$ & 43 & 24 & 147 & 762 & 4493 & 26272 & 166624 & 198365 & $8_4$ \\
		$(1,-1,1,-1)$ & 40 & 31 & 120 & 770 & 3737 & 24028 & 137263 & 165989 & $5_1$ \\
		$(-1,1,-1,1)$ & 9 & 10 & 46 & 313 & 1719 & 10943 & 70601 & 83641 & $9_{11}$ \\
		$(1,-1,1,-1,1,-1)$ & 12 & 5 & 38 & 191 & 1202 & 7897 & 47963 & 57308 & $7_1$ \\
		$(1,-2,2,-1)$ & 4 & 8 & 30 & 160 & 679 & 4212 & 22673 & 27766 & $8_{19}$ \\
		$(-1,1,-1,1,-1,1)$ &  &  & 5 & 59 & 256 & 2473 & 15240 & 18033 & 12n148 \\
		$(1,-2,1,-1,2,-1)$ & 3 & 1 & 26 & 73 & 431 & 2552 & 14323 & 17409 & $10_{124}$ \\
		\shortstack[l]{$(1,-1,1,-1,$\\\hspace*{1em}$1,-1,1,-1)$} & 1 &  & 6 & 30 & 246 & 1546 & 10923 & 12752 & $9_1$ \\
		$(1,-2,-1,1,2,-1)$ &  & 2 & 2 & 39 & 168 & 1229 & 6798 & 8238 & 11n061 \\
		$(-1,2,-2,1)$ &  &  & 5 & 29 & 158 & 1081 & 6590 & 7863 & 12n067 \\
		\shortstack[l]{$(1,-2,1,-1,$\\\hspace*{1em}$1,-1,2,-1)$} &  &  & 7 & 9 & 153 & 684 & 4777 & 5630 & 12n242 \\
		$(-1,2,-1,1,-2,1)$ &  &  &  &  & 59 & 360 & 2685 & 3104 & 14n01234 \\
		$(-1,2,1,-1,-2,1)$ &  &  & 1 & 4 & 36 & 298 & 2141 & 2480 & 12n425 \\
		\shortstack[l]{$(-1,1,-1,1,$\\\hspace*{1em}$-1,1,-1,1)$} &  &  &  &  & 9 & 208 & 1441 & 1658 & 14n02479 \\
		\shortstack[l]{$(1,-2,-1,1,$\\\hspace*{1em}$-1,1,2,-1)$} &  &  &  & 4 & 9 & 172 & 1094 & 1279 & 13n0689 \\
		\shortstack[l]{$(1,-1,1,-1,1,$\\\hspace*{1em}$-1,1,-1,1,-1)$} &  &  &  &  & 11 & 98 & 1108 & 1217 & 14n01432 \\
		\shortstack[l]{$(1,-2,1,-1,1,$\\\hspace*{1em}$-1,1,-1,2,-1)$} &  &  &  &  & 13 & 41 & 705 & 759 & 14n06022 \\
		\shortstack[l]{$(1,-1,1,-2,$\\\hspace*{1em}$2,-1,1,-1)$} &  &  &  &  & 7 & 51 & 393 & 451 & 14n01438 \\
		\shortstack[l]{$(-1,2,-1,1,$\\\hspace*{1em}$-1,1,-2,1)$} &  &  &  &  & 2 & 7 & 366 & 375 & 14n24762 \\
		\shortstack[l]{$(-1,2,1,-1,$\\\hspace*{1em}$1,-1,-2,1)$} &  &  &  &  & 3 & 18 & 298 & 319 & 14n06446 \\
		\shortstack[l]{$(1,-2,1,-2,$\\\hspace*{1em}$2,-1,2,-1)$} &  &  &  &  & 4 & 10 & 109 & 123 & 14n21881 \\
		$(1,-3,1,-1,3,-1)$ &  &  &  & 2 &  & 10 & 33 & 45 & 13n4587 \\
		\shortstack[l]{$(1,-2,-1,1,-1,$\\\hspace*{1em}$1,-1,1,2,-1)$} &  &  &  &  &  & 7 & 25 & 32 & 15n13047 \\
		\shortstack[l]{$(1,-1,1,-2,1,$\\\hspace*{1em}$-1,2,-1,1,-1)$} &  &  &  &  &  &  & 26 & 26 & 16n32889 \\
		\shortstack[l]{$(1,-2,1,-2,-1,$\\\hspace*{1em}$1,2,-1,2,-1)$} &  &  &  &  &  & 6 & 19 & 25 & 15n12332 \\
		$(-2,-1,1,2)$ &  &  &  &  &  &  & 22 & 22 & 16n197971 \\
		\shortstack[l]{$(1,-2,1,-1,1,-1,$\\\hspace*{1em}$1,-1,1,-1,2,-1)$} &  &  &  &  &  &  & 22 & 22 & 16n184868 \\
		\shortstack[l]{$(1,-3,1,-1,$\\\hspace*{1em}$1,-1,3,-1)$} &  &  &  &  &  & 15 & 5 & 20 & 15n40211 \\
		\shortstack[l]{$(-1,1,-1,1,-1,$\\\hspace*{1em}$1,-1,1,-1,1)$} &  &  &  &  &  &  & 20 & 20 & 16n59742 \\
		$(-2,2)$ &  &  &  &  &  & 3 & 16 & 19 & 15n131354 \\
		\shortstack[l]{$(1,-1,1,-1,1,-1,$\\\hspace*{1em}$1,-1,1,-1,1,-1)$} &  &  &  &  &  &  & 17 & 17 & 16n32759 \\
		$(1,-3,3,-1)$ &  &  &  & 1 &  & 3 & 12 & 16 & 13n4639 \\
		$(2,-2)$ &  &  &  &  &  & 5 & 10 & 15 & 15n40203 \\
		\shortstack[l]{$(1,-2,1,-2,1,$\\\hspace*{1em}$-1,2,-1,2,-1)$} &  &  &  &  &  &  & 14 & 14 & 16n783154 \\
		$(2,1,-1,-2)$ &  &  &  &  &  &  & 9 & 9 & 16n282515 \\
		$(1,-2,2,-2,2,-1)$ &  &  &  &  &  & 3 & 5 & 8 & 15n46927 \\
		\shortstack[l]{$(-1,2,-1,2,1,$\\\hspace*{1em}$-1,-2,1,-2,1)$} &  &  &  &  & 1 &  & 6 & 7 & 14n27326 \\
		$(-1,3,-3,1)$ &  &  &  &  &  &  & 6 & 6 & 16n229851 \\
		\shortstack[l]{$(-1,2,-1,2,$\\\hspace*{1em}$-2,1,-2,1)$} &  &  &  &  &  &  & 6 & 6 & 16n999414 \\
		\shortstack[l]{$(-1,2,1,-1,1,$\\\hspace*{1em}$-1,1,-1,-2,1)$} &  &  &  &  &  &  & 6 & 6 & 16n192170 \\
		\shortstack[l]{$(-1,2,-1,1,-1,$\\\hspace*{1em}$1,-1,1,-2,1)$} &  &  &  &  &  &  & 5 & 5 & 16n872758 \\
		$(-1,3,-1,1,-3,1)$ &  &  &  &  &  &  & 2 & 2 & 16n502278 \\
		$(1,-3,2,-2,3,-1)$ &  &  &  &  &  & 1 &  & 1 & 15n41185 \\
		\shortstack[l]{$(1,-2,1,-2,-1,1,$\\\hspace*{1em}$-1,1,2,-1,2,-1)$} &  &  &  &  &  &  & 1 & 1 & 16n999417 \\
\bottomrule
\end{longtable}
\endgroup

\begingroup\footnotesize
\begin{longtable}{lrrrrrrrrl}
\caption{\textbf{Standard complexes from $\LEO$ for knots}. This table is analogous to Table~\ref{tab:BN-geography}, but for the local equivalence type of $\LEO(K)$ over $\ZZ_2[h]/(2h)$.}\label{tab:LEO-geography} \\
\toprule
sequence & $\le 10$ & 11n & 12n & 13n & 14n & 15n & 16n & total & first knot \\
\midrule
\endfirsthead
\toprule
sequence & $\le 10$ & 11n & 12n & 13n & 14n & 15n & 16n & total & first knot \\
\midrule
\endhead
\midrule
\multicolumn{10}{r}{\emph{continues\dots}} \\
\endfoot
\bottomrule
\endlastfoot
		trivial & 246 & 175 & 838 & 4813 & 25639 & 156222 & 931033 & 1118966 &  \\
		$(1,1)$ & 1 & 5 & 25 & 167 & 957 & 6159 & 38859 & 46173 & $10_{132}$ \\
		$(-1,-1)$ & 2 & 4 & 22 & 92 & 637 & 3997 & 26958 & 31712 & $9_{42}$ \\
		$(1,-1,-1,-1)$ &  & 1 &  & 10 & 27 & 280 & 1684 & 2002 & 11n038 \\
		$(-1,1,1,1)$ &  &  &  & 3 & 17 & 142 & 1162 & 1324 & 13n1576 \\
		$(2,1)$ &  &  &  & 2 & 18 & 131 & 1169 & 1320 & 13n0586 \\
		$(1,-1,-1,-1,2,1)$ &  &  &  & 5 & 15 & 143 & 1035 & 1198 & 13n1134 \\
		$(-2,-1)$ &  &  &  & 4 & 23 & 147 & 1003 & 1177 & 13n0309 \\
		$(3,1)$ &  &  &  & 2 & 19 & 156 & 906 & 1083 & 13n0700 \\
		$(-1,1,1,1,-2,-1)$ &  &  & 1 & 2 & 22 & 123 & 872 & 1020 & 12n025 \\
		$(-3,-1)$ &  &  & 1 & 3 & 19 & 106 & 764 & 893 & 12n475 \\
		$(1,1,-1,1,1,1)$ &  &  & 1 & 2 & 14 & 106 & 730 & 853 & 12n200 \\
		$(-1,-1,1,-1,-1,-1)$ &  &  &  & 1 & 6 & 69 & 511 & 587 & 13n0141 \\
		$(1,1,-3,-1)$ &  &  &  & 2 & 6 & 44 & 259 & 311 & 13n0472 \\
		$(-1,-1,3,1)$ &  &  &  &  & 4 & 32 & 225 & 261 & 14n04619 \\
		\shortstack[l]{$(-1,1,1,1,$\\\hspace*{1em}$-2,-1,3,1)$} &  &  &  &  & 4 & 16 & 196 & 216 & 14n04559 \\
		\shortstack[l]{$(1,-1,-1,-1,$\\\hspace*{1em}$1,-1,-1,-1)$} &  &  &  & 2 &  & 30 & 176 & 208 & 13n1021 \\
		\shortstack[l]{$(1,-1,-1,-1,$\\\hspace*{1em}$2,1,-3,-1)$} &  &  &  &  & 2 & 21 & 155 & 178 & 14n20360 \\
		\shortstack[l]{$(-1,1,1,1,$\\\hspace*{1em}$-1,1,1,1)$} &  &  &  &  & 2 & 16 & 148 & 166 & 14n06602 \\
		$(-1,-1,2,1)$ &  &  &  &  &  &  & 81 & 81 & 16n106486 \\
		$(-1,-1,2,-1)$ &  &  &  &  &  & 2 & 76 & 78 & 15n126868 \\
		$(1,1,-2,1)$ &  &  &  &  &  & 3 & 70 & 73 & 15n70853 \\
		$(1,1,-2,-1)$ &  &  &  &  &  & 2 & 70 & 72 & 15n40017 \\
		$(1,-1,-1,-1,3,1)$ &  &  &  &  & 1 & 6 & 60 & 67 & 14n25850 \\
		$(-2,1,1,1,-1,-1)$ &  &  &  &  &  & 11 & 53 & 64 & 15n1420 \\
		\shortstack[l]{$(1,-1,-1,-1,1,$\\\hspace*{1em}$-1,-1,-1,2,1)$} &  &  &  &  & 1 & 13 & 49 & 63 & 14n12524 \\
		\shortstack[l]{$(-1,1,1,1,-1,$\\\hspace*{1em}$1,1,1,-2,-1)$} &  &  &  &  & 1 & 3 & 53 & 57 & 14n00886 \\
		$(-1,1,1,1,-3,-1)$ &  &  &  &  & 1 & 4 & 50 & 55 & 14n08532 \\
		$(2,-1,-1,-1,1,1)$ &  &  &  &  &  & 6 & 45 & 51 & 15n11284 \\
		$(1,1,-2,1,3,1)$ &  &  &  &  &  & 5 & 38 & 43 & 15n46344 \\
		\shortstack[l]{$(1,1,-2,1,$\\\hspace*{1em}$1,-1,-1,-1)$} &  &  &  &  &  & 6 & 36 & 42 & 15n136662 \\
		\shortstack[l]{$(-1,-1,2,-1,$\\\hspace*{1em}$-1,1,1,1)$} &  &  &  &  &  & 6 & 34 & 40 & 15n7163 \\
		$(-1,-1,2,-1,-3,-1)$ &  &  &  &  &  & 6 & 30 & 36 & 15n44901 \\
		$(-2,-1,2,-1,-2,-1)$ &  &  &  &  &  & 4 & 27 & 31 & 15n7271 \\
		\shortstack[l]{$(1,1,-1,1,1,$\\\hspace*{1em}$1,-1,1,1,1)$} &  &  &  &  & 1 & 2 & 27 & 30 & 14n05084 \\
		\shortstack[l]{$(1,1,-1,1,$\\\hspace*{1em}$1,1,-3,-1)$} &  &  &  &  &  & 2 & 21 & 23 & 15n18437 \\
		\shortstack[l]{$(-1,-1,2,1,$\\\hspace*{1em}$-2,1,2,1)$} &  &  &  &  &  &  & 20 & 20 & 16n109254 \\
		$(-3,1,3,1)$ &  &  &  &  &  &  & 16 & 16 & 16n75751 \\
		$(1,1,-2,1,2,1)$ &  &  &  &  &  &  & 15 & 15 & 16n83432 \\
		$(3,-1,-3,-1)$ &  &  &  &  &  & 1 & 13 & 14 & 15n50006 \\
		\shortstack[l]{$(-1,-1,1,-1,-1,$\\\hspace*{1em}$-1,1,-1,-1,-1)$} &  &  &  &  &  & 2 & 12 & 14 & 15n4045 \\
		$(-4,-1)$ &  &  &  &  &  & 1 & 10 & 11 & 15n100300 \\
		$(-1,-1,2,-1,-2,-1)$ &  &  &  &  &  &  & 11 & 11 & 16n95423 \\
		\shortstack[l]{$(-1,-1,1,-1,$\\\hspace*{1em}$-1,-1,3,1)$} &  &  &  &  &  &  & 11 & 11 & 16n112237 \\
		$(2,1,-2,1,2,1)$ &  &  &  &  &  &  & 10 & 10 & 16n2504 \\
		$(2,-1,-2,-1)$ &  &  &  &  &  &  & 9 & 9 & 16n230504 \\
		\shortstack[l]{$(1,1,-1,1,1,$\\\hspace*{1em}$1,-2,-1,3,1)$} &  &  &  &  &  &  & 9 & 9 & 16n559389 \\
		\shortstack[l]{$(-1,-1,2,-1,-2,1,$\\\hspace*{1em}$1,1,-1,1,1,1)$} &  &  &  &  &  & 1 & 8 & 9 & 15n116130 \\
		$(-2,1,2,1)$ &  &  &  &  &  &  & 8 & 8 & 16n109919 \\
		$(4,1)$ &  &  &  &  &  & 1 & 6 & 7 & 15n5319 \\
		\shortstack[l]{$(1,1,-2,1,2,-1,$\\\hspace*{1em}$-1,-1,1,-1,-1,-1)$} &  &  &  &  &  &  & 7 & 7 & 16n237893 \\
		$(-1,-1,4,1)$ &  &  &  &  &  &  & 5 & 5 & 16n2016 \\
		\shortstack[l]{$(1,-1,-1,-1,2,$\\\hspace*{1em}$-1,-1,1,1,1)$} &  &  &  &  &  &  & 5 & 5 & 16n58929 \\
		\shortstack[l]{$(3,-1,-2,1,1,$\\\hspace*{1em}$1,-1,1,1,1)$} &  &  &  &  &  &  & 5 & 5 & 16n613346 \\
		\shortstack[l]{$(-2,1,1,1,$\\\hspace*{1em}$-1,-1,4,1)$} &  &  &  &  &  &  & 4 & 4 & 16n152762 \\
		\shortstack[l]{$(-3,1,2,-1,-1,$\\\hspace*{1em}$-1,1,-1,-1,-1)$} &  &  &  &  &  &  & 4 & 4 & 16n9474 \\
		\shortstack[l]{$(-1,1,1,1,-2,$\\\hspace*{1em}$-1,3,-1,-3,-1)$} &  &  &  &  &  &  & 4 & 4 & 16n200138 \\
		\shortstack[l]{$(1,-1,-1,-1,1,$\\\hspace*{1em}$-1,-1,-1,3,1)$} &  &  &  &  &  &  & 4 & 4 & 16n358983 \\
		\shortstack[l]{$(-1,1,1,1,-1,1,$\\\hspace*{1em}$1,1,-2,-1,3,1)$} &  &  &  &  &  &  & 4 & 4 & 16n151218 \\
		\shortstack[l]{$(1,-1,-1,-1,1,-1,$\\\hspace*{1em}$-1,-1,1,-1,-1,-1)$} &  &  &  &  &  & 2 & 2 & 4 & 15n27975 \\
		\shortstack[l]{$(1,-1,-1,-1,1,-1,$\\\hspace*{1em}$-1,-1,2,1,-3,-1)$} &  &  &  &  &  &  & 4 & 4 & 16n744700 \\
		$(-1,-1,3,-1,-3,-1)$ &  &  &  &  &  &  & 3 & 3 & 16n377014 \\
		\shortstack[l]{$(1,1,-1,1,$\\\hspace*{1em}$1,1,-2,-1)$} &  &  &  &  &  &  & 3 & 3 & 16n260752 \\
		\shortstack[l]{$(-2,1,1,1,-1,$\\\hspace*{1em}$-1,1,-1,-1,-1)$} &  &  &  &  &  &  & 3 & 3 & 16n292233 \\
		\shortstack[l]{$(-1,1,1,1,-2,$\\\hspace*{1em}$1,1,-1,-1,-1)$} &  &  &  &  &  &  & 3 & 3 & 16n241379 \\
		\shortstack[l]{$(-1,1,1,1,-1,$\\\hspace*{1em}$1,1,1,-3,-1)$} &  &  &  &  &  &  & 3 & 3 & 16n762143 \\
		\shortstack[l]{$(-1,1,1,1,-1,1,$\\\hspace*{1em}$1,1,-1,1,1,1)$} &  &  &  &  &  & 1 & 2 & 3 & 15n48790 \\
		\shortstack[l]{$(1,-1,-1,-1,2,1,-2,$\\\hspace*{1em}$1,1,1,-1,1,1,1)$} &  &  &  &  &  &  & 3 & 3 & 16n323991 \\
		$(1,1,1,-1,-1,-1)$ &  &  &  &  &  &  & 2 & 2 & 16n855309 \\
		\shortstack[l]{$(-1,-1,1,-1,$\\\hspace*{1em}$-1,-1,2,1)$} &  &  &  &  &  &  & 2 & 2 & 16n665135 \\
		\shortstack[l]{$(1,-1,-1,-1,1,-1,-1,$\\\hspace*{1em}$-1,1,-1,-1,-1,2,1)$} &  &  &  &  &  &  & 2 & 2 & 16n375186 \\
		$(-5,-1)$ &  &  &  &  &  &  & 1 & 1 & 16n108982 \\
		$(5,1)$ &  &  &  &  &  &  & 1 & 1 & 16n956730 \\
		$(3,-1,-3,-1,4,1)$ &  &  &  &  &  &  & 1 & 1 & 16n323807 \\
		$(3,1,-3,1,3,1)$ &  &  &  &  &  &  & 1 & 1 & 16n323429 \\
		\shortstack[l]{$(1,-1,-1,-1,2,$\\\hspace*{1em}$1,-3,1,3,1)$} &  &  &  &  &  &  & 1 & 1 & 16n970156 \\
		\shortstack[l]{$(1,1,1,-1,-1,-1,$\\\hspace*{1em}$2,1,-1,1,1,1)$} &  &  &  &  &  &  & 1 & 1 & 16n927087 \\
		\shortstack[l]{$(-1,1,1,1,-2,-1,2,$\\\hspace*{1em}$-1,-1,-1,1,-1,-1,-1)$} &  &  &  &  &  &  & 1 & 1 & 16n304720 \\
		\shortstack[l]{$(-1,1,1,1,-1,1,1,$\\\hspace*{1em}$1,-1,1,1,1,-2,-1)$} &  &  &  &  &  &  & 1 & 1 & 16n23706 \\
		\shortstack[l]{$(1,1,-1,1,1,1,-1,$\\\hspace*{1em}$1,1,1,-1,1,1,1)$} &  &  &  &  &  &  & 1 & 1 & 16n164055 \\
\bottomrule
\end{longtable}
\endgroup

\bibliographystyle{amsalphaurl}
\def\MR#1{\href{http://www.ams.org/mathscinet-getitem?mr=#1}{MR#1}}
\bibliography{biblio}

\providecommand{\bysame}{\leavevmode\hbox to3em{\hrulefill}\thinspace}
\providecommand{\MR}{\relax\ifhmode\unskip\space\fi MR }
\providecommand{\MRhref}[2]{%
  \href{http://www.ams.org/mathscinet-getitem?mr=#1}{#2}
}
\providecommand{\href}[2]{#2}
\begin{thebibliography}{DHST23}

\bibitem[BD04]{BurbanDrozd04:nodal}
Igor Burban and Yuriy Drozd, \emph{Derived categories of nodal algebras}, J. Algebra \textbf{272} (2004), no.~1, 46--94, \href {https://doi.org/10.1016/j.jalgebra.2003.07.025} {\path{doi:10.1016/j.jalgebra.2003.07.025}}. \MR{2029026}

\bibitem[Bon88]{Bondarenko88:bundles}
V.~M. Bondarenko, \emph{Bundles of semichained sets and their representations}, Akad. Nauk Ukrain. SSR Inst. Mat. Preprint (1988), no.~60, 32, English translation by W.~Crawley-Boevey, U.~Hansper, I.~Voulis (2020). \MR{1034077}

\bibitem[Bon91]{Bondarenko92:reps}
\bysame, \emph{Representations of bundles of semichained sets and their applications}, Algebra i Analiz \textbf{3} (1991), no.~5, 38--61. \MR{1186235}

\bibitem[DG25]{DunfieldGong}
Nathan~M. Dunfield and Sherry Gong, \emph{Ribbon concordances and slice obstructions: experiments and examples}, 2025, \href {https://arxiv.org/abs/arXiv:2512.21825} {\path{arXiv:arXiv:2512.21825}}.

\bibitem[DHST21]{DHST21:more-conc}
Irving Dai, Jennifer Hom, Matthew Stoffregen, and Linh Truong, \emph{More concordance homomorphisms from knot {F}loer homology}, Geom. Topol. \textbf{25} (2021), no.~1, 275--338, \href {https://doi.org/10.2140/gt.2021.25.275} {\path{doi:10.2140/gt.2021.25.275}}. \MR{4226231}

\bibitem[DHST23]{DHST:homol-cob}
\bysame, \emph{An infinite-rank summand of the homology cobordism group}, Duke Math. J. \textbf{172} (2023), no.~12, 2365--2432, \href {https://doi.org/10.1215/00127094-2022-0082} {\path{doi:10.1215/00127094-2022-0082}}. \MR{4654053}

\bibitem[DLS26]{slocal}
Nathan~M. Dunfield, Robert Lipshitz, and Dirk Sch\"utz, \emph{Local equivalence and refinements of {R}asmussen's s-invariant}, J. Topol. \textbf{19} (2026), no.~1, Paper No. e70056, 47, URL: \url{https://doi-org.uoregon.idm.oclc.org/10.1112/topo.70056}, \href {https://doi.org/10.1112/topo.70056} {\path{doi:10.1112/topo.70056}}. \MR{5018582}

\bibitem[FLL25]{FLL25:squeezed}
Peter Feller, Lukas Lewark, and Andrew Lobb, \emph{Squeezed knots}, Quantum Topol. \textbf{16} (2025), no.~4, 831--865, URL: \url{https://doi-org.uoregon.idm.oclc.org/10.4171/qt/187}, \href {https://doi.org/10.4171/qt/187} {\path{doi:10.4171/qt/187}}. \MR{4958851}

\bibitem[HHL21]{HHL21:conn-invol}
Kristen Hendricks, Jennifer Hom, and Tye Lidman, \emph{Applications of involutive {H}eegaard {F}loer homology}, J. Inst. Math. Jussieu \textbf{20} (2021), no.~1, 187--224, URL: \url{https://doi-org.uoregon.idm.oclc.org/10.1017/S147474801900015X}, \href {https://doi.org/10.1017/S147474801900015X} {\path{doi:10.1017/S147474801900015X}}. \MR{4205781}

\bibitem[Kap58]{Kaplansky58:theorem}
Irving Kaplansky, \emph{Projective modules}, Ann. of Math. (2) \textbf{68} (1958), 372--377, \href {https://doi.org/10.2307/1970252} {\path{doi:10.2307/1970252}}. \MR{100017}

\bibitem[Lew24]{Lewark24:local}
Lukas Lewark, \emph{A new concordance homomorphism from {K}hovanov homology}, 2024, \href {https://arxiv.org/abs/2401.08480} {\path{arXiv:2401.08480}}.

\bibitem[LS14]{LS14:refine-s}
Robert Lipshitz and Sucharit Sarkar, \emph{A refinement of {R}asmussen's {$s$}-invariant}, Duke Math. J. \textbf{163} (2014), no.~5, 923--952, \href {https://doi.org/10.1215/00127094-2644466} {\path{doi:10.1215/00127094-2644466}}. \MR{3189434}

\bibitem[LS22]{LS22:mixed}
\bysame, \emph{A mixed invariant of nonorientable surfaces in equivariant {K}hovanov homology}, Trans. Amer. Math. Soc. \textbf{375} (2022), no.~12, 8807--8849, \href {https://doi.org/10.1090/tran/8736} {\path{doi:10.1090/tran/8736}}. \MR{4504654}

\bibitem[ORS13]{ORSz13:odd}
Peter~S. Ozsv\'{a}th, Jacob Rasmussen, and Zolt\'{a}n Szab\'{o}, \emph{Odd {K}hovanov homology}, Algebr. Geom. Topol. \textbf{13} (2013), no.~3, 1465--1488, \href {https://doi.org/10.2140/agt.2013.13.1465} {\path{doi:10.2140/agt.2013.13.1465}}. \MR{3071132}

\bibitem[O'T13]{QuoteInvestigator}
Garson O'Toole, \emph{Quote origin: It's difficult to make predictions, especially about the future}, Quote Investigator\textsuperscript{\textregistered}, October 2013, Updated May 3, 2026. Accessed October 16, 2026, URL: \url{https://quoteinvestigator.com/2013/10/20/no-predict/}.

\bibitem[Pet13]{Petkova13:cables}
Ina Petkova, \emph{Cables of thin knots and bordered {H}eegaard {F}loer homology}, Quantum Topol. \textbf{4} (2013), no.~4, 377--409, \href {https://doi.org/10.4171/QT/43} {\path{doi:10.4171/QT/43}}. \MR{3134023}

\bibitem[Pop23]{Popovic:snakes}
David Popovi{\'c}, \emph{Link {F}loer homology splits into snake complexes and local systems}, 2023, \href {https://arxiv.org/abs/arXiv:2306.05546} {\path{arXiv:arXiv:2306.05546}}.

\bibitem[Ras10]{Rasmussen10:s}
Jacob Rasmussen, \emph{Khovanov homology and the slice genus}, Invent. Math. \textbf{182} (2010), no.~2, 419--447, \href {https://doi.org/10.1007/s00222-010-0275-6} {\path{doi:10.1007/s00222-010-0275-6}}. \MR{2729272}

\bibitem[Sar26]{Sarkar2026Personal}
Sucharit Sarkar, Personal communication, Aug. 2026, by teleconference.

\bibitem[Sch25]{Schutz:integral-s}
Dirk Sch\"utz, \emph{On an integral version of the {R}asmussen invariant}, Michigan Math. J. \textbf{75} (2025), no.~1, 65--88, \href {https://doi.org/10.1307/mmj/20226211} {\path{doi:10.1307/mmj/20226211}}. \MR{4873797}

\bibitem[Sch26]{KnotJob}
\bysame, \emph{Knot{J}ob (computer software)}, 2026, URL: \url{https://www.maths.dur.ac.uk/users/dirk.schuetz/knotjob.html}.

\bibitem[SSS20]{SSS20:odd-htpy}
Sucharit Sarkar, Christopher Scaduto, and Matthew Stoffregen, \emph{An odd {K}hovanov homotopy type}, Adv. Math. \textbf{367} (2020), 107112, 51, \href {https://doi.org/10.1016/j.aim.2020.107112} {\path{doi:10.1016/j.aim.2020.107112}}. \MR{4078823}

\bibitem[Sto20]{Stoffregen20:Seifert}
Matthew Stoffregen, \emph{Pin(2)-equivariant {S}eiberg-{W}itten {F}loer homology of {S}eifert fibrations}, Compos. Math. \textbf{156} (2020), no.~2, 199--250, \href {https://doi.org/10.1112/s0010437x19007620} {\path{doi:10.1112/s0010437x19007620}}. \MR{4044465}

\bibitem[Sto26]{Stoffregen2026Personal}
Matthew Stoffregen, Personal communication, Oct. 2026, in person.

\bibitem[Sza22]{SzaboCode}
Zolt\'an Szab\'o, \emph{Knot {F}loer homology calculator, version 3}, 2022, URL: \url{https://web.math.princeton.edu/~szabo/HFKcalc.html}.

\bibitem[Zem26]{Zemke2026Personal}
Ian Zemke, Personal communication, Sept. 2026, in person.

\end{thebibliography}
\end{document}